\documentclass[12pt,leqno]{article}
\usepackage{amssymb}
\usepackage[mathscr]{eucal}
\usepackage{amsmath,amssymb,latexsym,theorem,bbm}
\usepackage{color}

\newcommand{\SC}{\scriptstyle}

\newcommand{\CC}{\mathsf{C}}
\newcommand{\DD}{\mathsf{D}}
\newcommand{\CCC}{\mathbb{C}}

\newcommand{\NN}{\mathbb{N}}

\newcommand{\RR}{\mathbb{R}}

\newcommand{\ZZ}{\mathbb{Z}}

\newcommand{\bA}{{\boldsymbol{A}}}

\newcommand{\bb}{{\boldsymbol{b}}}
\newcommand{\bB}{{\boldsymbol{B}}}

\newcommand{\tC}{\widetilde{C}}

\newcommand{\bI}{{\boldsymbol{I}}}

\newcommand{\bM}{{\boldsymbol{M}}}
\newcommand{\bp}{{\boldsymbol{p}}}
\newcommand{\bP}{{\boldsymbol{P}}}

\newcommand{\bR}{{\boldsymbol{R}}}
\newcommand{\bcR}{{\boldsymbol{\cR}}}

\newcommand{\bS}{{\boldsymbol{S}}}
\newcommand{\bT}{{\boldsymbol{T}}}
\newcommand{\bu}{{\boldsymbol{u}}}

\newcommand{\bV}{{\boldsymbol{V}}}

\newcommand{\bx}{{\boldsymbol{x}}}
\newcommand{\bX}{{\boldsymbol{X}}}

\newcommand{\bU}{{\boldsymbol{U}}}

\newcommand{\bgamma}{{\boldsymbol{\gamma}}}
\newcommand{\bxi}{{\boldsymbol{\xi}}}

\newcommand{\bbeta}{{\boldsymbol{\beta}}}
\newcommand{\bvare}{{\boldsymbol{\vare}}}

\newcommand{\bzero}{{\boldsymbol{0}}}

\newcommand{\cA}{{\mathcal A}}
\newcommand{\cB}{{\mathcal B}}

\newcommand{\cD}{{\mathcal D}}

\newcommand{\cF}{{\mathcal F}}

\newcommand{\cM}{{\mathcal M}}
\newcommand{\bcM}{\boldsymbol{\cM}}

\newcommand{\cP}{{\mathcal P}}

\newcommand{\cR}{{\mathcal R}}

\newcommand{\cU}{{\mathcal U}}

\newcommand{\bcU}{\boldsymbol{\cU}}

\newcommand{\cX}{{\mathcal X}}
\newcommand{\bcX}{\boldsymbol{\cX}}

\newcommand{\cY}{{\mathcal Y}}

\newcommand{\cW}{{\mathcal W}}
\newcommand{\bcW}{\boldsymbol{\cW}}

\newcommand{\bcY}{\boldsymbol{\cY}}

\newcommand{\dd}{\mathrm{d}}
\newcommand{\ee}{\mathrm{e}}
\newcommand{\ii}{\mathrm{i}}
\newcommand{\slu}{{\SC\mathrm{lu}}}

\newcommand{\INARp}{\textup{INAR($p$)}}

\newcommand{\EE}{\operatorname{\mathbb{E}}}
\newcommand{\PP}{\operatorname{\mathbb{P}}}
\newcommand{\OO}{\operatorname{O}}
\newcommand{\tr}{\operatorname{tr}}
\newcommand{\var}{\operatorname{Var}}
\newcommand{\cov}{\operatorname{Cov}}

\newcommand{\vare}{\varepsilon}

\renewcommand{\mid}{\,|\,}
\newcommand{\bmid}{\,\big|\,}

\renewcommand{\leq}{\leqslant}
\renewcommand{\geq}{\geqslant}

\newcommand{\stoch}{\stackrel{\PP}{\longrightarrow}}
\newcommand{\distr}{\stackrel{\cD}{\longrightarrow}}
\newcommand{\distrf}{\stackrel{\cD_f}{\longrightarrow}}
\newcommand{\distre}{\stackrel{\cD}{=}}
\newcommand{\mean}{\stackrel{L_1}{\longrightarrow}}

\newcommand{\lu}{\stackrel{\slu}{\longrightarrow}}
\newcommand{\as}{\stackrel{{\mathrm{a.s.}}}{\longrightarrow}}
\newcommand{\ase}{\stackrel{{\mathrm{a.s.}}}{=}}

\newcommand{\bbone}{\mathbbm{1}}
\newcommand{\ns}{{\lfloor ns\rfloor}}
\newcommand{\nt}{{\lfloor nt\rfloor}}
\newcommand{\ntk}{{\lfloor n_1t_1\rfloor}}
\newcommand{\ntn}{{\lfloor n_2t_2\rfloor}}
\newcommand{\nT}{{\lfloor nT\rfloor}}
\newcommand{\proofend}{\hfill\mbox{$\Box$}}

\numberwithin{equation}{section}

\theoremstyle{change} \theorembodyfont{\em}
\newtheorem{Lem}{Lemma.}[section]
\newtheorem{Thm}[Lem]{Theorem.}

\newtheorem{Cor}[Lem]{Corollary.}

\theorembodyfont{\rm}
\newtheorem{Rem}[Lem]{Remark.}

\begin{document}

\begin{center}
 {\bfseries\Large Asymptotic behaviour of critical decomposable 2-type}
  \\[2mm]
 {\bfseries\Large Galton-Watson processes with immigration}

\vspace*{3mm}

 {\sc\large
  M\'aty\'as $\text{Barczy}^{*,\diamond}$,
  \ D\'aniel $\text{Bezd\'any}^{**}$,
  \ \framebox[1.1\width]{Gyula $\text{Pap}$}}

\end{center}

\vskip0.2cm

\noindent
 * ELKH-SZTE Analysis and Applications Research Group,
   Bolyai Institute, University of Szeged,
   Aradi v\'ertan\'uk tere 1, H--6720 Szeged, Hungary.

\noindent
 ** Former master student of Bolyai Institute, University of Szeged,
    Aradi v\'ertan\'uk tere 1, H-6720 Szeged, Hungary.

\noindent e-mails: barczy@math.u-szeged.hu (M. Barczy),
                    bezdanydaniel@gmail.com (D. Bezd\'any).

\noindent $\diamond$ Corresponding author.

\renewcommand{\thefootnote}{}
\footnote{\textit{2020 Mathematics Subject Classifications\/}: 60J80, 60F17. }
\footnote{\textit{Key words and phrases\/}:
  multi-type Galton-Watson process with immigration, critical, decomposable, asymptotic behaviour, squared Bessel process, stationary distribution.}
\footnote{M\'aty\'as Barczy is supported by the Ministry of Innovation and Technology of Hungary
 from the National Research, Development and Innovation Fund, project no.\ TKP2021-NVA-09.}

\vspace*{0.1cm}

\begin{abstract}
In this paper the asymptotic behaviour of a critical 2-type Galton-Watson process with immigration is described
 when its offspring mean matrix is reducible, in other words, when the process is decomposable.
It is proved that, under second or fourth order moment assumptions on the offspring and immigration distributions,
 a sequence of appropriately scaled random step processes formed from a critical decomposable
 2-type Galton-Watson process with immigration converges weakly.
The limit process can be described using one or two independent squared Bessel processes and possibly the unique stationary
 distribution of an appropriate single-type subcritical Galton-Watson process with immigration.
Our results complete and extend the results of Foster and Ney (1978) for some strongly critical decomposable  2-type Galton-Watson processes with immigration.
\end{abstract}

\section{Introduction}
\label{section_intro}

The study of the limit behaviour of critical branching processes has a long tradition and history.
For critical branching processes without immigration, so-called conditioned limit theorems
 and for critical branching processes with immigration, unconditioned limit theorems are usually established.
Investigation of asymptotic properties of critical multi-type Galton-Watson processes with or without immigration
 goes back at least to the 60's and it is still an active area of research.
Below we give a review of some results in this field, with our main focus on critical decomposable multi-type Galton-Watson processes
 without or with immigration.
Then after giving the precise definition and basic properties of multi-type Galton-Watson processes with immigration (see Section \ref{section_multi_branching}),
 we present our main results on the asymptotic behaviour of critical decomposable 2-type Galton-Watson processes with immigration (see Section \ref{Section_conv_results}).
It is proved that, under second or fourth order moment assumptions on the offspring and immigration distributions,
 a sequence of appropriately scaled random step processes formed from a critical decomposable 2-type Galton-Watson process with immigration converges weakly.
The limit process can be described using one or two independent squared Bessel processes and possibly the unique stationary
 distribution of an appropriate single-type subcritical Galton-Watson process with immigration.
Our results complete and extend the results of Foster and Ney \cite[Section 9]{FosNey2} for some strongly critical decomposable 2-type Galton-Watson processes with immigration.

A multi-type Galton-Watson process with immigration is referred to as subcritical,
 critical, or supercritical if the spectral radius of its offspring mean matrix
 is less than 1, equal to 1 or greater than 1, respectively.
A multi-type Galton-Watson process with immigration is called indecomposable and decomposable
 if its offspring mean matrix is irreducible and reducible, respectively.
An indecomposable multi-type Galton-Watson process with immigration is called primitive (also called positively regular)
 if its offspring mean matrix is primitive.
For more details on these concepts, see Section \ref{section_multi_branching}.

For a review of some of the results on the asymptotic behaviour of critical single-type Galton-Watson processes with immigration,
 see, e.g., the introduction of Barczy et al.\ \cite{BarBezPap}.
Here we only mention the result of Wei and Winnicki \cite[Theorem 2.1]{WW} who proved weak convergence of a sequence of random step processes
 \ $(n^{-1} X_{\nt})_{t\geq 0}$, \ $n\geq 1$, \ as \ $n\to\infty$, \ formed from a critical single-type Galton-Watson process with
 immigration \ $(X_k)_{k\geq 0}$ \ under second order moment assumptions and characterized the limit process as a squared Bessel process
 (for more details, see Theorem \ref{thm:critical}).

Next, we make an overview of the existing results on the asymptotic behaviour of critical, decomposable multi-type
 Galton-Watson processes without immigration or with immigration, then those of critical, indecomposable or primitive
 branching processes, and finally we recall a result on the asymptotic behaviour of supercritical decomposable multi-type Galton-Watson
 processes without immigration, and on the weak convergence of a sequence of appropriately scaled (arbitrary) 2-type Galton-Watson processes with immigration
 towards a continuous state and continuous time branching process with immigration, respectively.

Under second order moment assumptions on the offspring distributions, Foster and Ney \cite{FosNey}
 described the asymptotic behaviour of the extinction probability of a critical decomposable multi-type Galton-Watson process without immigration
 and with a deterministic initial distribution.

Foster and Ney \cite{FosNey2} proved conditioned limit theorems for some special strongly critical decomposable multi-type
 Galton-Watson processes without immigration
 (see \cite[Theorem 2]{FosNey2}),
 and unconditioned limit theorems for some special strongly critical decomposable multi-type Galton-Watson processes with immigration
 (see \cite[Theorems 4 and 5]{FosNey2}).
They also specialized their results in case of 2 types, see \cite[Section 6 and page 42]{FosNey2}.
This special 2-type case in case of immigration corresponds to our case (2) (see \eqref{tablazat_esetek}).
For a special strongly critical decomposable 2-type Galton-Watson process with immigration $(\bX_k)_{k\geq 0} = ((X_{k,1},X_{k,2}))_{k\geq 0}$,
 under second order moment assumptions on the offspring distributions and first order moment assumptions on the immigration distributions,
 Foster and Ney \cite[Theorems 4 and 5]{FosNey2} showed that $(n^{-1} X_{n,1}, n^{-2} X_{n,2})$ converges in distribution as $n\to\infty$,
 and they characterized the limit distribution by its Laplace transform containing an integral as well (see \cite[formula (9.11)]{FosNey2}).
Note that here the right normalization is \ $n^{-i}$ \ for the \ $i^{\mathrm{th}}$-type, \ $i=1,2$.
In our Theorem \ref{main_2} we extend this result of Foster and Ney by proving weak convergence of a sequence of appropriately scaled random step
 processes formed from \ $(\bX_k)_{k\geq 0}$, \ and we characterize the limit process as a solution of a system of stochastic differential
 equations.
In Section \ref{Rem2}, we compare our results in Theorem \ref{main_2} and those of Foster and Ney \cite[Section 9]{FosNey2}
 in detail, where we give a closed formula of the Laplace transform of the above mentioned limit distribution in Foster and Ney
 \cite[formula (9.11)]{FosNey2} as well.
Foster and Ney \cite[Section 8]{FosNey2} also indicated some conjectures about the nature of possible results for critical decomposable
 multi-type Galton-Watson processes without immigration not supposing some of their restrictive hypotheses including strong criticality.
Our results in Section \ref{Section_conv_results} for critical decomposable 2-type Galton-Watson processes with immigration
 handle all the remaining cases not included in Foster and Ney \cite[Section 9]{FosNey2}, and instead of convergence of one-dimensional
 distributions we can prove weak convergence of a sequence of appropriate random step processes formed from the branching processes in question.

Sugitani \cite{Sug, Sug2} extended the results of Foster and Ney \cite{FosNey2} on conditional limit theorems for
 critical decomposable multi-type Galton-Watson processes without immigration (\cite[Theorems 2.2 and 2.4]{Sug2}),
 and some unconditional limit theorems were established as well (\cite[Theorems 2.1, 2.3 and 2.5]{Sug2}).

Zubkov \cite{Zub} proved some conditioned limit theorems for critical decomposable 2-type Galton-Watson processes without immigration
 such that the generating functions of the offspring distributions satisfy some regularity assumptions yielding that
 the offspring distributions do not have finite second moments.

Studying asymptotic properties of critical decomposable multi-type Galton-Watson processes with or without immigration
 is still attracting the attention of researchers.
Smadi and Vatutin \cite{SmaVat} considered a critical decomposable 2-type Galton-Watson process \ $(\bX_k)_{k\geq 0}$ \ without immigration
 such that the variance of the offspring distributions may be infinite.
Let \ $X_{m,n,i}$ \ be the number of individuals of type \ $i$ \ alive at time \ $m$ \ and having  descendants
 at time \ $n$, \ where \ $i=1,2$ \ and \ $m<n$, \ in other words,
 \ $X_{m,n,i}$ \ is the number of the type \ $i$ \ ancestors alive at generation \ $m$ \ of all the individuals
 of the population at generation \ $n$.
\ Smadi and Vatutin \cite{SmaVat} described the asymptotic behaviour of the conditional distribution
 of \ $(X_{m,n,1},X_{m,n,2})$ \ given that \ $\bX_n\ne \bzero$ \ as \ $m,n\to\infty$.
\ Here for each \ $n\in\NN$, \ $(X_{m,n,1},X_{m,n,2})_{m\in\{0,1,\ldots, n-1\}}$ \ can be thought as the family tree relating the individuals alive
 at time \ $n$, \ and \ $(X_{m,n,1},X_{m,n,2})$, \ $0\leq m \leq n$, \ $m,n\in\ZZ_+$, \ is sometimes called a  reduced branching process.
For strongly critical decomposable multi-type Galton-Watson processes without immigration, similar problems were considered and solved by
 Vatutin \cite{Vat1, Vat2}.

Next, we recall some results on the asymptotic behaviour of critical,
 indecomposable or primitive multi-type Galton-Watson processes without immigration or with immigration.

For a certain class of critical, primitive (also called positively regular) multi-type branching processes \ $(Z_n)_{n\geq 0}$ \ without immigration
 Mullikin \cite[Theorems 8 and 9]{Mul} characterized the limit of the conditional expectation and distribution of \ $n^{-1}Z_n$ \ given that
 \ $Z_n \ne 0$ \ as \ $n\to\infty$.
\ Mullikin's results \cite[Theorems 8 and 9]{Mul} are in fact much more general, a discrete time temporally homogeneous Markov process \ $(Z_n)_{n\geq 0}$ \
 was considered, where the range of \ $Z_n$ \ is a set of finitely additive, non-negative and integer-valued set functions on an abstract set
 (representing the set of possibly infinite number of types) furnished with a $\sigma$-algebra, and \ $Z_0$ \ is a given non-random functional.

Joffe and M\'etivier \cite[Theorem 4.3.1]{JofMet} studied a sequence \ $(\bX^{(n)}_k)_{k\geq 0}$, \ $n\geq 1$, \
 of critical multi-type Galton-Watson processes with the same offspring distributions having finite second moments,
 but without immigration and starting from a deterministic initial value
 \ $\bX^{(n)}_0$, \ supposing that the offspring mean matrix is primitive and \ $n^{-1}\bX^{(n)}_0$ \ converges to
 a non-zero (deterministic) limit as \ $n\to\infty$.
\ They determined the limiting behaviour of the martingale part
 \ $\big(n^{-1}\sum_{k=1}^{\nt} (\bX_k^{(n)} - \EE(\bX_k^{(n)} \mid \bX_0^{(n)},\ldots,\bX_{k-1}^{(n)}))\big)_{t\geq 0}$ \ as \ $n\to\infty$.
\ Joffe and M\'etivier \cite[Theorem 4.2.2]{JofMet} also studied a sequence of multi-type Galton-Watson processes without immigration
 \ $(\bX^{(n)}_k)_{k\geq 0}$, \ $n\geq 1$, \
 which is nearly critical of special type (see (i) of Theorem 4.2.2 in \cite{JofMet}), and, under second order moment assumptions and a Lindeberg-type condition,
 they proved that the sequence \ $(n^{-1}\bX^{(n)}_{\nt})_{t\geq 0}$ \ converges in distribution towards a diffusion process as \ $n\to\infty$.

Isp\'any and Pap \cite[Theorem 3.1]{IspPap2} described the asymptotic behaviour of a sequence of critical primitive (also called positively regular)
 multi-type Galton-Watson processes with immigration \ $(\bX^{(n)}_k)_{k\geq 0}$ \ sharing the same offspring and immigration distributions,
 but having possibly different initial distributions such that \ $n^{-1} \bX^{(n)}_0$ \ converges in distribution to \ $\cX\bu$ \
 as \ $n\to\infty$, \ where \ $\cX$ \ is a nonnegative random variable with distribution \ $\mu$ \ and
 \ $\bu$ \ is the Perron (right) eigenvector of the offspring mean matrix.
Under fourth order moment assumptions on the offspring and immigration distributions, they showed that
 \ $(n^{-1}\bX^{(n)}_{\nt})_{t\geq 0}$ \ converges in distribution as \ $n\to\infty$.
\ They characterized the limit process as \ $(\cX_t \bu)_{t\geq 0}$, \ where \ $(\cX_t)_{t\geq 0}$ \ is a squared Bessel process
 with initial distribution \ $\mu$.
\ Here it is interesting to point out the fact the limiting diffusion process \ $(\cX_t \bu)_{t\geq 0}$ \ is always one-dimensional
 in the sense that for all \ $t\geq 0$, \ the distribution of \ $\cX_t \bu$ \ is concentrated on the ray \ $[0,\infty)\cdot \bu$,
 \ while the original sequence of branching processes does not have this property.

For a critical indecomposable \ $p$-type Galton-Watson process \ $(\bX_k)_{k\geq 0}$ \ with immigration and starting from \ $\bX_0=\bzero$,
 \ Danka and Pap \cite{DanPap} obtained a generalization of Theorem 3.1 in Isp\'any and Pap \cite{IspPap2}.
In the indecomposable case the set of types \ $\{1,\ldots,p\}$ \ can be partitioned according to communication of types, namely,
 into \ $r$ \ nonempty mutually disjoint subsets \ $D_1,\ldots,D_r$ \ such that an individual of type \ $j$ \
 may not have offsprings of type \ $i$ \ unless there exists \ $\ell \in \{1,\ldots,r\}$ \ with \ $i\in D_{\ell-1}$ \
 and \ $j\in D_\ell$, \ where the subscripts are considered modulo \ $r$
 \ (for more details, see, e.g., Danka and Pap \cite[Section 2]{DanPap}).
This partitioning is unique up to cyclic permutation of the subsets, and the number \ $r$ \ is called the index of cyclicity
 (in other words, the index of imprimivity) of the mean matrix \ $\bA$.
\ Note that \ $r=1$ \ if and only if the matrix \ $\bA$ \ is primitive, i.e., the branching process in question is primitive (in other words positively regular).
Under second order moment assumptions on the offspring and immigration distributions
 for the given \ $p$-type Galton-Watson process \ $(\bX_k)_{k\geq 0}$ \ with immigration,
 using Theorem \ref{Conv2Funct}, Danka and Pap \cite[Theorem 3.1]{DanPap} determined the joint asymptotic behaviour of the random step processes
 \ $\big( (nr)^{-1} \bX_{r\nt + i-1}\big)_{t\geq 0}$, \ $n\in\NN$, \ $i\in\{1,\ldots,r\}$ \ towards
 the limiting diffusion processes \ $(\bA^{r-i+1} \bcY_t)_{t\geq 0}$, \ $i\in\{1,\ldots,r\}$ \ as \ $n\to\infty$.
\ Here the process \ $(\bcY_t)_{t\geq 0}$ \ is 1-dimensional in the sense that for each \ $t\geq 0$, \
 the distribution of \ $\bcY_t$ \ is concentrated on the ray \ $[0,\infty)\cdot \bu$, \ where \ $\bu$ \ is the Perron (right) eigenvector of
 the offspring mean matrix \ $\bA$.

To close the review of existing and connecting literature, we recall two more results that are somewhat connected.
It is interesting to note that Kesten and Stigum
 \cite[Theorems 2.1, 2.2 and 2.3]{KesSti} considered a supercritical decomposable multi-type Galton-Watson process
 $(\bX_n)_{n\geq 0}$ without immigration and with a fixed deterministic initial distribution,
 and they proved that appropriately normalizing $\bX_n$ (or its appropriate subsequence) it converges almost surely to a random limit vector as $n\to\infty$.
The normalizing factors in question always have the form $n^{-\gamma} \lambda^{-n}$, where $\gamma$ is a non-negative integer
 and $\lambda$ is a positive real number greater than or equal to one.
In some cases, they specialized their results to 2-type processes as well, see Kesten and Stigum \cite[page 321]{KesSti}.

Ma \cite[Theorem 2.1, (i)]{Ma} established sufficient conditions for the weak convergence of a sequence of
 (arbitrary, not necessarily critical or decomposable) 2-type Galton-Watson processes with immigration
 towards a continuous state and continuous time branching process with immigration
 using appropriate time and space scalings such that the time scaling in question depends on the immigration distributions
  (and being different from what we will consider in our limit theorems in Section \ref{Section_conv_results}).
Ma \cite{Ma} proved the convergence of the sequence of infinitesimal generators of the branching processes in question towards
 the infinitesimal generator of the limit process.

The paper is organized as follows.
In Section \ref{section_multi_branching} we recall the definition of multi-type Galton-Watson processes with immigration,
 their classification as subcritical, critical and supercritical ones, and the special classes of
 indecomposable, decomposable and primitive (also called positively regular) branching processes.
Section \ref{Section_conv_results} contains our main results on the asymptotic behaviour of critical decomposable 2-type Galton-Watson
 processes with immigration, see Theorems \ref{main_1}--\ref{main_5}.
The investigation of such processes can be reduced to five cases
 presented in \eqref{tablazat_esetek} according to the form of the offspring mean matrix.
We also explain how decomposable \ $2$-type Galton-Watson processes may model the sizes of a geographically structured population
 divided into two parts.
Under second or fourth order moment assumptions on the offspring and immigration distributions, in the above mentioned five cases,
 we describe the limit behaviour of a sequence of appropriately scaled random step processes
 formed from a critical decomposable 2-type Galton-Watson process with immigration.
The limit process can be described using either one or two independent squared Bessel processes,
 and possibly the unique stationary distribution of an appropriate single-type subcritical Galton-Watson process with immigration.
This is a new phenomenon compared to the existing results on critical indecomposable (specially primitive) multi-type Galton-Watson
 processes with immigration.
Concerning Theorem \ref{main_3}, we formulate a conjecture on the independence of the limit processes of the two coordinate processes.
We note that Theorem \ref{main_2} can be considered as a functional version of Theorems 4 and 5 in Foster and Ney \cite{FosNey2} for
 some strongly critical decomposable 2-type Galton-Watson processes with immigration.
For a detailed comparison of our results in Theorem \ref{main_2} and those of Foster and Ney \cite[Section 9]{FosNey2},
 see Section \ref{Rem2}, where, as a consequence of Theorem \ref{main_2}, we also give a
  functional generalization of the Corollary on page 42 in Foster and Ney \cite{FosNey2}, regarding
  the joint convergence of the appropriately normalized population size at time \ $n$ \ and total progeny up
 to time \ $n$ \ of a critical single-type Galton--Watson process with immigration as \ $n\to\infty$.
In Corollary \ref{Cor_relative_frequency} we describe the asymptotic behaviour of the relative frequency of individuals of type 2
 with respect to individuals of type 1 under the conditions of Theorem \ref{main_2} together with that the mean of immigration distribution of type 1
 individuals is positive and finite.
For different models, one can find similar results, e.g., in Jagers \cite[Corollary 1]{Jag} and in Yakovlev and Yanev \cite[Theorem 2]{YakYan}
 (for more details, see the paragraph before Corollary \ref{Cor_relative_frequency}).
Remark \ref{Rem1} is devoted to a discussion on the moment conditions in Theorems \ref{main_1}, \ref{main_4} and \ref{main_5},
  we explain why we suppose the finiteness of the fourth order moments of the offspring and immigration distributions in these theorems.
Section \ref{Section_prel_proof} contains some preliminaries for the proofs such as a formula for the powers of
 the offspring mean matrix and a useful decomposition of the process using martingale differences (see \eqref{Xdeco}).
Sections \ref{Proof1}--\ref{Proof5} are devoted to the proofs of Theorems \ref{main_1}--\ref{main_5} and Corollary \ref{Cor_relative_frequency}.
We close the paper with four appendices.
Appendix \ref{app:moments} contains some formulae and estimates for the first, second and fourth order moments of
 the coordinates of the branching process in question and those of the derived martingale differences, respectively.
These estimates are extensively used in the proofs.
In Appendix \ref{app:singletype_GWI} we present a result on the asymptotic behaviour of finite dimensional distributions
 of a single-type subcritical Galton-Watson process with immigration satisfying first order moment conditions, which may be known,
 but we could not address any reference for it, so we provided a proof as well.
We also recall the asymptotic behaviour of a single-type critical Galton-Watson process with immigration due to Wei and Winnicki \cite[Theorem 2.1]{WW}.
Appendix \ref{CMT} contains a version of the continuous mapping theorem.
In Appendix \ref{section_conv_step_processes} we recall a result about the convergence of random step processes towards
 a diffusion process due to Isp\'any and Pap \cite{IspPap}, this result is heavily used in our proofs.

\section{Multi-type Galton-Watson processes with immigration}\label{section_multi_branching}

Let \ $\ZZ_+$, \ $\NN$, \ $\RR$, \ $\RR_+$, \ $\RR_{++}$ \ and \ $\CCC$ \ denote the set
 of non-negative integers, positive integers, real numbers, non-negative real
 numbers, positive real numbers and complex numbers, respectively.
For \ $x,y\in\RR$, \ the minimum of \ $x$ \ and \ $y$ \ is denoted by \ $x\wedge y$.
\ The Euclidean norm on \ $\RR^d$ \ is denoted by \ $\Vert\cdot\Vert$, \ where \ $d\in\NN$.
\ The \ $d\times d$ \ identity matrix is denoted by \ $\bI_d$.
\ For a function \ $f:\RR\to\RR$, \ its positive part is denoted by \ $f^+$.
\ Every random variable will be defined on a fixed probability space
 \ $(\Omega, \cA, \PP)$.
\ Convergence in probability, convergence in \ $L_1$, \ convergence almost surely, equality in distribution and
 almost sure equality is denoted by \ $\stoch$, \ $\mean$, \ $\as$, \ $\distre$ \ and \ $\ase$, \ respectively.
We will use \ $\distrf$ \ for the weak convergence of the finite dimensional distributions,
 and \ $\distr$ \ for the weak convergence of \ $\RR^d$-valued stochastic processes with sample
 paths in \ $\DD(\RR_+,\RR^d)$, \ where $d\in\NN$ and \ $\DD(\RR_+,\RR^d)$ \ denotes the space of
 $\RR^d$-valued c\`adl\`ag functions defined on \ $\RR_+$ \ (for more details and notations, e.g., for \ $\lu$,
 \ see Appendix \ref{CMT}).
Given a non-empty set \ $I$, \ a stochastic process \ $(Y_t)_{t\in I}$ \ is called an i.i.d. process
 if the random variables \ $\{Y_t : t \in I\}$ \ are independent (i.e., for each \ $m\in\NN$ \ and each subset \ $\{t_1,\ldots,t_m\}\subset I$, \
 the random variables \ $Y_{t_1},\ldots,Y_{t_m}$ \ are independent) and identically distributed.
We note that in time series analysis, by a white noise (process), one usually means an uncorrelated, zero mean process having finite second moments,
 so according to our definition, an i.i.d.\ process is not necessarily a white noise (process).

We will investigate a certain 2-type Galton-Watson process with immigration.
First we recall the definition and first order moment formulae of
\ $p$-type Galton-Watson processes with immigration, where \ $p \in \NN$.

For each \ $k \in \ZZ_+$ \ and \ $i \in \{ 1, \dots, p \}$, \ the number of
 individuals of type \ $i$ \ in the \ $k^\mathrm{th}$ \ generation is denoted by
 \ $X_{k,i}$.
\ For simplicity, we suppose that the initial values are
 \ $X_{0,i} = 0$, \ $i \in \{ 1, \dots, p \}$.
\ By \ $\xi_{k,j,i,\ell}$ \ we denote the number of type \ $\ell$ \ offsprings
 produced by the \ $j^\mathrm{th}$ \ individual who is of type \ $i$
 \ belonging to the \ $(k-1)^\mathrm{th}$ \ generation.
The number of type \ $i$ \ immigrants in the \ $k^\mathrm{th}$ \ generation is denoted by \ $\vare_{k,i}$.
\ Consider the random vectors
 \[
   \bX_k := \begin{bmatrix}
             X_{k,1} \\
             \vdots \\
             X_{k,p}
            \end{bmatrix} , \qquad
   \bxi_{k,j,i} := \begin{bmatrix}
                    \xi_{k,j,i,1} \\
                    \vdots \\
                    \xi_{k,j,i,p}
                   \end{bmatrix} , \qquad
   \bvare_k := \begin{bmatrix}
                \vare_{k,1} \\
                \vdots \\
                \vare_{k,p}
               \end{bmatrix} .
 \]
Then we have
 \begin{equation}\label{MBPI(d)}
  \bX_k = \sum_{i=1}^p \sum_{j=1}^{X_{k-1,i}} \bxi_{k,j,i} + \bvare_k , \qquad
  k \in \NN,
 \end{equation}
 with \ $\bX_0=\bzero$ \ (and using the convention \ $\sum_{j=1}^0:=\bzero$).
\ Here
 \ $\big\{\bxi_{k,j,i}, \, \bvare_k
          : k, j \in \NN, \, i \in \{ 1, \dots, p \} \big\}$
 \ are supposed to be independent.
Moreover, \ $\big\{\bxi_{k,j,i} : k, j \in \NN\big\}$ \ for each
 \ $i \in \{1, \dots, p\}$, \ and \ $\{\bvare_k : k \in \NN\}$ \ are
 supposed to consist of identically distributed \ $\ZZ_+^p$-valued random vectors.
For notational convenience, let \ $\{\bxi_i : i \in \{1, \ldots, p\}\}$ \ and \ $\bvare$ \ be random vectors
 such that \ $\bxi_i \distre \bxi_{1,1,i}$ \ for all \ $i \in \{1, \ldots, p\}$ \ and \ $\bvare \distre \bvare_1$.

In all what follows we will suppose
 \begin{align}\label{help4}
   \EE(\|\bxi_i\|^2) < \infty, \quad  i=1, \dots, p, \qquad \text{and} \qquad  \EE (\|\bvare\|^2) < \infty.
 \end{align}
Introduce the notations
 \begin{gather*}
  \bA := \begin{bmatrix}
                \EE(\bxi_1) & \cdots & \EE(\bxi_p)
               \end{bmatrix} \in \RR^{p \times p}_+ , \qquad
  \bb := \EE(\bvare) \in \RR^p_+ , \\
  \bV^{(i)} := \var(\bxi_i) \in \RR^{p \times p} , \qquad i\in\{1,\ldots,p\}, \qquad
  \bV^{(0)} := \var(\bvare) \in \RR^{p \times p} .
 \end{gather*}
The matrix \ $\bA$ \ and the vector \ $\bb$ \ is called the offspring mean matrix and immigration mean vector
 of \ $(\bX_k)_{k\in\ZZ_+}$, \ respectively.
Note that some authors define the offspring mean matrix as \ $\bA^\top$.
\ For \ $k \in \ZZ_+$, \ let \ $\cF_k^\bX := \sigma( \bX_0,\bX_1 , \dots, \bX_k)$, \ where \ $\cF_0^\bX=\{\emptyset, \Omega\}$ \ (due to \ $\bX_0=\bzero$).
\ By \eqref{MBPI(d)}, we get
 \begin{equation}\label{mart}
  \EE(\bX_k \mid \cF_{k-1}^\bX)
  = \sum_{i=1}^p X_{k-1,i}  \EE(\bxi_i)  + \bb
  = \bA \bX_{k-1} + \bb , \qquad k\in\NN.
 \end{equation}
Consequently,
 \begin{equation*}
  \EE(\bX_k) = \bA \EE(\bX_{k-1}) + \bb , \qquad k \in \NN ,
 \end{equation*}
 and, since \ $\bX_0=\bzero$, \ we have
 \begin{equation*}
  \EE(\bX_k) = \sum_{j=0}^{k-1} \bA^j \bb , \qquad k \in \NN .
 \end{equation*}
Hence the offspring mean matrix \ $\bA$ \ plays a crucial role in the
 asymptotic behaviour of the sequence \ $(\EE(\bX_k))_{k\in\ZZ_+}$.
\ A \ $p$-type Galton-Watson process \ $(\bX_k)_{k\in\ZZ_+}$ \ with immigration
 is referred to respectively as \emph{subcritical}, \emph{critical} or
 \emph{supercritical} if \ $\varrho(\bA) < 1$, \ $\varrho(\bA) = 1$
 \ or \ $\varrho(\bA) > 1$, \ where \ $\varrho(\bA)$ \ denotes the
 spectral radius of the matrix \ $\bA$, \ i.e., the maximum of the absolute values of the eigenvalues
 of \ $\bA$ \ (see, e.g., Athreya and Ney \cite[V.3]{AN} or Quine \cite{Q}).

A multi-type Galton-Watson process \ $(\bX_k)_{k\in\ZZ_+}$ \ with immigration is called
 indecomposable and decomposable if its offspring mean matrix \ $\bA$ \ is irreducible and reducible, respectively.
We recall that the matrix \ $\bA$ \ is called reducible if there exist a permutation matrix
 \ $\bP \in \RR^{p \times p}$ \ and \ $q \in \{1, \ldots, p-1\}$ \ such that
 \[
   \bA
   = \bP \begin{bmatrix} \bR & \bzero \\ \bS & \bT \end{bmatrix} \bP^\top ,
 \]
 where \ $\bR \in \RR^{q \times q}$, \ $\bS \in \RR^{(p-q) \times q}$, \ $\bT \in \RR^{(p-q) \times (p-q)}$ \ and \ $\bzero \in \RR^{q \times (p-q)}$ \ is a null matrix.
The matrix \ $\bA$ \ is called irreducible if it is not reducible; see, e.g.,
 Horn and Johnson \cite[Definitions 6.2.21 and 6.2.22]{HJ}.
We do emphasize that no 1-by-1 matrix is reducible.
It is known that the matrix \ $\bA$ \ is irreducible if and only if for all \ $i,j\in\{1,\ldots,p\}$ \ there exists \ $n_{i,j}\in\NN$ \ such that
 the matrix entry \ $(\bA^{n_{i,j}})_{i,j}$ \ is positive.
An indecomposable multi-type Galton-Watson process \ $(\bX_k)_{k\in\ZZ_+}$ \ with immigration is called primitive (also called positively regular)
 if its offspring mean matrix \ $\bA$ \ is primitive, i.e., there exists an \ $n\in\NN$ \ such that
 the matrix entry \ $(\bA^n)_{i,j}$ \ is positive for each \ $i,j\in\{1,\ldots,p\}$.

\section{Convergence of random step processes}\label{Section_conv_results}

In what follows we consider a critical decomposable 2-type Galton-Watson process \ $(\bX_k)_{k\in\ZZ_+}$ \ with
 immigration starting from \ $\bX_0=\bzero$, \ and we suppose that the moment conditions \eqref{help4} hold.
Since \ $p=2$ \ and \ $\bA=(a_{i,j})_{i,j=1}^2\in\RR^{2\times 2}$ \ is reducible, we have \ $a_{1,2}=0$ \ or \ $a_{2,1}=0$.
\ Note also that if \ $(X_{k,1},X_{k,2})_{k\in\ZZ_+}$ \ is a decomposable 2-type Galton-Watson process with immigration having an offspring
 mean matrix with \ $(1,2)$-entry \ $0$, \ then \ $(X_{k,2},X_{k,1})_{k\in\ZZ_+}$ \ is also a decomposable 2-type Galton-Watson process with immigration
 having an offspring mean matrix with \ $(2,1)$-entry \ $0$.
\ Because of this, when dealing with decomposable $2$-type Galton-Watson processes with immigration
 it is enough to focus on those ones which have an offspring mean matrix with \ $(1,2)$-entry \ $0$.
\ So we may assume that the offspring mean matrix \ $\bA$ \ and the immigration mean vector \ $\bb$ \ take the following forms:
 \[
   \bA
       = \begin{bmatrix}
                 \EE(\bxi_1) & \EE(\bxi_2)
         \end{bmatrix}
   = \begin{bmatrix}
      a_{1,1} & 0 \\
      a_{2,1} & a_{2,2}
     \end{bmatrix}  \qquad \text{and} \qquad
   \bb = \begin{bmatrix} b_1 \\ b_2 \end{bmatrix},
 \]
 respectively, with \ $\varrho(\bA) = \max \{a_{1,1}, a_{2,2}\} = 1$.
\ Taking into account that \ $a_{1,2}=0$ \ implies \ $\PP(\xi_{1,1,2,1} =0)=1$,
 \ Equation \eqref{MBPI(d)} with \ $p=2$ \ takes the form
 \begin{align}\label{help1}
   \begin{bmatrix} X_{k,1} \\ X_{k,2} \end{bmatrix}
   = \sum_{j=1}^{X_{k-1,1}}
      \begin{bmatrix} \xi_{k,j,1,1} \\ \xi_{k,j,1,2} \end{bmatrix}
     + \sum_{j=1}^{X_{k-1,2}}
        \begin{bmatrix} 0 \\ \xi_{k,j,2,2} \end{bmatrix}
     + \begin{bmatrix} \vare_{k,1} \\ \vare_{k,2} \end{bmatrix} , \qquad
   k \in \NN ,
 \end{align}
 with \ $[X_{0,1}, X_{0,2}]=[0,0]$.
\ For a decomposable $2$-type Galton-Watson process \ $(\bX_k)_{k\in\ZZ_+}$ \ with immigration given by \eqref{help1},
 the individuals of type \ $1$ \ may produce individuals of types \ $1$ \ or $2$, \ and the individuals of
 type \ $2$ \ may produce individuals of type \ $2$ \ only.
This process may be viewed as a stochastic model of the sizes of a geographically structured population divided into two parts such that
 \begin{itemize}
    \item the individuals are located at one of the two parts, and the location of an individual is considered as its type,
    \item the newborn individuals of the part $1$ either stay at the part $1$ or migrate, just after their birth, to the part \ $2$,
    \item the newborn individuals of the part 2 stay at the part 2 (they do not migrate),
    \item at each step immigrants (newcomers) may join the part \ $i\in\{1,2\}$ \ and they become individuals of the part \ $i$,
    \item the offspring and immigration distributions depend on the parts on which the individuals are located, and the immigrants join, respectively.
  \end{itemize}
Jagers \cite{Jag} also pointed out that the reproduction of biological populations consisting
 of two types of individuals often displays the irreversibility property described above in the sense
 that individuals of one type might give birth to descendants of both kinds, whereas those of the other type
 can have descendants only of their own kind.
For example, if human diploid cells in a tumour are considered the first type in the cell population,
 and cells of higher diploidity are considered the second type, then, provided that endomitosis
 (a process where chromosomes duplicate but the cell does not subsequently divide, causing higher ploidity) is possible,
 the population of cells in this tumour has the irreversibility property in question.

 We can distinguish the following 5 cases \ ($a_{1,2}=0$ \ for each case):
\begin{align}\label{tablazat_esetek}
 \begin{tabular}{|c|c|c|c|}
  \hline
  \textup{(1)} & $a_{1,1} = 1$ & $a_{2,1} = 0$ & $a_{2,2} = 1$ \\
  \hline
  \textup{(2)} & $a_{1,1} = 1$ & $a_{2,1} \in \RR_{++}$ & $a_{2,2} = 1$ \\
  \hline
  \textup{(3)} & $a_{1,1} = 1$ & $a_{2,1} = 0$ & $a_{2,2} \in [0, 1)$ \\
  \hline
  \textup{(4)} & $a_{1,1} = 1$ & $a_{2,1} \in \RR_{++}$ & $a_{2,2} \in [0, 1)$ \\
  \hline
  \textup{(5)} & $a_{1,1} \in [0, 1)$ & $a_{2,1} \in \RR_+$ & $a_{2,2} = 1$ \\
  \hline
 \end{tabular}
\end{align}
For abbreviation, we can write the above five cases in matrix form as follows:
 \[
    \begin{bmatrix}
      1 & 0 \\
      0 & 1
     \end{bmatrix}_1 , \qquad
    \begin{bmatrix}
      1 & 0 \\
       ++ & 1
     \end{bmatrix}_2 , \qquad
    \begin{bmatrix}
      1 & 0 \\
      0 & 1_-
     \end{bmatrix}_3 , \qquad
    \begin{bmatrix}
      1 & 0 \\
       ++ & 1_-
     \end{bmatrix}_4 , \qquad
    \begin{bmatrix}
      1_- & 0 \\
      + & 1
     \end{bmatrix}_5 .
 \]
We remark that in the literature the cases (1) and (2) are called strongly critical due to \ $a_{1,1} = a_{2,2}=1$, \
 the other cases (3), (4) and (5) are critical, but not strongly critical, see, e.g., Foster and Ney \cite[page 13]{FosNey2}.

Note that the first coordinate process \ $(X_{k,1})_{k\in\ZZ_+}$ \ of \ $(\bX_k)_{k\in\ZZ_+}$ \ satisfies
 \begin{align}\label{GWI_def}
   X_{k,1} = \sum_{j=1}^{X_{k-1,1}} \xi_{k,j,1,1} + \vare_{k,1} , \qquad k \in \NN ,
 \end{align}
 hence \ $(X_{k,1})_{k\in\ZZ_+}$ \ is a single-type Galton-Watson process with immigration,
 which is critical in cases (1)--(4) and is subcritical in case (5) (due to \ $\EE(\xi_{1,1,1,1})=a_{1,1}$).

If the process \ $(X_{k,1})_{k\in\ZZ_+}$ \ given in \eqref{GWI_def} is critical, i.e.,
 \ $\EE(\xi_{1,1,1,1}) = a_{1,1} = 1$, \ then, by a result of Wei and Winnicki \cite{WW} (see also Theorem \ref{thm:critical}), we have
 \begin{align}\label{help_Wei_Winnicki1}
   (n^{-1} X_{\nt,1})_{t\in\RR_+} \distr (\cX_{t,1})_{t\in\RR_+} \qquad
  \text{as \ $n \to \infty$,}
 \end{align}
 where the limit process \ $(\cX_{t,1})_{t\in\RR_+}$ \ is the pathwise unique strong solution of the stochastic differential equation (SDE)
 \[
   \dd \cX_{t,1} = b_1 \, \dd t + \sqrt{v^{(1)}_{1,1} \, \cX_{t,1}^+} \, \dd \cW_{t,1} , \qquad t\in\RR_+, \qquad \cX_{0,1} = 0 ,
 \]
 where \ $(\cW_{t,1})_{t\in\RR_+}$ \ is a standard Wiener process, \ $b_1= \EE(\vare_{1,1})$ \ and \ $v^{(1)}_{1,1}:=\var(\xi_{1,1,1,1})$.
\ The process \ $(\cX_{t,1})_{t\in\RR_+}$ \ is called a squared Bessel process.

If the process \ $(X_{k,1})_{k\in\ZZ_+}$ \ given in \eqref{GWI_def} is subcritical,
 i.e., \ $\EE(\xi_{1,1,1,1}) = a_{1,1} \in [0, 1)$, \ then the Markov chain
 \ $(X_{k,1})_{k\in\ZZ_+}$ \ admits a unique stationary distribution \ $\mu_1$ \ (for its existence and generator function, see Appendix \ref{app:singletype_GWI})
  and, by Lemma \ref{lemma:subcritical}, we have
 \[
   (X_{\nt,1})_{t\in\RR_{++}} \distrf (\cX_{t,1})_{t\in\RR_{++}} \qquad
  \text{as \ $n \to \infty$,}
 \]
 where \ $(\cX_{t,1})_{t\in\RR_{++}}$ \ is an i.i.d.\ process (see the first paragraph of Section \ref{section_multi_branching})
 such that for each $t \in \RR_{++}$, the distribution of $\cX_{t,1}$ is $\mu_1$.
\ We note that the index set for the weak convergence of finite dimensional distributions above is \ $\RR_{++}$ \ and not \ $\RR_+$, \
 since \ $X_{0,1}=0$ \ not having the stationary distribution \ $\mu_1$ \ unless \ $\PP(\vare_{1,1} = 0) = 1$ \
 (for more details, see Appendix \ref{app:singletype_GWI}).
If \ $\PP(\vare_{1,1} = 0) = 1$, \ then \ $\PP(X_{n,1} = 0)=1$, $n\in\ZZ_+$ (due to $X_{0,1}=0$), and in this case
 the index set in question can be chosen as \ $\RR_+$ \ as well.

If \ $a_{2,1}=0$ \ holds as well, then \ $\PP(\xi_{1,1,1,2} =0)=1$ \ and
 \eqref{help1} yields that in this case the second coordinate process \ $(X_{k,2})_{k\in\ZZ_+}$ \ satisfies
 \[
   X_{k,2} = \sum_{j=1}^{X_{k-1,2}} \xi_{k,j,2,2} + \vare_{k,2} , \qquad k \in \NN .
 \]
Hence if \ $a_{2,1}=0$ \ holds as well, then \ $(X_{k,2})_{k\in\ZZ_+}$ \ is a single-type Galton-Watson process
 with immigration, which is critical in cases (1) and (5), and is subcritical in case (3)
 due to \ $\EE(\xi_{1,1,2,2})=a_{2,2}$.

Next we present our results on the asymptotic behaviour of \ $(\bX_k)_{k\in\ZZ_+}$ \ in the five cases (1)--(5) of its offspring mean matrix \ $\bA$.
\ The matrices \ $\bV^{(1)}$ \ and \ $\bV^{(2)}$ \ (introduced in Section \ref{section_multi_branching})
 in case of \ $p=2$ \ will be written in the form \ $\bV^{(1)} =: (v^{(1)}_{i,j})_{i,j=1}^2$ \ and \ $\bV^{(2)} =: (v^{(2)}_{i,j})_{i,j=1}^2$, \ respectively.

\begin{Thm}\label{main_1}
Let \ $(\bX_k)_{k\in\ZZ_+}$ \ be a critical decomposable 2-type Galton-Watson process with immigration such that \ $\bX_0=\bzero$, \
 the moment conditions \ $\EE(\|\bxi_i\|^4) < \infty$, \ $i=1,2$, \ and \ $\EE (\|\bvare\|^4) < \infty$ \
 hold and its offspring mean matrix $\bA$ satisfies \textup{(1)} of \eqref{tablazat_esetek}.
\ Then we have
 \begin{equation}\label{Conv_X_1}
  \Biggl(\begin{bmatrix}
          n^{-1} X_{\nt,1} \\
          n^{-1} X_{\nt,2}
         \end{bmatrix}\Biggr)_{t\in\RR_+}
  \distr
  \Biggl(\begin{bmatrix}
          \cX_{t,1} \\
          \cX_{t,2}
         \end{bmatrix}\Biggr)_{t\in\RR_+}
  \qquad \text{as \ $n \to \infty$,}
 \end{equation}
 where the limit process is the pathwise unique strong solution of the SDE
 \begin{equation}\label{SDE_X_1}
  \begin{cases}
   \dd\cX_{t,1}
   = b_1 \, \dd t
     + \sqrt{v^{(1)}_{1,1} \, \cX_{t,1}^+} \,
       \dd \cW_{t,1} , \\[2mm]
   \dd\cX_{t,2}
   = b_2 \, \dd t
     + \sqrt{v^{(2)}_{2,2} \, \cX_{t,2}^+} \,
       \dd \cW_{t,2} ,
  \end{cases}
  \qquad t \in\RR_+ ,
 \end{equation}
 with initial value \ $(\cX_{0,1} , \cX_{0,2}) = (0,0)$, \ where \ $(\cW_{t,1})_{t\in\RR_+}$ \
 and \ $(\cW_{t,2})_{t\in\RR_+}$ \ are independent standard Wiener processes yielding the independence of
 \ $(\cX_{t,1})_{t\in\RR_+}$ \ and \ $(\cX_{t,2})_{t\in\RR_+}$ \ as well.
\end{Thm}

\begin{Thm}\label{main_2}
Let \ $(\bX_k)_{k\in\ZZ_+}$ \ be a critical decomposable 2-type Galton-Watson process with immigration such that \ $\bX_0=\bzero$, \
 the moment condition \eqref{help4} holds and its offspring mean matrix \ $\bA$ \ satisfies \textup{(2)} of \eqref{tablazat_esetek}.
Then we have
 \begin{equation*}
  \Biggl(\begin{bmatrix}
          n^{-1} X_{\nt,1} \\
          n^{-2} X_{\nt,2}
         \end{bmatrix}\Biggr)_{t\in\RR_+}
  \distr
  \Biggl(\begin{bmatrix}
          \cX_{t,1} \\
          \cX_{t,2}
         \end{bmatrix}\Biggr)_{t\in\RR_+}
  \qquad \text{as \ $n \to \infty$,}
 \end{equation*}
 where the limit process is the pathwise unique strong solution of the SDE
 \begin{equation}\label{SDE_X_2}
  \begin{cases}
   \dd\cX_{t,1}
   = b_1 \, \dd t
     + \sqrt{v^{(1)}_{1,1} \, \cX_{t,1}^+} \,
       \dd \cW_{t,1} , \\[2mm]
   \dd\cX_{t,2}
   = a_{2,1} \cX_{t,1} \, \dd t ,
  \end{cases}
  \qquad t \in\RR_+ ,
 \end{equation}
 with initial value \ $(\cX_{0,1}, \cX_{0,2}) = (0,0)$, \ where \ $(\cW_{t,1})_{t\in\RR_+}$ \ is a standard Wiener process.
\end{Thm}

In Section \ref{Rem2}, we will compare our results in Theorem \ref{main_2} and those of Foster and Ney \cite[Section 9]{FosNey2} in detail.
Here we only note that Theorem \ref{main_2} can be considered as a functional version of Theorems 4 and 5 in Foster and Ney \cite{FosNey2}
 for some strongly critical decomposable 2-type Galton-Watson processes with immigration.
In Section \ref{Rem2}, as a consequence of Theorem \ref{main_2}, we also give a functional generalization of the Corollary on page 42 in Foster and Ney \cite{FosNey2},
which concerns the joint convergence of the appropriately normalized population size and total progeny
 of a critical single-type Galton--Watson process with immigration as \ $n\to\infty$.
\ We also give a stochastic representation of the limit process.

In the next corollary we describe the asymptotic behaviour of the relative frequency of individuals of type 2
 with respect to individuals of type 1 under the conditions of Theorem \ref{main_2} together with \ $b_1\in\RR_{++}$.
\ For different models, one can find similar results, e.g., in Jagers \cite[Corollary 1]{Jag}
 for supercritical decomposable age-dependent 2-type Galton-Watson processes without immigration,
 and in Yakovlev and Yanev \cite[Theorem 2]{YakYan} for some primitive multi-type Galton-Watson processes without immigration.

\begin{Cor}\label{Cor_relative_frequency}
Let us suppose that the conditions of Theorem \ref{main_2} together with \ $b_1\in\RR_{++}$ \ hold.
Then for all \ $t\in\RR_{++}$, \ we have
 \[
  n^{-1} \bbone_{\{ X_{\nt,1}\ne 0 \}}\frac{X_{\nt,2}}{X_{\nt,1}}
     \distr a_{2,1}\frac{\int_0^t \cX_{s,1}\,\dd s}{\cX_{t,1}}
     \qquad \text{as \ $n\to\infty$,}
 \]
 where \ $(\cX_{t,1})_{t\in\RR_+}$ \ is the pathwise unique strong solution of the first SDE in \eqref{SDE_X_2}
 with initial value \ $\cX_{0,1}=0$.
\end{Cor}

\begin{Thm}\label{main_3}
Let \ $(\bX_k)_{k\in\ZZ_+}$ \ be a critical decomposable 2-type Galton-Watson process with immigration such that \ $\bX_0=\bzero$, \
 the moment condition \eqref{help4} holds and its offspring mean matrix \ $\bA$ \ satisfies \textup{(3)} of \eqref{tablazat_esetek}.
Then we have
 \begin{equation*}
  (n^{-1} X_{\nt,1})_{t\in\RR_+}
  \distr
  (\cX_{t,1})_{t\in\RR_+}
  \qquad \text{as \ $n \to \infty$,}
 \end{equation*}
 where the limit process is the pathwise unique strong solution of the SDE
 \begin{equation}\label{SDE_X_3}
  \dd\cX_{t,1}
  = b_1 \, \dd t
    + \sqrt{v^{(1)}_{1,1} \, \cX_{t,1}^+} \,
    \dd \cW_{t,1} ,
  \qquad t \in\RR_+ ,
 \end{equation}
 with initial value \ $\cX_{0,1} = 0$, \ where \ $(\cW_{t,1})_{t\in\RR_+}$ \ is a standard Wiener process.
Further, the Markov chain \ $(X_{k,2})_{k\in\ZZ_+}$ \ admits a unique stationary distribution \ $\mu_2$ \ (for its existence and generator function, see the beginning of Appendix \ref{app:singletype_GWI}) and
 \begin{equation}\label{Conv_X_3_2}
  (X_{\nt,2})_{t\in\RR_{++}}
  \distrf
  (\cX_{t,2})_{t\in\RR_{++}}
  \qquad \text{as \ $n \to \infty$,}
 \end{equation}
 where \ $(\cX_{t,2})_{t\in\RR_{++}}$ \ is an i.i.d. process such that for each \ $t \in \RR_{++}$, \ the distribution of \ $\cX_{t,2}$ \ is \ $\mu_2$.
\ Moreover,
 \begin{align}\label{help13_cov1}
  \lim_{n_1\to\infty} \sup_{t_1,t_2\in\RR_+} \sup_{n_2\in\NN}
            \Big\vert \cov(n_1^{-1} X_{\lfloor n_1t_1\rfloor, 1},  X_{\lfloor n_2t_2\rfloor, 2} )\Big\vert
     = 0,
 \end{align}
 and
 \begin{align}\label{help13_cov1_mod}
  \lim_{n_2\to\infty} \cov(n_1^{-1} X_{\lfloor n_1t_1\rfloor, 1},  X_{\lfloor n_2t_2\rfloor, 2} )= 0,
  \qquad t_1,t_2\in\RR_+,  \; n_1\in\NN.
 \end{align}
\end{Thm}

We note that the index set for the weak convergence of finite dimensional distributions
 in \eqref{Conv_X_3_2} is \ $\RR_{++}$ \ and not \ $\RR_+$, \ since \ $X_{0,2} = 0$ \ not having the stationary distribution \ $\mu_2$
 \ unless \ $\PP(\vare_{1,2} = 0) = 1$ \ (for more details, see Appendix \ref{app:singletype_GWI}).

Note that, under the conditions of Theorem \ref{main_3}, if the two coordinates \ $\vare_{1,1}$ \ and \ $\vare_{1,2}$ \ of \ $\bvare_1$ \ are independent,
 then \ $(\cX_{t,1})_{t\in\RR_+}$ \ and \ $(\cX_{t,2})_{t\in\RR_{++}}$ \ are independent in Theorem \ref{main_3}, since
 in this special case the two coordinate processes \ $(X_{k,1})_{k\in\ZZ_+}$ \ and \ $(X_{k,2})_{k\in\ZZ_+}$ \ of \ $(\bX_k)_{k\in\ZZ_+}$ \
 are independent.
Motivated by this, \eqref{help13_cov1} and \eqref{help13_cov1_mod}, under the conditions of Theorem \ref{main_3}, we conjecture that
 \[
   \Biggl(\begin{bmatrix}
          n^{-1} X_{\nt,1} \\
          X_{\nt,2}
         \end{bmatrix}\Biggr)_{t\in\RR_{++}}
   \distrf
   \Biggl(\begin{bmatrix}
          \cX_{t,1} \\
          \cX_{t,2}
         \end{bmatrix}\Biggr)_{t\in\RR_{++}}
   \qquad \text{as \ $n \to \infty$,}
 \]
 where the driving process \ $(\cW_{t,1})_{t\in\RR_+}$ \ of  \ $(\cX_{t,1})_{t\in\RR_+}$ \ is independent of \ $(\cX_{t,2})_{t\in\RR_{++}}$, \
 yielding the independence of \ $(\cX_{t,1})_{t\in\RR_+}$ \ and \ $(\cX_{t,2})_{t\in\RR_{++}}$ \ as well.

\begin{Thm}\label{main_4}
Let \ $(\bX_k)_{k\in\ZZ_+}$ \ be a critical decomposable 2-type Galton-Watson process with immigration such that \ $\bX_0=\bzero$, \
 the moment conditions \ $\EE(\|\bxi_i\|^4) < \infty$, \ $i=1,2$, \ and \ $\EE (\|\bvare\|^4) < \infty$ \
 hold and its offspring mean matrix \ $\bA$ \ satisfies \textup{(4)} of \eqref{tablazat_esetek}.
Then we have
 \begin{equation*}
  \Biggl(\begin{bmatrix}
          n^{-1} X_{\nt,1} \\
          n^{-1} X_{\nt,2}
         \end{bmatrix}\Biggr)_{t\in\RR_+}
  \distr
  \Biggl(\begin{bmatrix}
          \cX_{t,1} \\
          \cX_{t,2}
         \end{bmatrix}\Biggr)_{t\in\RR_+}
  \qquad \text{as \ $n \to \infty$,}
 \end{equation*}
 where the limit process is the pathwise unique strong solution of the SDE
 \begin{equation}\label{SDE_X_4}
  \begin{cases}
   \dd\cX_{t,1}
   = b_1 \, \dd t
     + \sqrt{v^{(1)}_{1,1} \, \cX_{t,1}^+} \,
       \dd \cW_{t,1} , \\[2mm]
   \dd\cX_{t,2}
   = \frac{a_{2,1}}{1-a_{2,2}} \, \dd\cX_{t,1} ,
  \end{cases}
  \qquad t \in\RR_+ ,
 \end{equation}
 with initial value \ $(\cX_{0,1},\cX_{0,2}) = (0,0)$, \ where \ $(\cW_{t,1})_{t\in\RR_+}$ \ is a standard Wiener process.
\end{Thm}

\begin{Rem}
If the conditions of Theorem \ref{main_4} hold together with \ $a_{2,1}=1$ \ and \ $a_{2,2}=0$, \ then
 \ $X^{(1)}_{\nt,2} = \sum_{j=1}^{\nt-1} (M_{j,1} + b_1) = \sum_{j=1}^{\nt -1} (X_{j,1} - X_{j-1,1}) = X_{\nt-1,1}$ \ and
 \ $X^{(2)}_{\nt,2} = M_{\nt,2} + b_2$ \ for \ $n\in\NN$ \ and \ $t\in\RR_+$ \ (see \eqref{Mk1} and \eqref{Xdeco}).
\ So \eqref{help10} and Theorem \ref{main_4} yield that
 \begin{align*}
 \Biggl(\begin{bmatrix}
          n^{-1} X_{\nt,1} \\
          n^{-1} X_{\nt,2}
         \end{bmatrix}\Biggr)_{t\in\RR_+}
  = \Biggl(\begin{bmatrix}
          n^{-1} X_{\nt,1} \\
          n^{-1} X_{\nt-1,1} + n^{-1}(M_{\nt,2} + b_2)
         \end{bmatrix}\Biggr)_{t\in\RR_+}
  \distr
  \Biggl(\begin{bmatrix}
          \cX_{t,1} \\
          \cX_{t,1}
         \end{bmatrix}\Biggr)_{t\in\RR_+}
 \end{align*}
 as \ $n \to \infty$, \ where \ $(\cX_{t,1})_{t\in\RR_+}$ \ is the pathwise unique strong solution of the first SDE in \eqref{SDE_X_4}
 with initial value \ $\cX_{0,1}=0$.
\proofend
\end{Rem}

\begin{Thm}\label{main_5}
Let \ $(\bX_k)_{k\in\ZZ_+}$ \ be a critical decomposable 2-type Galton-Watson process with immigration such that \ $\bX_0=\bzero$, \
 the moment conditions \ $\EE(\|\bxi_i\|^4) < \infty$, \ $i=1,2$, \ and $\EE (\|\bvare\|^4) < \infty$
 hold and its offspring mean matrix $\bA$ satisfies \textup{(5)} of \eqref{tablazat_esetek}.
Then the Markov chain \ $(X_{k,1})_{k\in\ZZ_+}$ \ admits a unique stationary distribution \ $\mu_1$
 \ (for its existence and generator function, see the beginning of Appendix \ref{app:singletype_GWI}) and
 \begin{equation}\label{Conv_X_5_1}
  (X_{\nt,1})_{t\in\RR_{++}}
  \distrf
  (\cX_{t,1})_{t\in\RR_{++}}
  \qquad \text{as \ $n \to \infty$,}
 \end{equation}
 where \ $(\cX_{t,1})_{t\in\RR_{++}}$ \ is an i.i.d. process such that for each \ $t \in \RR_{++}$, \ the distribution of \ $\cX_{t,1}$ \ is \ $\mu_1$.
\ Further, we have
 \begin{equation*}
  (n^{-1} X_{\nt,2})_{t\in\RR_+}
  \distr
  (\cX_{t,2})_{t\in\RR_+}
  \qquad \text{as \ $n \to \infty$,}
 \end{equation*}
 where the limit process is the pathwise unique strong solution of the SDE
 \begin{equation}\label{SDE_X_5}
  \dd\cX_{t,2}
  = \biggl(\frac{a_{2,1}}{1-a_{1,1}} b_1 + b_2\biggr) \dd t
    + \sqrt{v^{(2)}_{2,2} \, \cX_{t,2}^+} \,
    \dd \cW_{t,2} ,
  \qquad t \in\RR_+ ,
 \end{equation}
 with initial value \ $\cX_{0,2} = 0$, \ where \ $(\cW_{t,2})_{t\in\RR_+}$ \ is a standard Wiener process.
Moreover,
 \begin{align}\label{help14_cov2}
   \lim_{n_2\to\infty} \sup_{t_1,t_2\in\RR_+} \sup_{n_1\in\NN}
            \Big\vert \cov(X_{\lfloor n_1t_1\rfloor, 1}, n_2^{-1}  X_{\lfloor n_2t_2\rfloor, 2} )\Big\vert
     = 0,
 \end{align}
 and
 \begin{align}\label{help14_cov2_mod}
   \lim_{n_1\to\infty}  \cov(X_{\lfloor n_1t_1\rfloor, 1}, n_2^{-1}  X_{\lfloor n_2t_2\rfloor, 2} )
     = 0, \qquad t_1,t_2\in\RR_+,\; n_2\in\NN.
 \end{align}
\end{Thm}

We remark that the index set for the weak convergence of finite dimensional distributions
 in \eqref{Conv_X_5_1} is \ $\RR_{++}$ \ and not \ $\RR_+$, \ since \ $X_{0,1} = 0$ \ not having the stationary distribution \ $\mu_1$ \
 unless \ $\PP(\vare_{1,1} = 0) = 1$ \ (for more details, see Appendix \ref{app:singletype_GWI}).
Further, note that the parameters \ $a_{1,1}$ \ and \ $b_1$ \ related to the first coordinate process
 \ $(X_{k,1})_{k\in\ZZ_+}$ \ appear in the drift coefficient of the SDE \eqref{SDE_X_5} for \ $(\cX_{t,2})_{t\in\RR_+}$,
 \ which is the limit process corresponding to the second coordinate process \ $(X_{k,2})_{k\in\ZZ_+}$.
\ It can be considered as a consequence of the decomposition \ $X_{k,2} = \sum_{\ell=1}^k (M_{\ell,2} + a_{2,1} X_{\ell-1,1} + b_2)$, \ $k\in\NN$
 \ (see \eqref{help17}), where \ $k^{-1} \sum_{\ell=1}^k X_{\ell-1,1}$ \ converges in probability to \ $b_1/(1-a_{1,1})$ \ as \ $k\to\infty$ \
 (see \eqref{Cond11_5}).
Moreover, note that if \ $a_{2,1}=0$ \ in Theorem  \ref{main_5} and if we switch the two coordinate processes, then
 we get back Theorem \ref{main_3} under fourth order moment assumptions on the offspring and immigration distributions.
The question of joint convergence of the two coordinate processes in Theorem \ref{main_5} remains open.

In the next remark we discuss the role of fourth order moment conditions in Theorems \ref{main_1}, \ref{main_4} and \ref{main_5}.

\begin{Rem}\label{Rem1}
We suspect that the moment conditions in Theorems \ref{main_1} and \ref{main_5} might be relaxed to
 \ $\EE(\|\bxi_i\|^2) < \infty$, \ $i=1,2$, \ and \ $\EE (\|\bvare\|^2) < \infty$ \ using the method of the proof of Theorem 3.1
 in Barczy et al.\ \cite{BarIspPap0}.
For Theorem \ref{main_5} in the special case \ $a_{2,1}=0$, \ it follows by Theorem  \ref{main_3} (by switching the two coordinate processes).
In fact, the fourth order moment assumptions in the proofs of Theorems \ref{main_1} and \ref{main_5} are used only for checking the
 conditional Lindeberg condition, namely, condition (iii) of Theorem \ref{Conv2DiffThm}, in order to prove convergence
 of some random step processes {towards a diffusion process.}
For single-type critical Galton-Watson processes with immigration, a detailed exposition of a proof of the conditional Lindeberg condition
 in question under second order moment assumptions can be found, e.g., in Barczy et al.\ \cite{BarBezPap}.
The fourth order moment conditions in Theorem \ref{main_4} come into play in another way, namely,
 via the estimation of tail probabilities of the maximum of a stable AR(1) process with heteroscedastic innovations
 \ $(M_{k,i})_{k\in\NN}$, \ $i=1,2$, \ that are martingale differences formed from the coordinate processes of the branching process in question.
Our technique is not suitable for relaxing them to second order ones, and we do not know any other technique (for more details,
 see the proof of Theorem \ref{main_4}).
\proofend
\end{Rem}

\section{Comparison of Theorem  \ref{main_2} and some results of Foster and Ney \cite{FosNey2}}\label{Rem2}

Under the conditions of Theorem \ref{main_2} together with \ $v^{(1)}_{1,1}\in\RR_{++}$, \
 Foster and Ney \cite[Theorems 4 and 5, and formula (9.11)]{FosNey2} proved that
 \begin{align}\label{FN_case2}
  \big( n^{-1} X_{n,1}, n^{-2} X_{n,2}\big) \distr (Y_1,Y_2) \qquad \text{as \ $n\to\infty$,}
 \end{align}
 where the Laplace transform of \ $(Y_1,Y_2)$ \ takes the form
 \begin{align}\label{help_Laplace_FN}
  \begin{split}
   &\EE\left( \ee^{-s_1 Y_1 -s_2 Y_2}\right) \\
   &\qquad =\exp\Bigg\{  -b_1 \int_0^1 \sqrt{\frac{2a_{1,2} s_1}{v^{(1)}_{1,1}}}
                   \frac{\frac{1}{2}v^{(1)}_{1,1} s_1 + \sqrt{\frac{1}{2}v^{(1)}_{1,1} a_{1,2}s_2}\tanh\Big(\tau \sqrt{\frac{1}{2}v^{(1)}_{1,1}a_{1,2}s_2}\Big)}
                        { \frac{1}{2}v^{(1)}_{1,1} s_1 \tanh\Big(\tau \sqrt{\frac{1}{2}v^{(1)}_{1,1} a_{1,2}s_2}\Big) + \sqrt{\frac{1}{2}v^{(1)}_{1,1} a_{1,2}s_2}  }\,\dd \tau\Bigg\}
  \end{split}
 \end{align}
 for \ $s_1\in\RR_+$ and $s_2\in\RR_{++}$, where we recall \ $b_1=\EE(\vare_{1,1})$ \ and \ $v^{(1)}_{1,1}=\var(\xi_{1,1,1,1})$.
\ They also derived the Laplace transforms of the marginal distributions \ $Y_1$ \ and \ $Y_2$, \ respectively.
Namely,
 \[
   \EE(\ee^{-s_1 Y_1}) = \left(1 + \frac{1}{2}v^{(1)}_{1,1}s_1\right)^{-\frac{2b_1}{v^{(1)}_{1,1}}},\qquad s_1\in\RR_+,
 \]
 yielding that \ $Y_1$ \ is Gamma-distributed with parameters \ $\frac{v^{(1)}_{1,1}}{2}$ \ and \ $\frac{2b_1}{v^{(1)}_{1,1}}$, \ and
 \begin{align*}
   \EE(\ee^{-s_2 Y_2}) =  \left( \cosh\left(\sqrt{\frac{1}{2} v^{(1)}_{1,1} a_{1,2} s_2}\right) \right)^{-\frac{2b_1}{v^{(1)}_{1,1}}},
     \qquad s_2\in\RR_+.
 \end{align*}
As a consequence of Theorem \ref{main_2}, the distribution of \ $(Y_1,Y_2)$ \ coincides with that of
 \ $(\cX_{1,1}, \cX_{1,2})=(\cX_{1,1}, a_{2,1}\int_0^1 \cX_{u,1}\,\dd u)$,
 \ where \ $(\cX_{t,1})_{t\in\RR_+}$ \ is given as the pathwise unique strong solution of the first SDE in \eqref{SDE_X_2}
 with initial value \ $\cX_{0,1}=0$.
\ We check that
  \begin{align}\label{help_Laplace_FN_sajat}
   \begin{split}
   \EE\Big(\ee^{-s_1 Y_1 - s_2 Y_2}\Big)
    & =\EE\Big(\ee^{-s_1 \cX_{1,1} - s_2 a_{2,1}\int_0^1 \cX_{u,1} \,\dd u}\Big)\\
    & = \left( \cosh\left(\sqrt{\frac{ v_{1,1}^{(1)} a_{2,1} s_2}{2}}\right)
        + \frac{s_1\sqrt{v_{1,1}^{(1)}}}{\sqrt{2 a_{2,1}s_2}} \sinh\left(\sqrt{\frac{ v_{1,1}^{(1)} a_{2,1} s_2}{2}}\right) \right)^{-\frac{2b_1}{v_{1,1}^{(1)}}},
   \end{split}
 \end{align}
  where \ $s_1\in\RR_+$ \ and \ $s_2\in\RR_{++}$, \ by giving a closed formula for the Laplace transform
  \ $\EE( \ee^{-s_1 Y_1 -s_2 Y_2})$ \ given in \eqref{help_Laplace_FN}.
First, recall that for all \ $\nu\in[-1,\infty)$, \ the pathwise unique strong solution of the SDE
 \begin{equation*}
  \begin{cases}
    \dd \cY_t = (2\nu + 2)\,\dd t + \sqrt{4 \cY_t^+}\,\dd \cW_t,\qquad t\in\RR_+,\\[1mm]
     \cY_0 = 0,
   \end{cases}
  \end{equation*}
 is called a squared Bessel process with parameter \ $\nu$, \ where \ $(\cW_t)_{t\in\RR_+}$ \ is a standard Wiener process,
 and for \ $\nu\in(-1,0)$, \ we have
 \begin{align}\label{Laplace_squared_Bessel_process}
   \EE\Big(\ee^{-\alpha \cY_t - \frac{\beta^2}{2} \int_0^t \cY_s\,\dd s}\Big)
     = \left( \cosh(\beta t) + \frac{2\alpha}{\beta} \sinh(\beta t) \right)^{-\nu-1}, \qquad t\in\RR_+,
 \end{align}
 where \ $\alpha\in\RR_+$ \ and \ $\beta\ne 0$, \ $\beta\in\RR$, \ see, e.g., Borodin and Salminen \cite[pages 76 and 135]{BorSal}.
Since \ $v_{1,1}^{(1)}\in\RR_{++}$, \ we can introduce \ $\widetilde \cX_{t,1}:= \frac{4}{v_{1,1}^{(1)}} \cX_{t,1}$, \ $t\in\RR_+$, \
 where \ $(\cX_{t,1})_{t\in\RR_+}$ \ is given as the pathwise unique strong solution of the first SDE in \eqref{SDE_X_2}
 with initial value \ $\cX_{0,1}=0$.
Then \ $(\widetilde \cX_{t,1})_{t\in\RR_+}$ \ is the pathwise unique strong solution of the SDE
 \begin{equation*}
 \begin{cases}
   \dd \widetilde \cX_{t,1} = \frac{4}{v_{1,1}^{(1)}} b_1\,\dd t + \sqrt{4\widetilde \cX_{t,1}}\,\dd \cW_t,\qquad t\in\RR_+, \\[1mm]
   \widetilde \cX_{0,1} = 0,
 \end{cases}
 \end{equation*}
 so \ $(\widetilde \cX_{t,1})_{t\in\RR_+}$ \ is a squared Bessel process with parameter \ $\frac{2b_1}{v_{1,1}^{(1)}}-1$,
 and using \eqref{Laplace_squared_Bessel_process} we have
 \[
   \EE\Big(\ee^{-\widetilde s_1\widetilde \cX_{t,1} - \widetilde s_2 \int_0^t \widetilde \cX_{u,1} \,\dd u}\Big)
     = \left( \cosh(\sqrt{2\widetilde s_2} t)
        + \frac{2 \widetilde s_1}{\sqrt{2\widetilde s_2}} \sinh(\sqrt{2\widetilde s_2}t) \right)^{-\frac{2b_1}{v_{1,1}^{(1)}}}, \qquad t\in\RR_+,
 \]
 where \ $t\in\RR_+$, \ $\widetilde s_1\in\RR_+$ \ and \ $\widetilde s_2\in\RR_{++}$.
This readily implies \eqref{help_Laplace_FN_sajat}.

For historical fidelity, we remark that Foster and Ney \cite{FosNey2} proved \eqref{FN_case2} using the finiteness of \ $\EE(\Vert \bvare\Vert)$
 \ instead of that of \ $\EE(\Vert \bvare\Vert^2)$ \ as we supposed in Theorem \ref{main_2}.

Finally, we give a functional generalization of the Corollary on page 42 in Foster and Ney \cite{FosNey2}.
Let us suppose that the conditions of Theorem \ref{main_2} hold.
\ Then, by the proof of this theorem (see \eqref{Conv_I_2}),
 \begin{equation}\label{help12}
  \Biggl(\begin{bmatrix}
          n^{-1} X_{\nt,1} \\
          n^{-2} \sum_{j=1}^{\nt -1}X_{j,1}
         \end{bmatrix}\Biggr)_{t\in\RR_+}
  \distr
  \Biggl(\begin{bmatrix}
          \cX_{t,1} \\
          \int_0^t \cX_{s,1} \, \dd s
         \end{bmatrix}\Biggr)_{t\in\RR_+}
 \end{equation}
 as \ $n \to \infty$, \ where \ $(\cX_{t,1})_{t\in\RR_+}$ \ is the pathwise unique strong solution of the first SDE in \eqref{SDE_X_2}
 with initial value \ $\cX_{0,1}=0$, \ since \ $X^{(1)}_{\nt,2} = \sum_{j=1}^{\nt -1} X_{j,1}$, \ $n\in\NN$, \ $t\in\RR_+$,
 where \ $X^{(1)}_{\nt,2}$ \ is defined by \eqref{Xdeco} and the equality in question follows by \eqref{help16}.
Note that, by \eqref{GWI_def}, \ $(X_{k,1})_{k\in\ZZ_+}$ \ is a single-type critical Galton-Watson process with
 immigration due to \ $a_{1,1}=1$, \ and \ $\sum_{j=1}^{\nt -1}X_{j,1}$ \ is the total progeny of individuals of type 1
 up to time \ $\nt -1$.
\ So \eqref{help12} gives us a functional generalization of the Corollary on page 42 in Foster and Ney \cite{FosNey2}
 together with a stochastic representation of the limit process as well.
We mention that in the considered special case, \ $\sum_{j=1}^{\nt -1}X_{j,1} = a_{2,1}^{-1}(X_{\nt,2} - X_{\nt,2}^{(2)})$, \ $t\in\RR_+$, \
 where \ $X_{\nt,2}^{(2)}$ \ is given by \eqref{Xdeco} and, by \eqref{sup_J_2},
 \ $\sup_{t\in[0,T]} |n^{-2} X^{(2)}_{\nt,2}| \stoch 0$ \ as \ $n \to \infty$ \ for all \ $T \in \RR_{++}$.
\ In the end, we remark that one can alternatively derive \eqref{help12} directly from the statement of Theorem \ref{main_2}.
Namely, if \ $\PP(\xi_{k,j,1,2} = 1)=1$, \ $\PP(\xi_{k,j,2,2} = 1)=1$
 and $\PP(\vare_{k,2} = 0)=1$ for each $k,j\in\NN$, then $X_{k,2} = X_{k-1,1} + X_{k-1,2} = \sum_{j=1}^{k-1} X_{j,1}$ for each $k\in\NN$
 almost surely, where $(X_{k,1})_{k\in\ZZ_+}$ is a single-type Galton-Watson process with immigration and $a_{2,1}=1$.
Consequently, if, in addition, \ $\EE(\xi_{1,1,1,1})=a_{1,1}=1$, \ then Theorem \ref{main_2} directly yields \eqref{help12}.

\section{Preliminaries for the proofs: decompositions}\label{Section_prel_proof}

Let us introduce the sequence
 \begin{equation}\label{Mk1}
  \begin{split}
  \begin{bmatrix} M_{k,1} \\ M_{k,2} \end{bmatrix}
  &:= \bM_k := \bX_k - \EE(\bX_k \mid \cF_{k-1}^{\bX})
   = \bX_k - \bA \bX_{k-1} - \bb \\
  &= \begin{bmatrix}
      X_{k,1} - a_{1,1} X_{k-1,1} - b_1 \\
      X_{k,2} - a_{2,1} X_{k-1,1} - a_{2,2} X_{k-1,2} - b_2
     \end{bmatrix} , \qquad
   k \in \NN ,
  \end{split}
 \end{equation}
 of martingale differences with respect to the filtration
 \ $(\cF_k^\bX)_{k\in\ZZ_+}$, \ where we used \eqref{mart}.
From \eqref{Mk1}, we obtain the recursion
 \begin{equation*}
  \bX_k = \bA \bX_{k-1} + \bM_k + \bb , \qquad k \in \NN ,
 \end{equation*}
 which together with \ $\bX_0=\bzero$ \ implies
 \begin{equation}\label{X}
   \bX_k = \sum_{j=1}^k \bA^{k-j} (\bM_j + \bb) , \qquad
   k \in \NN.
 \end{equation}
Indeed, since \ $\bX_0=\bzero$, \ we have \ $\bX_1 = \bM_1 + \bb$, \ and, by induction, for all \ $k\in\NN$,
 \begin{align*}
  \bX_{k+1} =\bA \bX_k + \bM_{k+1}+\bb
        &   =\bA \sum_{j=1}^k \bA^{k-j} (\bM_j + \bb) + \bM_{k+1}+\bb \\
        &   = \sum_{j=1}^k \bA^{k+1-j} (\bM_j + \bb) + \bM_{k+1}+\bb
            = \sum_{j=1}^{k+1} \bA^{k+1-j} (\bM_j + \bb) .
 \end{align*}
For each \ $\ell \in \NN$, \ we have
 \begin{equation}\label{m^ell}
   \bA^\ell
   = \begin{bmatrix}
      a_{1,1}^\ell & 0 \\
      a_{2,1} \sum\limits_{i=1}^{\ell} a_{2,2}^{i-1} a_{1,1}^{\ell-i}
       & a_{2,2}^\ell
     \end{bmatrix} ,
 \end{equation}
 where
 \[
   \sum\limits_{i=1}^{\ell} a_{2,2}^{i-1} a_{1,1}^{\ell-i}
   = \begin{cases}
      \frac{a_{1,1}^\ell - a_{2,2}^\ell }{a_{1,1}-a_{2,2}}
       & \text{if \ $a_{1,1} \ne a_{2,2}$,} \\[2mm]
      \ell a_{1,1}^{\ell-1} & \text{if \ $a_{1,1} = a_{2,2}$.}
     \end{cases}
 \]
Indeed, by induction, for all \ $\ell\in\NN$, \ we have
 \begin{align*}
   \bA^{\ell+1} & = \bA^\ell\bA
                 = \begin{bmatrix}
                      a_{1,1}^\ell & 0 \\
                      a_{2,1} \sum\limits_{i=1}^{\ell} a_{2,2}^{i-1} a_{1,1}^{\ell-i}
                       & a_{2,2}^\ell
                  \end{bmatrix}
                  \begin{bmatrix}
                     a_{1,1} & 0 \\
                     a_{2,1} & a_{2,2}
                  \end{bmatrix}
                 =  \begin{bmatrix}
                      a_{1,1}^{\ell+1} & 0 \\
                       a_{2,1} a_{1,1}\sum\limits_{i=1}^{\ell} a_{2,2}^{i-1} a_{1,1}^{\ell-i} + a_{2,1}a_{2,2}^\ell
                       & a_{2,2}^{\ell+1}
                  \end{bmatrix}  \\
                & =  \begin{bmatrix}
                    a_{1,1}^{\ell+1} & 0 \\
                    a_{2,1} \sum\limits_{i=1}^{\ell+1} a_{2,2}^{i-1} a_{1,1}^{\ell+1-i}
                     & a_{2,2}^{\ell+1}
                   \end{bmatrix},
 \end{align*}
 and, if \ $a_{1,1}\ne a_{2,2}$, \ then
 \begin{align*}
  \sum_{i=1}^{\ell} a_{2,2}^{i-1} a_{1,1}^{\ell-i}
     = \frac{a_{1,1}^\ell}{a_{2,2}}
       \sum_{i=1}^\ell \left(\frac{a_{2,2}}{a_{1,1}}\right)^i
     = \frac{a_{1,1}^\ell}{a_{2,2}}\cdot \frac{a_{2,2}}{a_{1,1}}
        \cdot \frac{\left(\frac{a_{2,2}}{a_{1,1}}\right)^\ell - 1}{\frac{a_{2,2}}{a_{1,1}} - 1}
     = \frac{a_{1,1}^\ell - a_{2,2}^\ell}{a_{1,1} - a_{2,2}} .
 \end{align*}
Note that \eqref{m^ell} holds for \ $\ell=0$ \ as well with the convention \ $\sum_{i=1}^0:=0$.
\ Consequently, by \eqref{X}, we get a decomposition
 \begin{equation}\label{Xdeco}
  \begin{bmatrix}
   X_{k,1} \\
   X_{k,2}
  \end{bmatrix}
  = \begin{bmatrix}
     X_{k,1} \\
     a_{2,1} X_{k,2}^{(1)} + X_{k,2}^{(2)}
    \end{bmatrix} , \qquad k \in \NN ,
 \end{equation}
 where
 \begin{align*}
  X_{k,1} &= \sum_{j=1}^k a_{1,1}^{k-j} (M_{j,1} + b_1) ,
 \end{align*}
 \begin{align*}
  X_{k,2}^{(1)}
  &:= \begin{cases}
       \sum\limits_{j=1}^k \frac{a_{1,1}^{k-j}-a_{2,2}^{k-j}}{a_{1,1}-a_{2,2}} (M_{j,1} + b_1)
       & \text{if \ $a_{1,1} \ne a_{2,2}$,} \\
       \sum\limits_{j=1}^k (k - j) a_{1,1}^{k-j-1} (M_{j,1} + b_1)
       & \text{if \ $a_{1,1} = a_{2,2}$,}
      \end{cases} , \\
  X_{k,2}^{(2)}
  &:= \sum_{j=1}^k a_{2,2}^{k-j} (M_{j,2} + b_2) .
 \end{align*}

\section{Proof of Theorem \ref{main_1}}
\label{Proof1}

The SDE \eqref{SDE_X_1} has a pathwise unique strong solution \ $(\bcX_t:=(\cX_{t,1},\cX_{t,2})^\top )_{t\in\RR_+}$ \ for all initial values
 \ $\bcX_0 = \bx \in \RR^2$, \ and if \ $\bx \in \RR_+^2$, \ then \ $\bcX_t \in \RR_+^2$ \ almost surely for all \ $t \in \RR_+$ \
  since \ $b_1, b_2, v^{(1)}_{1,1}, v^{(2)}_{2,2} \in \RR_+$, \ see, e.g., Ikeda and Watanabe \cite[Chapter IV, Example 8.2]{IkeWat}.
Since \ $a_{1,1} = 1$, \ by \eqref{help_Wei_Winnicki1}, we have \ $(n^{-1} X_{\nt,1})_{t\in\RR_+} \distr (\cX_{t,1})_{t\in\RR_+}$ \ as \ $n \to \infty$,
 \ where \ $(\cX_{t,1})_{t\in\RR_+}$ \ satisfies the first equation of the SDE \eqref{SDE_X_1} with initial value \ $\cX_{0,1}=0$.
\ Similarly, since \ $a_{2,1}=0$ \ and \ $a_{2,2}=1$, \ as it was explained in Section \ref{Section_conv_results}, the second coordinate process \ $(X_{k,2})_{k\in\ZZ_+}$ \ is
 a critical single-type Galton-Watson process with immigration, so \ $(n^{-1} X_{\nt,2})_{t\in\RR_+} \distr (\cX_{t,2})_{t\in\RR_+}$ \ as \ $n \to \infty$,
 \ where \ $(\cX_{t,2})_{t\in\RR_+}$ \ satisfies the second equation of the SDE \eqref{SDE_X_1} with initial value \ $\cX_{0,2}=0$.
\ However, we need to prove joint convergence of \ $(n^{-1} X_{\nt,1})_{t\in\RR_+}$ \ and \ $(n^{-1} X_{\nt,2})_{t\in\RR_+}$ \ as \ $n\to\infty$.

Using \ $a_{1,1} = 1$, \ $a_{2,1}=0$,  \ $a_{2,2} = 1$ \ and \eqref{Mk1},
 we obtain that the sequence \ $(\bM_k)_{k\in\NN}$ \ of martingale differences with respect to the filtration \ $(\cF_k^\bX)_{k\in\ZZ_+}$ \ takes the form
 \[
   \bM_k = \bX_k - \bX_{k-1} - \bb , \qquad k \in \NN .
 \]
Consider the random step processes
 \begin{align}\label{help15}
   \bcM_t^{(n)} := \begin{bmatrix} \cM_{t,1}^{(n)} \\ \cM_{t,2}^{(n)} \end{bmatrix}
   := \frac{1}{n} \sum_{k=1}^\nt \bM_k
   = \frac{1}{n} \bX_\nt - \frac{\nt}{n} \bb , \qquad t \in \RR_+ , \qquad n \in \NN ,
 \end{align}
 where we used that \ $\bX_0 = \bzero$.
\ We show that
 \begin{equation}\label{conv_bM}
  (\bcM_t^{(n)})_{t\in\RR_+} \distr (\bcM_t)_{t\in\RR_+} \qquad
  \text{as \ $n \to \infty$,}
 \end{equation}
 where the limit process \ $\bcM_t = (\cM_{t,1}, \cM_{t,2})^\top$, \ $t \in \RR_+$, \ is the pathwise unique strong solution of the SDE
 \begin{equation}\label{SDE_bM}
  \begin{cases}
   \dd \cM_{t,1}
   = \sqrt{v^{(1)}_{1,1} (\cM_{t,1} + b_1 t)^+ } \, \dd \cW_{t,1} , \qquad t\in\RR_+,\\
   \dd \cM_{t,2}
   = \sqrt{v^{(2)}_{2,2} (\cM_{t,2} + b_2 t)^+ } \, \dd \cW_{t,2} , \qquad t\in\RR_+
  \end{cases}
 \end{equation}
 with initial value \ $(\cM_{0,1}, \cM_{0,2} ) = (0,0)$.
\ In order to prove \eqref{conv_bM}, we want to apply Theorem \ref{Conv2DiffThm} for \ $d = r = 2$, \ $\bcU = \bcM$,
\ $\bU_k^{(n)} = n^{-1} \bM_k$, \ $n,k\in\NN$, \ $\bU_0^{(n)} =\bzero$, \ $n\in\NN$, \
 $\cF_k^{(n)} = \cF_k^\bX$, \ $n \in \NN$, \ $k \in \ZZ_+$ \ (yielding \ $\bcU^{(n)} = \bcM^{(n)}$, \ $n\in\NN$),
 \ and with coefficient functions \ $\bbeta : \RR_+ \times \RR^2 \to \RR^2$ \ and \ $\bgamma : \RR_+ \times \RR^2 \to \RR^{2\times2}$ \ of the SDE \eqref{SDE_bM} given by \ $\bbeta(t, \bx) = \bzero$ \ and
 \[
   \bgamma(t, \bx)
   = \begin{bmatrix}
      \sqrt{(x_1 + b_1 t)^+} & 0 \\
      0 & \sqrt{(x_2 + b_2 t)^+}
     \end{bmatrix}
     \begin{bmatrix}
      \sqrt{v^{(1)}_{1,1}} & 0 \\
      0 & \sqrt{v^{(2)}_{2,2}}
     \end{bmatrix}
 \]
 for \ $t \in \RR_+$ \ and \ $\bx = (x_1, x_2)^\top \in \RR^2$.
\ First we check that the SDE \eqref{SDE_bM} has a pathwise unique strong solution \ $(\bcM_t^{(\bx)})_{t\in\RR_+}$ \ for all initial values \ $\bcM_0^{(\bx)} = \bx \in \RR^2$.
\ Observe that if \ $(\bcM_t^{(\bx)})_{t\in\RR_+}$ \ is a strong solution of the SDE \eqref{SDE_bM} with initial value
 \ $\bcM_0^{(\bx)} = \bx \in \RR^2$, \ then, by It\^o's formula, the process \ $(\cP_{t,1}, \cP_{t,2})^\top := \bcM_t^{(x)} + \bb t$, \ $t \in \RR_+$, \
 is a pathwise unique strong solution of the SDE
 \begin{equation}\label{SDE_bP}
  \begin{cases}
   \dd \cP_{t,1}
   = b_1 \, \dd t
     + \sqrt{v^{(1)}_{1,1} \, \cP_{t,1}^+} \, \dd \cW_{t,1} , \qquad t\in\RR_+, \\
   \dd \cP_{t,2}
   = b_2 \, \dd t
     + \sqrt{v^{(2)}_{2,2} \, \cP_{t,2}^+} \, \dd \cW_{t,2}, \qquad t\in\RR_+,
  \end{cases}
 \end{equation}
 with initial value \ $(\cP_{0,1}, \cP_{0,2})^\top = \bx$.
\ Conversely, if \ $(\cP_{t,1}^{(\bp)}, \cP_{t,2}^{(\bp)})^\top$, \ $t \in \RR_+$, \ is a strong solution of the SDE \eqref{SDE_bP} with initial value \ $(\cP_{0,1}^{(\bp)}, \cP_{0,2}^{(\bp)})^\top = \bp \in \RR^2$, \ then, by It\^o's formula, the process \ $\bcM_t := (\cP_{t,1}^{(\bp)}, \cP_{t,2}^{(\bp)})^\top - \bb t$, \ $t \in \RR_+$, \ is a strong solution of the SDE \eqref{SDE_bM} with initial value \ $\bcM_0 = \bp$.
\ The equations in \eqref{SDE_bP} are the same as in \eqref{SDE_X_1}.
Consequently, as it was explained at the beginning of the proof,
 the SDE \eqref{SDE_bP} and hence the SDE \eqref{SDE_bM} as well admits a unique strong solution with arbitrary initial value in \ $\RR^2$,
 \ and \ $(\bcM_t + \bb t)_{t\in\RR_+} \distre (\bcX_t)_{t\in\RR_+}$.

The convergence \ $\bU_0^{(n)}\distr \bzero$ \ as \ $n\to\infty$, \ and condition (i) of Theorem \ref{Conv2DiffThm} trivially holds
 (since \ $\EE(\bM_k \mid \cF_{k-1}^\bX)=\bzero$, \ $k\in\NN$, \ and \ $\bbeta(t, \bx) = \bzero$, \ $t\in\RR_+$, \ $\bx\in\RR^2$).
\ Now, we show that conditions (ii) and (iii) of Theorem \ref{Conv2DiffThm} hold.
We have to check that for each \ $T \in \RR_{++}$,
 \begin{gather} \label{Condb1}
  \sup_{t\in[0,T]}
   \biggl\Vert\frac{1}{n^2}
          \sum_{k=1}^\nt
           \EE(\bM_k \bM_k^\top \mid \cF_{k-1}^\bX)
          - \int_0^t \bcR^{(n)}_s \,\bV_\bxi \, \dd s \biggr\Vert
  \stoch 0 \qquad \text{as \ $n\to\infty$,} \\
  \frac{1}{n^2}
  \sum_{k=1}^\nT
   \EE(\|\bM_k\|^2 \bbone_{\{\|\bM_k\|>n\theta\}} \mid \cF_{k-1}^\bX)
  \stoch 0  \qquad \text{as \ $n\to\infty$ \ for all \ $\theta \in \RR_{++}$,} \label{Condb2}
 \end{gather}
 where the process \ $(\bcR^{(n)}_s)_{s \in \RR_+}$ \ and the matrix \ $\bV_\bxi$ \ are
 defined by
 \begin{gather*}
  \bcR^{(n)}_s
  := \begin{bmatrix}
      (\cM_{s,1}^{(n)} + b_1 t)^+ & 0 \\
      0 & (\cM_{s,2}^{(n)} + b_2 t)^+
     \end{bmatrix} , \qquad s \in \RR_+ , \qquad n \in \NN , \\
  \bV_\bxi
  := \begin{bmatrix}
      v^{(1)}_{1,1} & 0 \\
      0 & v^{(2)}_{2,2}
     \end{bmatrix} .
 \end{gather*}
Indeed, for \ $t\in\RR_+$ \ and \ $\bx\in\RR^2$,
 \begin{align*}
  \gamma(t,\bx) \gamma(t,\bx)^\top
  & = \begin{bmatrix}
        \sqrt{(x_1 + b_1 t)^+} & 0 \\
        0 & \sqrt{(x_2 + b_2 t)^+}
       \end{bmatrix}\!\!
       \begin{bmatrix}
      v^{(1)}_{1,1} & 0 \\
      0 & v^{(2)}_{2,2}
     \end{bmatrix}\!\!
     \begin{bmatrix}
        \sqrt{(x_1 + b_1 t)^+} & 0 \\
        0 & \sqrt{(x_2 + b_2 t)^+}
       \end{bmatrix} \\
  & = \begin{bmatrix}
       v^{(1)}_{1,1} (x_1 + b_1 t)^+ & 0 \\
        0 & v^{(2)}_{2,2}(x_2 + b_2 t)^+
       \end{bmatrix}
    = \bV_\bxi
       \begin{bmatrix}
        (x_1 + b_1 t)^+ & 0 \\
        0 & (x_2 + b_2 t)^+
       \end{bmatrix}.
 \end{align*}

For each \ $s \in \RR_+$ \ and \ $n \in \NN$, \ we have
 \[
   \bcM_s^{(n)} + \bb s
   = \frac{1}{n} \bX_\ns + \frac{ns-\ns}{n} \bb ,
 \]
 thus
 \[
   \bcR^{(n)}_s
   = \begin{bmatrix}
      \cM_{s,1}^{(n)} + b_1 t & 0 \\
      0 & \cM_{s,2}^{(n)} + b_2 t
     \end{bmatrix} , \qquad s \in \RR_+ , \qquad n \in \NN ,
 \]
 and hence
 \begin{align*}
  \int_0^t \bcR^{(n)}_s \, \dd s
  &= \frac{1}{n^2} \sum_{k=0}^{\nt-1} \begin{bmatrix} X_{k,1} & 0 \\ 0 & X_{k,2} \end{bmatrix}
     + \frac{nt-\nt}{n^2} \begin{bmatrix} X_{\nt,1} & 0 \\ 0 & X_{\nt,2} \end{bmatrix} \\
  &\quad
     + \frac{\nt+(nt-\nt)^2}{2n^2}
        \begin{bmatrix} b_1 & 0 \\ 0 & b_2 \end{bmatrix} ,
   \qquad t \in \RR_+ , \qquad n \in \NN ,
 \end{align*}
 as, e.g., in the proof of Theorem 1.1 in Barczy et al.\ \cite{BarBezPap}.
By Lemma \ref{Moments},
 \[
   \frac{1}{n^2}
   \sum_{k=1}^\nt
    \EE(\bM_k \bM_k^\top \mid \cF_{k-1}^\bX)
   = \frac{\nt}{n^2} \bV^{(0)}
     + \frac{1}{n^2}
       \sum_{k=1}^\nt (X_{k-1,1} \bV^{(1)} + X_{k-1,2} \bV^{(2)})
 \]
 for all \ $t \in \RR_+$ \ and \  $n \in \NN$.
\ Since \ $\xi_{1,1,2,1} \ase 0$ \ and \ $\xi_{1,1,1,2} \ase 0$ \ (due to \ $a_{1,2} = a_{2,1}=0$), \ we have
 \ $v^{(2)}_{1,1} = v^{(2)}_{1,2} = v^{(2)}_{2,1} = 0$ \ and \ $v^{(1)}_{2,2} = v^{(1)}_{1,2} = v^{(1)}_{2,1} = 0$, \ and consequently
 \begin{align*}
  &X_{k-1,1} \bV^{(1)} + X_{k-1,2} \bV^{(2)}
   = X_{k-1,1} \begin{bmatrix}
                v^{(1)}_{1,1} & 0 \\
                0 & 0
               \end{bmatrix}
     + X_{k-1,2} \begin{bmatrix}
                  0 & 0 \\
                  0 & v^{(2)}_{2,2}
                 \end{bmatrix} \\
  &= \begin{bmatrix} X_{k-1,1} & 0 \\ 0 & X_{k-1,2} \end{bmatrix}
     \begin{bmatrix}
      v^{(1)}_{1,1} & 0 \\
      0 & v^{(2)}_{2,2}
     \end{bmatrix}
   = \begin{bmatrix} X_{k-1,1} & 0 \\ 0 & X_{k-1,2} \end{bmatrix}
     \bV_\bxi , \qquad k\in\NN.
 \end{align*}
So
 \begin{align*}
  & \frac{1}{n^2}
   \sum_{k=1}^\nt
    \EE(\bM_k \bM_k^\top \mid \cF_{k-1}^\bX)
     - \int_0^t \bcR^{(n)}_s \bV_\bxi \, \dd s   \\
  &\qquad  = \frac{\nt}{n^2}\bV^{(0)}
      - \frac{nt-\nt}{n^2}
       \begin{bmatrix} X_{\nt,1} & 0 \\ 0 & X_{\nt,2} \end{bmatrix}
       \bV_\bxi
       - \frac{\nt + (nt - \nt)^2}{2n^2} \begin{bmatrix} b_1 & 0 \\ 0 & b_2 \end{bmatrix}
         \bV_\bxi
 \end{align*}
 for \ $t\in\RR_+$ \ and \ $n\in\NN$.
\ Hence, in order to show \eqref{Condb1}, by Slutsky's lemma and taking into account the facts that
 for each \ $T\in\RR_{++}$,
 \[
   \sup_{t\in[0,T]}
      \frac{\lfloor nt\rfloor + (nt - \lfloor nt\rfloor)^2}{n^2}
    \leq \sup_{t\in[0,T]}  \frac{\lfloor nt\rfloor + 1}{n^2} \to 0
      \qquad \text{as \ $n\to\infty$,}
 \]
 and \ $\sup_{t\in[0,T]}\frac{\lfloor nt\rfloor}{n^2}\bV^{(0)}\as \bzero$ \ as \ $n\to\infty$, \ it suffices to prove
 that for each \ $T\in\RR_{++}$, \ we have
 \begin{equation}\label{Condb11}
  \frac{1}{n^2}
  \sup_{t \in [0,T]}
   \| (nt-\nt ) \bX_\nt\|
  \leq
  \frac{1}{n^2}
  \sup_{t \in [0,T]}
   \|\bX_\nt\|
  \stoch 0 \qquad \text{as \ $n \to \infty$.}
 \end{equation}
For each \ $k \in \NN$, \ we have \ $\bX_k = \bX_{k-1} + \bM_k + \bb$, \ and thus, using \ $\bX_0 = \bzero$,
 \begin{equation*}
  \bX_k = \sum_{j=1}^k \bM_j + k \bb ,
 \end{equation*}
 hence, for each \ $t \in \RR_+$ \ and \  $n \in \NN$, \ we get
 \[
   \|\bX_\nt\| \leq \sum_{j=1}^\nt \|\bM_j\| + \nt \|\bb\| .
 \]
Consequently, in order to prove \eqref{Condb11}, it suffices to show that for each \ $T\in\RR_{++}$,
 \[
   \frac{1}{n^2} \sum_{j=1}^\nT \|\bM_j\| \stoch 0 \qquad \text{as \ $n \to \infty$.}
 \]
By Lemma \ref{EEX1}, \ $\EE(X_{k,i}) = \OO(k)$ \ for \ $k\in\NN$, \ $i=1,2$, \ hence, by Lemma \ref{Moments}, we get
 \begin{align*}
  &\EE(\|\bM_j\|)
   \leq \sqrt{\EE(\|\bM_j\|^2)}
   = \sqrt{\EE(\bM_j^\top \bM_j)}
   = \sqrt{\EE(\tr(\bM_j^\top \bM_j))}
   = \sqrt{\EE(\tr(\bM_j \bM_j^\top))} \\
  &= \sqrt{\tr(\EE(\bM_j \bM_j^\top))}
   = \sqrt{\tr\big(\bV^{(0)} + \EE(X_{j-1,1}) \bV^{(1)} + \EE(X_{j-1,2}) \bV^{(2)}\big)}
   = \OO(j^{1/2}), \;\; j\in\NN.
 \end{align*}
Thus for each \ $T\in\RR_{++}$,
 \[
   \EE\Biggl(\frac{1}{n^2} \sum_{j=1}^\nT \|\bM_j\|\Biggr)
   = \frac{1}{n^2} \sum_{j=1}^\nT \OO(j^{1/2})
   = \OO(n^{-1/2}) \qquad \text{for \ $n \in \NN$,}
 \]
 and consequently we obtain \eqref{Condb11}, and hence \eqref{Condb1}.

Next, we check condition \eqref{Condb2}.
We show that for each \ $T\in\RR_{++}$ \ and \ $\theta\in\RR_{++}$,
 \[
   \frac{1}{n^2}
  \sum_{k=1}^\nT
   \EE(\|\bM_k\|^2 \bbone_{\{\|\bM_k\|>n\theta\}} \mid \cF_{k-1}^\bX)
   \mean 0 \qquad \text{as \ $n\to\infty$.}
 \]
By Markov's inequality and Lemma \ref{EEX1}, for each \ $T\in\RR_{++}$ \ and \ $\theta\in\RR_{++}$, \ we have
 \begin{align*}
  &\EE\left(  \frac{1}{n^2}
        \sum_{k=1}^\nT
         \EE(\|\bM_k\|^2 \bbone_{\{\|\bM_k\|>n\theta\}} \mid \cF_{k-1}^\bX)
     \right)
   = \frac{1}{n^2}
      \sum_{k=1}^\nT \EE(\|\bM_k\|^2 \bbone_{\{\|\bM_k\|>n\theta\}}) \\
  & \leq  \frac{1}{n^2}
          \sum_{k=1}^\nT \EE\left( \frac{\|\bM_k\|^4}{n^2\theta^2} \right)
   \leq \frac{2}{n^4 \theta^2} \sum_{k=1}^\nT \EE(M_{j,1}^4 + M_{j,2}^4 )
    =   \frac{1}{n^4\theta^2} \sum_{k=1}^\nT \OO(k^2)
    = \OO(n^{-1})\to 0
 \end{align*}
 as \ $n\to\infty$.

Using \eqref{conv_bM} and Lemma \ref{lemma:kallenberg}, we can prove \eqref{Conv_X_1}.
For each \ $n\in\NN$, \ by \eqref{help15}, we have \ $(n^{-1}\bX_{\nt})_{t\in\RR_+} = \Psi^{(n)}(\bcM^{(n)})$, \
 where the mapping \ $\Psi^{(n)} : \DD(\RR_+, \RR^2) \to \DD(\RR_+, \RR^2)$ \ is given by
 \[
   (\Psi^{(n)}(f))(t)
    := f\biggl(\frac{\nt}{n}\biggr) + \frac{\nt}{n} \bb
 \]
 for \ $f \in \DD(\RR_+, \RR^2)$ \ and \ $t \in \RR_+$.
Further, using that \ $(\bcM_t +\bb t)_{t\in\RR_+} \distre (\bcX_t)_{t\in\RR_+}$, \ we have
 \ $\bcX \distre \Psi(\bcM)$, \ where the mapping
 \ $\Psi : \DD(\RR_+, \RR^2) \to \DD(\RR_+, \RR^2)$ \ is given by
 \[
   (\Psi(f))(t) := f(t) + \bb t , \qquad
   f \in \DD(\RR_+, \RR^2) , \qquad t \in \RR_+.
 \]
The mappings \ $\Psi^{(n)}$, \ $n\in\NN$, \ and \ $\Psi$ \ are measurable, which can be checked in the same way
 as in Step 4/(a) in Barczy et al. \cite{BarBezPap} replacing \ $\DD(\RR_+,\RR)$ \ by \ $\DD(\RR_+,\RR^2)$
 \ in the argument given there.
One can also check that the set \ $C := \CC(\RR_+, \RR^2)$ \ satisfies \ $C \in \cB(\DD(\RR_+, \RR^2))$, \ $\PP(\bcM \in C) = 1$, \ and
 \ $\Psi^{(n)}(f^{(n)}) \to \Psi(f)$ \ in \ $\DD(\RR_+, \RR^2)$ \ as \ $n \to \infty$ \ if \ $f^{(n)} \to f$ \ in \ $\DD(\RR_+, \RR^2)$ \ as \ $n \to \infty$
 \ with \ $f \in C$, \ $f^{(n)}\in \DD(\RR_+, \RR^2)$, \ $n\in\NN$.
\ Namely, one can follow the same argument as in Step 4/(b) in Barczy et al.\ \cite{BarBezPap} replacing
 \ $\DD(\RR_+,\RR)$ \ by \ $\DD(\RR_+,\RR^2)$, \ and  \ $\CC(\RR_+,\RR)$ \ by \ $\CC(\RR_+,\RR^2)$, \ respectively,
 in the argument given there.
So we can apply Lemma \ref{lemma:kallenberg}, and we obtain
 \ $(n^{-1} \bX_\nt)_{t\in\RR_+} = \Psi^{(n)}(\bcM^{(n)}) \distr \Psi(\bcM)$ \ as \ $n \to \infty$, \ where
 \ $( (\Psi(\bcM))(t))_{t\in\RR_+}  = (\bcM_t + \bb t)_{t\in\RR_+} \distre (\bcX_t)_{t\in\RR_+}$, \ as desired.

\section{Proofs of Theorem \ref{main_2} and Corollary \ref{Cor_relative_frequency} }
\label{Proof2}

Since \ $a_{1,1} = 1$, \ by \eqref{help_Wei_Winnicki1}, we have \ $(n^{-1} X_{\nt,1})_{t\in\RR_+} \distr (\cX_{t,1})_{t\in\RR_+}$ \ as
 \ $n \to \infty$, \ where \ $(\cX_{t,1})_{t\in\RR_+}$ \ satisfies the first equation of the SDE \eqref{SDE_X_2} with \ $\cX_{0,1}=0$.
\ By \eqref{Xdeco}, we have the decomposition
 \[
   \bcX^{(n)}_t
    = \begin{bmatrix}
       \cX^{(n)}_{t,1} \\
       \cX^{(n)}_{t,2}
      \end{bmatrix}
   := \begin{bmatrix}
       n^{-1} X_{\nt,1} \\
       n^{-2} X_{\nt,2}
      \end{bmatrix}
   = \begin{bmatrix}
      n^{-1} X_{\nt,1} \\
      a_{2,1} \, n^{-2} X^{(1)}_{\nt,2} + n^{-2} X^{(2)}_{\nt,2}
     \end{bmatrix} , \qquad
   t \in \RR_+ , \quad n \in \NN ,
 \]
 where, since \ $a_{1,1}=a_{2,2}=1$, \
 \[
    X^{(1)}_{\nt,2} = \sum_{j=1}^{\lfloor nt\rfloor} (\lfloor nt\rfloor - j) (M_{j,1} + b_1)
 \]
 and
 \[
   X^{(2)}_{\nt,2} = \sum_{j=1}^{\lfloor nt\rfloor} (M_{j,2} + b_2).
 \]
Since \ $a_{1,1} = 1$, \ by \eqref{Mk1}, we obtain \ $M_{j,1} + b_1 = X_{j,1} - X_{j-1,1}$ \ for all \ $j \in \NN$, \ and hence, since \ $X_{0,1}=0$,
 \begin{align}\label{help16}
  \begin{split}
   n^{-2} X^{(1)}_{\nt,2}
    &= n^{-2} \sum_{j=1}^\nt (\nt - j) (X_{j,1} - X_{j-1,1})
      = n^{-2} \sum_{j=1}^\nt \sum_{k=1}^{\nt - j} (X_{j,1} - X_{j-1,1})\\
    & = n^{-2} \sum_{k=1}^{\nt - 1} \sum_{j=1}^{\nt -k} (X_{j,1} - X_{j-1,1})
      = n^{-2} \sum_{k=1}^{\nt - 1} (X_{\nt -k,1} - X_{0,1}) \\
    & = n^{-2} \sum_{j=1}^{\nt-1} X_{j,1}
      = \int_0^{\nt/n} \cX^{(n)}_{s,1} \, \dd s
  \end{split}
 \end{align}
 for all \ $t \in \RR_+$ \ and \ $n \in \NN$, \
 where the last equality follows by
 \begin{align*}
   \int_0^{\nt/n} \cX^{(n)}_{s,1} \, \dd s
    & = n^{-1}  \int_0^{\nt/n} X_{\ns,1} \, \dd s
      = n^{-1}  \sum_{j=1}^\nt \int_{\frac{j-1}{n}}^{\frac{j}{n}} X_{\ns,1} \, \dd s
      = n^{-1}  \sum_{j=1}^\nt \int_{\frac{j-1}{n}}^{\frac{j}{n}} X_{j-1,1} \, \dd s \\
    & = n^{-1}  \sum_{j=1}^\nt \frac{1}{n}  X_{j-1,1}
      = n^{-2}  \sum_{j=1}^\nt X_{j-1,1}
      = n^{-2}  \sum_{j=1}^{\nt-1} X_{j,1}.
 \end{align*}
By the continuous mapping theorem, we have
 \begin{equation}\label{Conv_I_2}
  \Biggl(\begin{bmatrix}
          n^{-1} X_{\nt,1} \\
          n^{-2} X^{(1)}_{\nt,2}
         \end{bmatrix}\Biggr)_{t\in\RR_+}
  = \Biggl(\begin{bmatrix}
            n^{-1} X_{\nt,1} \\
            \int_0^{\nt/n} \cX^{(n)}_{s,1} \, \dd s
           \end{bmatrix}\Biggr)_{t\in\RR_+}
  \distr
  \Biggl(\begin{bmatrix}
          \cX_{t,1} \\
          \int_0^t \cX_{s,1} \, \dd s
         \end{bmatrix}\Biggr)_{t\in\RR_+}
 \end{equation}
 as \ $n \to \infty$.
Indeed, \eqref{Conv_I_2} follows by an application of Lemma \ref{Conv2Funct} with the choices
 \ $d=1$, \ $q=2$ \ and \ $\bcU^{(n)}_t:=\cX^{(n)}_{t,1}$, \ $\bcU_t:=\cX_{t,1}$,
 \[
  (\Phi_n(f))(t):=\begin{bmatrix}
                    f(t) \\
                    \int_0^{\frac{\nt}{n}}  f(s)\,\dd s\\
                  \end{bmatrix},
    \qquad
  (\Phi(f))(t):=\begin{bmatrix}
                    f(t) \\
                    \int_0^t  f(s)\,\dd s\\
                  \end{bmatrix}
 \]
 for \ $n\in\NN$, \ $t\in\RR_+$, \ and \ $f\in \DD(\RR_+,\RR)$.
Next, we check that the conditions of Lemma \ref{Conv2Funct} with the above choices are satisfied, so we have right to apply Lemma \ref{Conv2Funct}.
\ The mapping \ $\Phi$ \ is continuous and hence measurable, see, e.g., Ethier and Kurtz \cite[Problem 3.11.8]{EthKur}.
For each \ $n\in\NN$, \ the mapping \ $\Phi_n$ \ is measurable as well, since \ $\Phi_n = \widetilde\Phi_n \circ \Phi$, \ where
 \ $\widetilde\Phi_n: \DD(\RR_+,\RR^2)\to \DD(\RR_+,\RR^2)$, \ $(\widetilde\Phi_n(g))(t):=\Big(g_1(t),g_2\Big(\frac{\lfloor nt\rfloor}{n}\Big) \Big)$, \ $t\in\RR_+$,
 \ $g=(g_1,g_2)\in \DD(\RR_+,\RR^2)$, \ is measurable (checked below), and composition of measurable mappings is measurable.
Using that the finite dimensional sets in \ $\DD(\RR_+,\RR^2)$ \ generate the Borel \ $\sigma$-algebra on \ $\DD(\RR_+,\RR^2)$ \
 (see, e.g., Jacod and Shiryaev \cite[Chapter VI, Theorem 1.14, part c)]{JacShi}), to check the measurability of \ $\widetilde\Phi_n$ \
 it is enough to verify that the mapping \ $\pi_t\circ \widetilde\Phi_n$ \ is measurable for each \ $t\in\RR_+$, \
 where \ $\pi_t: \DD(\RR_+,\RR^2)\to \RR^2$, \ $\pi_t(g):=g(t)$, \ $g\in \DD(\RR_+,\RR^2)$, \ is the natural projection onto \ $t$.
\ Since for each \ $n\in\NN$ \ and \ $t\in\RR_+$, \ $\pi_t\circ \widetilde\Phi_n : \DD(\RR_+,\RR^2)\to \RR^2$, \
 $(\pi_t \circ \widetilde\Phi_n)(g) = \Big( g_1(t), g_2\Big(\frac{\lfloor nt\rfloor}{n}\Big)\Big)^\top$, \ $g=(g_1,g_2)\in \DD(\RR_+,\RR^2)$,
 \ the first and second coordinate functions of \ $\pi_t \circ \widetilde\Phi_n$ \ can be identified with the natural projections of
 \ $\DD(\RR_+,\RR)$ \ onto the coordinate \ $t$ \ and \ $\frac{\lfloor nt\rfloor}{n}$, \ respectively,
 which are measurable (see, e.g., Billingsley \cite[Theorem 16.6, part (i)]{Bil}),
 yielding that \ $\pi_t \circ \widetilde\Phi_n$ \ is measurable as well.
We check that \ $C_{\Phi,(\Phi_n)_{n\in\NN}} = \CC(\RR_+,\RR)$.
\ For this we need to verify that \ $\Phi_n(f_n)\lu \Phi(f)$ \ as \ $n\to\infty$ \ whenever \ $f_n\lu f$ \ as \ $n\to\infty$ \
 with \ $f\in \CC(\RR_+,\RR)$ \ and \ $f_n\in \DD(\RR_+,\RR)$, \ $n\in\NN$.
\ For each \ $t\in\RR_+$, \ we have
 \begin{align*}
  \Vert (\Phi_n(f_n))(t) - (\Phi(f))(t)\Vert
    & \leq \vert f_n(t) - f(t)\vert
           + \left\vert \int_0^{\frac{\nt}{n}} f_n(s)\,\dd s - \int_0^t f(s)\,\dd s  \right\vert\\
    &\leq \vert f_n(t) - f(t)\vert
          + \int_0^{\frac{\nt}{n}} \vert f_n(s) - f(s)\vert \,\dd s
           + \int_{\frac{\nt}{n}}^t \vert f(s) \vert \,\dd s ,
 \end{align*}
 and hence for each \ $T\in\RR_{++}$,
 \begin{align*}
   \sup_{t\in[0,T]}
    \Vert (\Phi_n(f_n))(t) - (\Phi(f))(t) \Vert
     & \leq \sup_{t\in[0,T]}
            \vert f_n(t) - f(t)\vert
           +  \left(\sup_{t\in[0,T]} \vert f_n(t) - f(t)\vert\right) \frac{\nT}{n} \\
     &\phantom{\leq\,}
           + \sup_{t\in[0,T]} \vert f(t) \vert
              \sup_{t\in[0,T]} \left(t - \frac{\nt}{n}\right)
          \to 0
 \end{align*}
 as \ $n\to\infty$, \ since \ $f_n\lu f$ \ as \ $n\to\infty$ \ and \ $\sup_{t\in[0,T]} \vert f(t) \vert<\infty$ \
 (due to \ $f\in \CC(\RR_+,\RR)$).
All in all, the conditions of Lemma \ref{Conv2Funct} are satisfied with our choices, and hence we get \eqref{Conv_I_2}.

The aim of the following discussion is to show that
 \begin{equation}\label{sup_J_2}
  \sup_{t\in[0,T]} |n^{-2} X^{(2)}_{\nt,2}| \stoch 0 \qquad
  \text{as \ $n \to \infty$ \ for all \ $T \in \RR_{++}$.}
 \end{equation}
For each \ $t\in\RR_+$,
 \[
   |X^{(2)}_{\nt,2}|
   = \Biggl|\sum_{j=1}^\nt (M_{j,2} + b_2) \Biggr|
   \leq \Biggl|\sum_{j=1}^\nt M_{j,2}\Biggr| + b_2 \nt,
 \]
 hence, in order to check \eqref{sup_J_2}, it suffices to prove
 \[
   n^{-2} \sup_{t\in[0,T]} \Biggl|\sum_{j=1}^\nt M_{j,2}\Biggr|
   = n^{-2} \max_{k\in\{1,\ldots,\nT\}} \Biggl|\sum_{j=1}^k M_{j,2}\Biggr|
   \stoch 0 \qquad
   \text{as \ $n \to \infty$ \ for all \ $T \in \RR_{++}$,}
 \]
 which is equivalent to
 \begin{equation}\label{sup_J2_2}
  n^{-4} \max_{k\in\{1,\ldots,\nT\}} \Biggl(\sum_{j=1}^k M_{j,2}\Biggr)^2 \stoch 0
  \qquad \text{as \ $n \to \infty$ \ for all \ $T \in \RR_{++}$.}
 \end{equation}
Applying Doob's maximal inequality (see, e.g., Revuz and Yor
 \cite[Chapter II, Corollary (1.6)]{RevYor}) for the martingale
 \ $\sum_{j=1}^k M_{j,2}$, \ $k \in \NN$ \ (with the filtration
 \ $(\cF_k^\bX)_{k\in\NN}$), we obtain
 \begin{align*}
   \EE\left(\max_{k\in\{1,\ldots,\nT\}}  \left(\sum_{j=1}^k M_{j,2} \right)^2\right)
     &\leq 4 \EE\left(\left(\sum_{j=1}^\nT M_{j,2}\right)^2\right)\\
     &= 4\EE\left( \sum_{j=1}^\nT M_{j,2}^2 + 2 \sum_{j=1}^{\nT-1} \sum_{\ell=j+1}^\nT M_{j,2} M_{\ell,2}\right)\\
     &= 4\sum_{j=1}^\nT \EE(M_{j,2}^2)
        + 8 \sum_{j=1}^{\nT-1} \sum_{\ell=j+1}^\nT \EE(M_{j,2}M_{\ell,2})
      = 4 \sum_{j=1}^\nT \EE(M_{j,2}^2),
 \end{align*}
 since for each \ $j=1,\ldots,\nT-1$ \ and \ $\ell=j+1,\ldots,\nT$, \ we have
 \begin{align*}
  \EE(M_{j,2}M_{\ell,2})
    & = \EE\big( \EE(M_{j,2}M_{\ell,2} \mid \cF_{\ell-1}^\bX) \big)
      = \EE\big( M_{j,2} \EE(M_{\ell,2} \mid \cF_{\ell-1}^\bX) \big)\\
    & = \EE\big( M_{j,2} \EE(X_{\ell,2}  - \EE(X_{\ell,2} \mid \cF_{\ell-1}^\bX ) \mid \cF_{\ell-1}^\bX) \big)
     = \EE(M_{j,2}\cdot 0)
     =0.
 \end{align*}
Using Lemma \ref{EEX2}, we get \ $\sum_{j=1}^\nT \EE(M_{j,2}^2) = \sum_{j=1}^\nT \OO(j) = \OO(n^2)$ \ for \ $n\in\NN$ \ and \ $T\in\RR_+$,
 \ and consequently we obtain
 \[
   n^{-4} \max_{k\in\{1,\ldots,\nT\}} \left(\sum_{j=1}^k M_{j,2}\right)^2
      \mean 0 \qquad \text{as \ $n\to\infty$ \ for all \ $T\in\RR_{++}$,}
 \]
 yielding \eqref{sup_J2_2}, and hence \eqref{sup_J_2}.
Now, by Lemma VI.3.31 in Jacod and Shiryaev \cite{JacShi} (a kind of Slutsky's lemma for stochastic processes with trajectories
 in \ $\DD(\RR_+, \RR^d)$), \ convergences \eqref{Conv_I_2} and \eqref{sup_J_2} yield convergence
 \[
   (\bcX^{(n)}_t)_{t\in\RR_+} \distr (\bcX_t)_{t\in\RR_+} \qquad
   \text{as \ $n \to \infty$,}
 \]
 where the process \ $(\bcX_t)_{t\in\RR_+}$ \ is given by
 \[
   \bcX_t = \begin{bmatrix}
             \cX_{t,1} \\
             a_{2,1} \int_0^t \cX_{s,1} \, \dd s
            \end{bmatrix} , \qquad t \in \RR_+ .
 \]
By It\^o's formula, we obtain that \ $(\bcX_t)_{t\in\RR_+}$ \ satisfies the SDE \eqref{SDE_X_2} with initial value \ $\bcX_0=(0,0)$,
 \ thus we conclude the statement of Theorem \ref{main_2}.
\proofend

\bigskip

\noindent{\bf Proof of Corollary \ref{Cor_relative_frequency}.}
Since \ $b_1\in\RR_{++}$, \ we have \ $\PP(\cX_{t,1}\in\RR_{++})=1$ \ for each \ $t\in\RR_{++}$.
\ Indeed, if \ $v^{(1)}_{1,1}=0$, \ then \ $\cX_{t,1}= b_1 t$, \ $t\in\RR_+$, \ and if \ $v^{(1)}_{1,1}\in\RR_{++}$,
 \ then \ $\cX_{t,1}$ \ is Gamma distributed with parameters \ $2/v^{(1)}_{1,1}$ \ and \ $2b_1/v^{(1)}_{1,1}$
 \ for each \ $t\in\RR_{++}$, \ see, e.g., Ikeda and Watanabe \cite[Chapter IV, Example 8.2]{IkeWat}.
Let \ $g:\RR^2\to\RR$ \ be defined by
 \[
   g(x,y):= \bbone_{\{x\ne 0\}} \frac{y}{x}
          = \begin{cases}
             \frac{y}{x} & \text{if \ $x\ne 0$ \ and \ $y\in\RR$,}\\
             0 & \text{if \ $x=0$ \ and \ $y\in\RR$.}
           \end{cases}
 \]
Then \ $g$ \ is continuous on the set \ $(\RR\setminus \{0\})\times \RR$ \ and the distribution of \ $(\cX_{t,1},\cX_{t,2})$ \
 is concentrated on this set, since \ $\PP(\cX_{t,1}\in\RR_{++})=1$.
\ By Theorem \ref{main_2} and the continuous mapping theorem (see, e.g., Billingsley \cite[Theorem 5.1]{Bil0}),
 \begin{align*}
   & n^{-1}\bbone_{\{ X_{\nt,1}\ne 0 \}}\frac{X_{\nt,2}}{X_{\nt,1}}
      = g(n^{-1} X_{\nt,1}, n^{-2} X_{\nt,2}) \\
  &\qquad  \distr g(\cX_{t,1},\cX_{t,2})
             = \bbone_{\{ \cX_{t,1}\ne 0 \}} a_{2,1} \frac{\int_0^t \cX_{s,1}\,\dd s}{\cX_{t,1}}
             =  a_{2,1} \frac{\int_0^t \cX_{s,1}\,\dd s}{\cX_{t,1}}
      \qquad \text{as \ $n\to\infty$,}
 \end{align*}
 as desired.
\proofend

\section{Proof of Theorem \ref{main_3}}
\label{Proof3}

Since \ $a_{1,1} = 1$, \ by \eqref{help_Wei_Winnicki1},
 we have \ $(n^{-1} X_{\nt,1})_{t\in\RR_+} \distr (\cX_{t,1})_{t\in\RR_+}$ \ as \ $n \to \infty$, \ where \ $(\cX_{t,1})_{t\in\RR_+}$ \ satisfies the SDE \eqref{SDE_X_3}
 with \ $\cX_{0,1}=0$.

Since \ $a_{2,1} = 0$, \ as it was explained in Section \ref{Section_conv_results},
 the second coordinate process \ $(X_{k,2})_{k\in\ZZ_+}$ \ satisfies
 \[
   X_{k,2} = \sum_{j=1}^{X_{k-1,2}} \xi_{k,j,2,2} + \vare_{k,2} , \qquad k \in \NN .
 \]
Since \ $a_{2,2} \in [0, 1)$, \ the Markov chain \ $(X_{k,2})_{k\in\ZZ_+}$ \ admits a unique stationary distribution \ $\mu_2$ \
 (for its existence and generator function, see the beginning of Appendix \ref{app:singletype_GWI}), and, by Lemma \ref{lemma:subcritical},
  we have \ $(X_{\nt,2})_{t\in\RR_{++}} \distrf (\cX_{t,2})_{t\in\RR_{++}}$ \ as \ $n \to \infty$, \ where \ $(\cX_{t,2})_{t\in\RR_{++}}$ \ is an i.i.d. process such that for each \ $t \in \RR_{++}$, \ the distribution of \ $\cX_{t,2}$ \ is \ $\mu_2$.

Further, using \eqref{Xdeco} with \ $a_{1,1}=1$, \ \eqref{Cov}, and \ $\EE(\bM_k \mid \cF^{\bX}_{k-1})=\bzero$, \ $k\in\NN$ \
 (yielding \ $\cov(M_{j,1},M_{k,2})=0$, \ $j\ne k$, \ $j,k\in\NN$), \ we have for all \ $t_1,t_2\in\RR_+$,
 \begin{align*}
  &\cov\big(n_1^{-1} X_{\ntk,1}, X_{\ntn,2}\big)\\
  &\qquad = \cov\left( n_1^{-1}\left( \sum_{j=1}^{\ntk} M_{j,1} + \ntk b_1\right),
           \sum_{j=1}^{\ntn} a_{2,2}^{\ntn -j} M_{j,2} + \frac{a_{2,2}^{\ntn}-1}{a_{2,2}-1}b_2 \right)\\
 &\qquad = n_1^{-1}\cov\left(  \sum_{j=1}^{\ntk} M_{j,1}, \sum_{j=1}^{\ntn} a_{2,2}^{\ntn -j} M_{j,2} \right)
       = n_1^{-1} \sum_{j=1}^{\ntk\wedge \ntn} a_{2,2}^{\ntn -j} \cov(M_{j,1}, M_{j,2})\\
 &\qquad = n_1^{-1} \sum_{j=1}^{\ntk\wedge \ntn} a_{2,2}^{\ntn -j} \Big( v_{1,2}^{(0)} + \EE(X_{j-1,1}) v_{1,2}^{(1)} + \EE(X_{j-1,2}) v_{1,2}^{(2)}  \Big).
 \end{align*}
Since \ $a_{1,2}=a_{2,1}=0$, \ we have \ $\PP(\xi_{1,1,2,1}=0)=1$ \ and \ $\PP(\xi_{1,1,1,2}=0)=1$, \ and hence \ $v^{(2)}_{1,2}=0$ \
 and \ $v^{(1)}_{1,2}=0$, \ yielding that for all \ $t_1,t_2\in\RR_+$,
 \[
   \cov\big(n_1^{-1} X_{\ntk,1}, X_{\ntn,2}\big) = \frac{v^{(0)}_{1,2}}{n_1}  \sum_{j=1}^{\ntk\wedge \ntn} \!\!a_{2,2}^{\ntn -j}
                                           = \frac{v^{(0)}_{1,2}}{n_1}
                                             a_{2,2}^{\ntn - \ntk\wedge \ntn}\frac{1 - a_{2,2}^{\ntk\wedge \ntn}}{1-a_{2,2}},
 \]
 which yields \eqref{help13_cov1_mod} (due to \ $a_{2,2}\in[0,1)$ \ and \ $X_{0,2}=0$).
\ Using again \ $a_{2,2}\in[0,1)$, \ we have
 \begin{align*}
    \sup_{t_1,t_2\in\RR_+} \sup_{n_2\in\NN} \vert \cov\big(n_1^{-1} X_{\ntk,1}, X_{\ntn,2}\big) \vert
      &\leq \frac{v^{(0)}_{1,2}}{1-a_{2,2}} \frac{1}{n_1} \sup_{t_1,t_2\in\RR_+} \sup_{n_2\in\NN} a_{2,2}^{\ntn - \ntk\wedge \ntn}\\
      &\leq \frac{v^{(0)}_{1,2}}{1-a_{2,2}} \frac{1}{n_1} \to 0\qquad \text{as \ $n_1\to\infty$,}
 \end{align*}
 which yields \eqref{help13_cov1}.

\section{Proof of Theorem \ref{main_4}}
\label{Proof4}

Since \ $a_{1,1} = 1$, \ by \eqref{help_Wei_Winnicki1},
 we have \ $(n^{-1} X_{\nt,1})_{t\in\RR_+} \distr (\cX_{t,1})_{t\in\RR_+}$ \ as \ $n \to \infty$, \
 where \ $(\cX_{t,1})_{t\in\RR_+}$ \ satisfies the first equation of the SDE \eqref{SDE_X_4} with \ $\cX_{0,1}=0$.
\ By \eqref{Xdeco}, we have the decomposition
 \begin{align}\label{help10}
   \bcX^{(n)}_t
   = \begin{bmatrix}
       \cX_{t,1}^{(n)} \\
       \cX_{t,2}^{(n)}
      \end{bmatrix}
   := \begin{bmatrix}
       n^{-1} X_{\nt,1} \\
       n^{-1} X_{\nt,2}
      \end{bmatrix}
   = \begin{bmatrix}
       n^{-1} X_{\nt,1} \\
       a_{2,1} \, n^{-1} X^{(1)}_{\nt,2} + n^{-1} X^{(2)}_{\nt,2}
      \end{bmatrix} , \qquad
   t \in \RR_+ , \quad n \in \NN ,
 \end{align}
 where, since \ $a_{1,1}=1$ \ and \ $a_{2,2}\in[0,1)$,
 \begin{align*}
  X_{\nt,2}^{(1)}
  :=  \sum_{j=1}^\nt \frac{1-a_{2,2}^{\nt-j}}{1-a_{2,2}} (M_{j,1} + b_1)
  \qquad \text{and} \qquad
 X_{\nt,2}^{(2)}
  := \sum_{j=1}^\nt a_{2,2}^{\nt-j} (M_{j,2} + b_2) .
 \end{align*}
The aim of the following discussion is to show that
 \begin{equation}\label{sup_J_4}
  \sup_{t\in[0,T]} |n^{-1} X^{(2)}_{\nt,2}| \stoch 0 \qquad
  \text{as \ $n \to \infty$ \ for all \ $T \in \RR_{++}$.}
 \end{equation}
We have
 \begin{align*}
   |n^{-1} X^{(2)}_{\nt,2}|
   \leq \frac{1}{n}
        \Biggl|\sum_{j=1}^\nt
               a_{2,2}^{\nt-j} M_{j,2}\Biggr|
        + \frac{b_2}{n}
          \frac{1-a_{2,2}^\nt}
               {1-a_{2,2}}
   \leq \frac{1}{n}
        \Biggl|\sum_{j=1}^\nt
               a_{2,2}^{\nt-j} M_{j,2}\Biggr|
        + \frac{b_2}{n(1-a_{2,2})},
 \end{align*}
 hence, in order to check \eqref{sup_J_4}, it suffices to prove
 \begin{equation}\label{sup_J_42}
  n^{-1} \sup_{t\in[0,T]}
          \Biggl|\sum_{j=1}^\nt
                  a_{2,2}^{\nt-j} M_{j,2}\Biggr|
  = n^{-1} \max_{k\in\{1,\ldots,\nT\}} |V_{k,2}|
  \stoch 0
 \end{equation}
 as \ $n \to \infty$ \ for all \ $T \in \RR_{++}$, \ where
 \[
   V_{k,2}
   := \sum_{j=1}^k a_{2,2}^{k-j} M_{j,2} , \qquad
   k \in \NN .
 \]
Note that
 \begin{align}\label{help2}
   V_{k,2} = a_{2,2} V_{k-1,2} + M_{k,2} , \qquad k \in \NN ,
 \end{align}
 where \ $V_{0,2} := 0$, \ hence \ $(V_{k,2})_{k\in\ZZ_+}$ \ is a stable AR(1) process with heteroscedastic innovation \ $(M_{k,2})_{k\in\NN}$.
\ For all \ $\delta > 0$, \ by Markov's inequality, we have
 \begin{align*}
  \PP\Bigl( n^{-1} \max_{k\in\{1,\ldots,\nT\}} |V_{k,2}| > \delta\Bigr)
  &= \PP\Bigl(\max_{k\in\{1,\ldots,\nT\}} V_{k,2}^4 > \delta^4 n^4 \Bigr)
   \leq \sum_{k=1}^\nT \PP(V_{k,2}^4 > \delta^4 n^4) \\
  &\leq \delta^{-4} n^{-4} \sum_{k=1}^\nT \EE(V_{k,2}^4)
   = \delta^{-4} n^{-4} \sum_{k=1}^\nT \OO(k^2)
   = \OO(n^{-1}) ,\qquad n\in\NN,
 \end{align*}
 where we applied \ $\EE(V_{k,2}^4) = \OO(k^2)$, \ $k\in\NN$ \ (see Lemma \ref{EEX4}), thus we obtain
 \eqref{sup_J_42}, and hence \eqref{sup_J_4}.
Here we note that in the previous application of Markov's inequality, we took the fourth moment of \ $\vert V_{k,2}\vert$, \
 since even if we took its second moment, then we could not argue that \ $n^{-2} \sum_{k=1}^\nT \EE(V_{k,2}^2)\to 0$ \ as \ $n\to\infty$
 \ due to the fact that \ $\EE(V_{k,2}^2) = \OO(k)$, \ $k\in\NN$ \ (which can be checked similarly as \ $\EE(V_{k,2}^4) = \OO(k^2)$, \ $k\in\NN$,
 \ checked in the proof of Lemma \ref{EEX4}).

Recall that
 \[
   n^{-1} X^{(1)}_{\nt,2}
   = n^{-1}
     \sum_{j=1}^\nt
      \frac{1-a_{2,2}^{\nt-j}}
           {1-a_{2,2}}
      (M_{j,1} + b_1) , \qquad t \in \RR_+ , \quad n \in \NN .
 \]
Since \ $a_{1,1} = 1$, \ by \eqref{Mk1}, we obtain \ $M_{j,1} + b_1 = X_{j,1} - X_{j-1,1}$ \ for all \ $j \in \NN$, \ and
 hence, using that \ $X_{0,1}=0$, \ we have \ $\sum_{j=1}^k (M_{j,1} + b_1) = \sum_{j=1}^k (X_{j,1} - X_{j-1,1}) = X_{k,1}$ \ for all \ $k \in \NN$.
\ Consequently, we get
 \[
   n^{-1} X^{(1)}_{\nt,2}
   = \frac{n^{-1}X_{\nt,1}}{1-a_{2,2}}
     - \frac{n^{-1}}{1-a_{2,2}}
       \sum_{j=1}^\nt a_{2,2}^{\nt-j} (M_{j,1} + b_1)
 \]
 for all \ $t \in \RR_+$ \ and \ $n \in \NN$.
\ In a similar way as above for \eqref{sup_J_42}, using Lemma \ref{EEX4}, we have
 \begin{align}\label{help7}
   n^{-1} \sup_{t\in[0,T]}
           \Biggl|\sum_{j=1}^\nt
                   a_{2,2}^{\nt-j} M_{j,1}\Biggr|
   = n^{-1} \max_{k\in\{1,\ldots,\nT\}} |V_{k,1}|
   \stoch 0
 \end{align}
 as \ $n \to \infty$ \ for all \ $T \in \RR_{++}$, \ where
 \[
   V_{k,1}
   := \sum_{j=1}^k a_{2,2}^{k-j} M_{j,1} , \qquad
   k \in \NN ,
 \]
 which satisfies
 \begin{align}\label{help3}
  V_{k,1} = a_{2,2} V_{k-1,1} + M_{k,1}, \qquad k\in\NN,
 \end{align}
 where \ $V_{0,1}:=0$.
\ Moreover, since \ $a_{2,2}\in[0,1)$, \ we have
 \begin{align*}
   n^{-1} \sup_{t\in[0,T]}
           \Biggl|\sum_{j=1}^\nt
                   a_{2,2}^{\nt-j}\Biggr|
  & = n^{-1} \sup_{t\in[0,T]}  \Biggl|\sum_{j=0}^{\nt-1}  a_{2,2}^j \Biggr|
    = n^{-1} \sum_{j=0}^{\nT-1} a_{2,2}^j \\
  & \leq n^{-1} \sum_{\ell=0}^\infty a_{2,2}^\ell
   = \frac{n^{-1}}{1-a_{2,2}}
   \to 0 \qquad \text{as \ $n \to \infty$.}
 \end{align*}
By Lemma VI.3.31 in Jacod and Shiryaev \cite{JacShi} (a kind of Slutsky's lemma
 for stochastic processes with trajectories in \ $\DD(\RR_+, \RR^d)$) \ and \eqref{help10},
 the above convergences yield
 \[
   (\bcX^{(n)}_t)_{t\in\RR_+} \distr (\bcX_t)_{t\in\RR_+} \qquad
   \text{as \ $n \to \infty$,}
 \]
 where the process \ $(\bcX_t)_{t\in\RR_+}$ \ is given by
 \[
   \bcX_t = \begin{bmatrix}
             \cX_{t,1} \\
             \frac{a_{2,1}}
                  {1-a_{2,2}}
             \cX_{t,1}
            \end{bmatrix} , \qquad t \in \RR_+ .
 \]

\section{Proof of Theorem \ref{main_5}}
\label{Proof5}

Since \ $a_{1,1} \in [0, 1)$, \ the Markov chain \ $(X_{k,1})_{k\in\ZZ_+}$ \ admits a unique stationary distribution \ $\mu_1$
 (for its existence and generator function, see the beginning of Appendix \ref{app:singletype_GWI}), and, by Theorem \ref{thm:critical},
 \ $(X_{\nt,1})_{t\in\RR_{++}} \distrf (\cX_{t,1})_{t\in\RR_{++}}$ \ as \ $n \to \infty$, \ where \ $(\cX_{t,1})_{t\in\RR_{++}}$ \ is an i.i.d. process such that for each \ $t \in \RR_{++}$, \ the distribution of \ $\cX_{t,1}$ \ is \ $\mu_1$.

Since \ $a_{2,2} = 1$, \ by \eqref{Mk1}, we have
 \[
   M_{k,2} = X_{k,2} - a_{2,1} X_{k-1,1} - X_{k-1,2} - b_2 , \qquad k \in \NN ,
 \]
 hence, using that \ $X_{0,2}=0$,
 \begin{align}\label{help17}
   \cX_{t,2}^{(n)}
   := n^{-1} X_{\nt,2}
   = \frac{1}{n} \sum_{k=1}^\nt (X_{k,2} - X_{k-1,2})
   = \frac{1}{n} \sum_{k=1}^\nt U_{k,2} , \qquad t \in \RR_+ , \qquad n \in \NN ,
 \end{align}
 where
 \[
   U_{k,2} := M_{k,2} + a_{2,1} X_{k-1,1} + b_2 , \qquad k \in \NN .
 \]
We show that
 \begin{equation}\label{conv_X2_5}
  (\cX_{t,2}^{(n)})_{t\in\RR_+} \distr (\cX_{t,2})_{t\in\RR_+} \qquad
  \text{as \ $n \to \infty$,}
 \end{equation}
 where the limit process \ $(\cX_{t,2})_{t\in\RR_+}$ \ is the unique strong solution of the SDE \eqref{SDE_X_5} with \ $\cX_{0,2}=0$.
\ In order to prove \eqref{conv_X2_5}, we want to apply Theorem \ref{Conv2DiffThm} for \ $d = r = 1$,
 \ $(\bcU_t)_{t\in\RR_+} = (\cX_{t,2})_{t\in\RR_+}$, \ $\bU_k^{(n)} = n^{-1} U_{k,2}$, \ $n,k\in\NN$,
 \ $\bU_0^{(n)} = 0$, \ $n\in\NN$, \ $\cF_k^{(n)} = \cF_k^\bX$ \ for \ $n \in \NN$ \ and \ $k \in \ZZ_+$, \
 and with coefficient functions \ $\bbeta : \RR_+ \times \RR \to \RR$ \ and \ $\bgamma : \RR_+ \times \RR \to \RR$ \ of the SDE \eqref{SDE_X_5} given by
 \[
   \bbeta(t, x) = \frac{a_{2,1}}{1-a_{1,1}} b_1 + b_2 , \qquad
   \bgamma(t, x) = \sqrt{v^{(2)}_{2,2}\, x^+} , \qquad t \in \RR_+ , \quad x \in \RR .
 \]
The SDE \eqref{SDE_X_5} has a pathwise unique strong solution \ $(\cX_{t,2}^{(x)})_{t\in\RR_+}$ \ for all initial values \ $\cX_{0,2}^{(x)} = x \in \RR$, \ and if \ $x \in \RR_+$, \ then \ $\cX_{t,2}^{(x)} \in \RR_+$ \ almost surely for all \ $t \in \RR_+$ \ since
 \[
   \frac{a_{2,1}}{1-a_{1,1}} b_1 + b_2 \in \RR_+ , \qquad
   v^{(2)}_{2,2} \in \RR_+ ,
 \]
 see, e.g., Ikeda and Watanabe \cite[Chapter IV, Example 8.2]{IkeWat}.

Now, we show that conditions (i), (ii) and (iii) of Theorem \ref{Conv2DiffThm} hold.
We have to check that for each \ $T \in \RR_{++}$,
 \begin{gather}\label{Cond1_5}
  \sup_{t\in[0,T]}
   \biggl|\frac{1}{n}
          \sum_{k=1}^\nt
           \EE(U_{k,2} \mid \cF_{k-1}^\bX)
          - \biggl(\frac{a_{2,1}}{1-a_{1,1}} b_1 + b_2\biggr) t\biggr|
  \stoch 0 \qquad \text{as \ $n \to \infty$,} \\
  \sup_{t\in[0,T]}
   \biggl|\frac{1}{n^2}
          \sum_{k=1}^\nt
           \var(U_{k,2} \mid \cF_{k-1}^\bX)
          - \int_0^t (\cX^{(n)}_{s,2})^+ v^{(2)}_{2,2} \, \dd s\biggr|
  \stoch 0 \qquad \text{as \ $n \to \infty$,} \label{Cond2_5} \\
  \frac{1}{n^2}
  \sum_{k=1}^\nT
   \EE(U_{k,2}^2 \bbone_{\{|U_{k,2}|>n\theta\}} \mid \cF_{k-1}^\bX)
  \stoch 0  \qquad \text{as \ $n\to\infty$ \ for all \ $\theta \in \RR_{++}$.} \label{Cond3_5}
 \end{gather}

For each \ $k \in \NN$, \ we have \ $\EE(U_{k,2} \mid \cF_{k-1}^\bX) = a_{2,1} X_{k-1,1} + b_2$, \
 and \ $\sup_{t\in[0,T]} \Big\vert \frac{\nt}{n} -t \Big\vert\to 0$ \ as \ $n\to\infty$ \ for each \ $T\in\RR_{++}$,
 \ hence, in order to show \eqref{Cond1_5}, it suffices to prove that  for each \ $T\in\RR_{++}$,
 \begin{equation}\label{Cond11_5}
  \sup_{t\in[0,T]}
   \biggl|\frac{1}{n} \sum_{k=1}^\nt X_{k-1,1}
          - \frac{b_1}{1-a_{1,1}} t\biggr|
  \stoch 0 \qquad \text{as \ $n \to \infty$.}
 \end{equation}
By \eqref{Mk1}, we have
 \[
   M_{k,1} = X_{k,1} - a_{1,1} X_{k-1,1} - b_1 , \qquad \text{i.e.,}\qquad X_{k,1} = a_{1,1} X_{k-1,1}+M_{k,1} + b_1, \qquad k \in \NN ,
 \]
 hence, using that \ $X_{0,1}=0$, \ we have
 \[
   X_{k,1} = \sum_{j=1}^k a_{1,1}^{k-j} (M_{j,1} + b_1) , \qquad k \in \NN ,
 \]
 thus
 \begin{align*}
  \sum_{k=1}^\nt X_{k-1,1}
  &= \sum_{k=1}^\nt \sum_{j=1}^{k-1} a_{1,1}^{k-1-j} (M_{j,1} + b_1)
   = \sum_{j=1}^{\nt-1} \sum_{k=j+1}^\nt a_{1,1}^{k-1-j} (M_{j,1} + b_1) \\
  &= \sum_{j=1}^{\nt-1} \frac{1-a_{1,1}^{\nt-j}}{1-a_{1,1}} (M_{j,1} + b_1) \\
  &= \frac{1}{1-a_{1,1}} \sum_{j=1}^{\nt-1} (1-a_{1,1}^{\nt-j}) M_{j,1}
     + \frac{b_1}{1-a_{1,1}} \sum_{j=1}^{\nt-1} (1-a_{1,1}^{\nt-j}) \\
  &= \frac{1}{1-a_{1,1}} \sum_{j=1}^{\nt-1} (1-a_{1,1}^{\nt-j}) M_{j,1}
     + \frac{b_1}
            {1-a_{1,1}}
       \biggl(\nt - 1
              - \frac{a_{1,1}-a_{1,1}^\nt}
                     {1-a_{1,1}}\biggr) .
 \end{align*}
Hence for each \ $t\in\RR_+$ \ and \ $n\in\NN$,
 \begin{align*}
   \frac{1}{n} \sum_{k=1}^\nt X_{k-1,1} - \frac{b_1}{1-a_{1,1}} t
     & = \frac{1}{1-a_{1,1}} \cdot \frac{1}{n} \sum_{j=1}^{\nt-1} M_{j,1}
         - \frac{1}{1-a_{1,1}} \cdot \frac{1}{n} \sum_{j=1}^{\nt-1} a_{1,1}^{\nt-j} M_{j,1} \\
     &\phantom{=\;} + \frac{b_1}{1-a_{1,1}}\cdot \frac{1}{n}  \left( \nt -nt - 1  - \frac{a_{1,1}-a_{1,1}^\nt}{1-a_{1,1}} \right).
 \end{align*}
Here, since \ $a_{1,1}\in[0,1)$ \ and \ $\vert \nt - nt\vert\leq 1$, \ $t\in\RR_+$, \ we have for each \ $T\in\RR_{++}$,
 \[
    \sup_{t\in[0,T]} n^{-1}  \left( \nt -nt - 1  - \frac{a_{1,1}-a_{1,1}^\nt}{1-a_{1,1}} \right) \to 0 \qquad \text{as \ $n\to\infty$.}
 \]

Next we check that
 \begin{align}\label{help8}
 n^{-1} \sup_{t\in[0,T]}
           \Biggl|\sum_{j=1}^{\nt-1} M_{j,1}\Biggr|
   \stoch 0 \qquad \text{as \ $n\to\infty$}
 \end{align}
 for each \ $T\in\RR_{++}$, \ which is equivalent to
 \[
  n^{-2} \sup_{t\in[0,T]}
           \Biggl(\sum_{j=1}^{\nt-1} M_{j,1}\Biggr)^2
  = n^{-2} \max_{k\in\{1,\ldots,\nT-1\}} \Biggl(\sum_{j=1}^k M_{j,1}\Biggr)^2
   \stoch 0\qquad \text{as \ $n\to\infty$}
 \]
 for each \ $T\in\RR_{++}$.
\ Applying Doob's maximal inequality (see, e.g., Revuz and Yor
 \cite[Chapter II, Corollary (1.6)]{RevYor}) for the martingale
 \ $\sum_{j=1}^k M_{j,1}$, \ $k \in \NN$ \ (with the filtration \ $(\cF_k^\bX)_{k\in\NN}$), we obtain
 \begin{align*}
   \EE\left[\max_{k\in\{1,\ldots,\nT-1\}}  \left(\sum_{j=1}^k M_{j,1} \right)^2\right]
   & \leq 4 \EE\left[\left(\sum_{j=1}^{\nT-1} M_{j,1}\right)^2\right]
    = 4 \sum_{j=1}^{\nT-1} \EE(M_{j,1}^2)\\
   & = 4 \sum_{j=1}^{\nT-1} \OO(1) = \OO(n) ,
 \end{align*}
 where we used Lemma \ref{EEX5}, thus we obtain
 \[
   n^{-2} \max_{k\in\{1,\ldots,\nT -1\}} \left(\sum_{j=1}^k M_{j,1}\right)^2
      \mean 0 \qquad \text{as \ $n\to\infty$ \ for all \ $T\in\RR_{++}$,}
 \]
 yielding \eqref{help8}.

In a similar way as in the proof of \eqref{help7} (replacing \ $a_{2,2}$ \ by \ $a_{1,1}$), \ we prove that
 \[
   n^{-1} \sup_{t\in[0,T]}
           \Biggl|\sum_{j=1}^{\nt-1}
                   a_{1,1}^{\nt-j} M_{j,1}\Biggr|
          \leq n^{-1} \sup_{t\in[0,T]}
           \Biggl|\sum_{j=1}^{\nt-1}
                   a_{1,1}^{\nt-1-j} M_{j,1}\Biggr|
            = n^{-1} \max_{k\in\{1,\ldots,\nT-1\}} |\widetilde V_{k,1}|
   \stoch 0
 \]
 as \ $n \to \infty$ \ for all \ $T \in \RR_{++}$, \ where
 \[
   \widetilde V_{k,1}
     := \sum_{j=1}^k a_{1,1}^{k-j} M_{j,1} , \qquad k \in \NN ,
 \]
 and we used that \ $a_{1,1}\in[0,1)$.
\ Note that
 \begin{align}\label{help2_new}
   \widetilde V_{k,1} = a_{1,1} \widetilde V_{k-1,1} + M_{k,1} , \qquad k \in \NN ,
 \end{align}
 where \ $\widetilde V_{0,1} := 0$, \ hence \ $(\widetilde V_{k,1})_{k\in\ZZ_+}$ \ is a stable AR(1) process with heteroscedastic innovation \ $(M_{k,1})_{k\in\NN}$.
\ For all \ $\delta > 0$, \ by Markov's inequality, we have
 \begin{align*}
  \PP\Bigl( n^{-1} \max_{k\in\{1,\ldots,\nT-1\}} |\widetilde V_{k,1}| > \delta\Bigr)
   &= \PP\Bigl(\max_{k\in\{1,\ldots,\nT-1\}} \widetilde V_{k,1}^2 > \delta^2 n^2 \Bigr)
     \leq \sum_{k=1}^{\nT-1} \PP(\widetilde V_{k,1}^2 > \delta^2 n^2) \\
   &\leq \delta^{-2} n^{-2} \sum_{k=1}^{\nT-1} \EE(\widetilde V_{k,1}^2)
    = \delta^{-2} n^{-2} \sum_{k=1}^{\nT-1} \OO(1)
    = \OO(n^{-1})
 \end{align*}
 for \ $n\in\NN$, \ where we applied \ $\EE(\widetilde V_{k,1}^2) = \OO(1)$, \ $k\in\NN$ \ (see Lemma \ref{EEX5}).
Thus we obtain \eqref{Cond11_5}, and hence \eqref{Cond1_5}.

Now we turn to prove \eqref{Cond2_5}.
For each \ $s \in \RR_+$ \ and \ $n \in \NN$, \ we have \ $(\cX^{(n)}_{s,2})^+ = \cX^{(n)}_{s,2}$ \
 (due to the fact that \ $X_k$ \ is non-negative for each \ $k\in\ZZ_+$), \ and
 \begin{align*}
  &\int_0^t (\cX^{(n)}_{s,2})^+ \, \dd s
   = \int_0^t n^{-1} X_{\ns,2} \, \dd s
   = \sum_{k=0}^{\nt-1} \int_{k/n}^{(k+1)/n} n^{-1} X_{k,2} \, \dd s
     + \int_{\nt/n}^t n^{-1} X_{\nt,2} \, \dd s \\
  &= \frac{1}{n^2} \sum_{k=0}^{\nt-1} X_{k,2}
     + \frac{1}{n} \biggl(t - \frac{\nt}{n}\biggr) X_{\nt,2}
   = \frac{1}{n^2} \sum_{k=0}^{\nt-1} X_{k,2}
     + \frac{nt-\nt}{n^2} X_{\nt,2}
 \end{align*}
 for all \ $t \in \RR_+$ \ and \  $n \in \NN$.
\ Since
 \[
   U_{k,2} - \EE(U_{k,2}\mid \cF_{k-1}^\bX)
    = M_{k,2} + a_{2,1} X_{k-1,1} + b_2 - a_{2,1}X_{k-1,1} - b_2
    = M_{k,2},\qquad k\in\NN,
 \]
 by Lemma \ref{Moments}, we have for each \ $t \in \RR_+$ \ and \ $n \in \NN$,
 \[
   \frac{1}{n^2}
    \sum_{k=1}^\nt
     \var(U_{k,2} \mid \cF_{k-1}^\bX)
   = \frac{1}{n^2}
     \sum_{k=1}^\nt
      \EE(M_{k,2}^2 \mid \cF_{k-1}^\bX)
   = \frac{1}{n^2}
     \sum_{k=1}^\nt
      (v^{(0)}_{2,2} + v^{(1)}_{2,2} X_{k-1,1} + v^{(2)}_{2,2} X_{k-1,2}) .
 \]
Hence
 \begin{align*}
   &\frac{1}{n^2} \sum_{k=1}^\nt \var(U_{k,2} \mid \cF_{k-1}^\bX)
     -  \int_0^t (\cX^{(n)}_{s,2})^+ v_{2,2}^{(2)}\, \dd s \\
   &\qquad = \frac{\nt}{n^2} v_{2,2}^{(0)} +  v_{2,2}^{(1)}\frac{1}{n^2} \sum_{k=1}^\nt X_{k-1,1}
     - v_{2,2}^{(2)} \frac{nt-\nt}{n^2}  X_{\nt,2}, \qquad t\in\RR_+,\qquad n\in\NN.
 \end{align*}
By \eqref{Cond11_5}, we have for each \ $T\in\RR_{++}$,
 \begin{equation}\label{Cond12_5}
   n^{-2} \sup_{t\in[0,T]} \sum_{k=1}^\nt X_{k-1,1}
  \stoch 0 \qquad \text{as \ $n \to \infty$,}
 \end{equation}
 hence, in order to show \eqref{Cond2_5}, it suffices to prove
 \begin{equation}\label{Cond21_5}
  n^{-2} \sup_{t\in[0,T]} \vert (nt - \nt) X_{\nt,2} \vert
  \leq n^{-2} \sup_{t\in[0,T]} X_{\nt,2}
  \stoch 0 \quad \text{as \ $n \to \infty$ \ for each \ $T\in\RR_{++}$.}
 \end{equation}
For each \ $t \in \RR_+$ \ and \ $n \in \NN$, \
 \[
   n^{-2} X_{\nt,2}
   = \frac{1}{n^2} \sum_{k=1}^\nt U_{k,2}
   = \frac{1}{n^2} \sum_{k=1}^\nt (M_{k,2} + a_{2,1} X_{k-1,1} + b_2) ,
 \]
 hence
 \[
   n^{-2} X_{\nt,2}
    =  n^{-2} \vert X_{\nt,2} \vert
   \leq \frac{1}{n^2} \sum_{k=1}^\nt |M_{k,2}|
        + \frac{a_{2,1}}{n^2} \sum_{k=1}^\nt X_{k-1,1}
        + \frac{\nt}{n^2} b_2 .
 \]
Taking into account \eqref{Cond12_5}, in order to show \eqref{Cond21_5}, it suffices to show
 \begin{equation}\label{Cond22_5}
   \frac{1}{n^2} \sum_{k=1}^\nT |M_{k,2}| \stoch 0 \qquad \text{as \ $n \to \infty$ \ for each \ $T\in\RR_{++}$.}
 \end{equation}
We have for each \ $T\in\RR_{++}$,
 \[
   \EE\biggl(\frac{1}{n^2} \sum_{k=1}^\nT |M_{k,2}|\biggr)
   = \frac{1}{n^2} \sum_{k=1}^\nT \EE(|M_{k,2}|)
   \leq \frac{1}{n^2} \sum_{k=1}^\nT \sqrt{\EE(M_{k,2}^2)}
   \to 0 \qquad \text{as \ $n \to \infty$}
 \]
 since, by Lemma \ref{EEX5}, \ $\EE(M_{k,2}^2) = \OO(k)$, \ $k \in\NN$.
\ This yields \eqref{Cond22_5} and hence \eqref{Cond21_5}, implying \eqref{Cond2_5}, as desired.

Next, we check condition \eqref{Cond3_5}.
We show that for each \ $T\in\RR_{++}$ \ and \ $\theta\in\RR_{++}$,
 \begin{align}\label{help9}
   \frac{1}{n^2}
  \sum_{k=1}^\nT
   \EE( U_{k,2}^2 \bbone_{\{\vert U_{k,2}\vert >n\theta\}} \mid \cF_{k-1}^\bX)
   \mean 0 \qquad \text{as \ $n\to\infty$.}
 \end{align}
By Markov's inequality, for each \ $T\in\RR_{++}$ \ and $\theta\in\RR_{++}$, \ we have
 \begin{align*}
  &\EE\left(  \frac{1}{n^2}
  \sum_{k=1}^\nT
   \EE( U_{k,2}^2 \bbone_{\{\vert U_{k,2}\vert >n\theta\}} \mid \cF_{k-1}^\bX)
     \right)
   = \frac{1}{n^2}
      \sum_{k=1}^\nT \EE( U_{k,2}^2 \bbone_{\{ \vert U_{k,2} \vert >n\theta\}}) \\
  & \leq  \frac{1}{n^2}
          \sum_{k=1}^\nT \EE\left( \frac{U_{k,2}^4}{n^2\theta^2} \right)
   = \frac{1}{n^4\theta^2} \sum_{k=1}^\nT \EE( (M_{k,2} + a_{2,1} X_{k-1,1} + b_2)^4 ) \\
  & \leq \frac{3^3}{n^4 \theta^2}
        \sum_{k=1}^\nT\Big( \EE(M_{k,2}^4) + a_{2,1}^4 \EE(X_{k-1,1}^4) + b_2^4\Big).
 \end{align*}
Here, by Lemma \ref{EEX5},
 \begin{align*}
 \frac{1}{n^4} \sum_{k=1}^\nT \EE(M_{k,2}^4)
      = \frac{1}{n^4} \sum_{k=1}^\nT \OO(k^2)
      = \OO(n^{-1})
      \to 0 \qquad \text{as \ $n\to\infty$.}
 \end{align*}
Further, since the first coordinate process \ $(X_{k,1})_{k\in\ZZ_+}$ \ is a subcritical Galton-Watson process with immigration
 having offspring and immigration distributions with finite fourth moments, by Sz\H ucs \cite[Theorem 4]{Szu}
 (for an alternative proof, see also Kevei and Wiandt \cite{KevWia}), the unique stationary distribution
 \ $\mu_1$ \ of \ $(X_{k,1})_{k\in\ZZ_+}$ \ has a finite fourth moment as well.
Consequently, by Chung \cite[Part I, Chapter 15, Theorem 3]{Chu}, we have
 \[
   \frac{1}{n} \sum_{k=1}^n \EE(X_{k,1}^4) \to c_{4,\mu} \qquad \text{as \ $n\to\infty$,}
 \]
 where \ $c_{4,\mu}\in\RR_+$ \ denotes the fourth moment of \ $\mu_1$.
\ Hence
 \[
   \frac{1}{n^4} \sum_{k=1}^n \EE(X_{k,1}^4) \to 0 \qquad \text{as \ $n\to\infty$ \ for each \ $T\in\RR_{++}$.}
 \]
Putting parts together, we have \eqref{help9}, as desired.

Finally, using \eqref{Xdeco} with \ $a_{2,2}=1$, \ \eqref{Cov}, \eqref{help5} (which holds in case of (5) as well, since only the fact that
 \ $a_{1,2}=0$ \ was used for deriving it) and \ $\EE(\bM_k \mid \cF^{\bX}_{k-1})=\bzero$, \ $k\in\NN$ \
 (yielding \ $\cov(M_{j,1},M_{k,1})=0$ \ and \ $\cov(M_{j,1},M_{k,2})=0$ \ for \ $j\ne k$, \ $j,k\in\NN$), \ we have for all \ $t_1,t_2\in\RR_+$,
 \begin{align*}
  &\cov\big(X_{\ntk,1}, n_2^{-1} X_{\ntn,2}\big) \\
   &= \cov\!\left( \sum_{j=1}^{\ntk} a_{1,1}^{\ntk -j} (M_{j,1} + b_1),
           a_{2,1} n_2^{-1} \sum_{j=1}^{\ntn} \frac{1-a_{1,1}^{\ntn -j}}{1-a_{1,1}}(M_{j,1} + b_1)
           + n_2^{-1} \sum_{j=1}^{\ntn} (M_{j,2} + b_2) \right)
 \end{align*}
 \begin{align*}
   &= a_{2,1} n_2^{-1} \sum_{j=1}^{\ntk\wedge\ntn} a_{1,1}^{\ntk -j} \frac{1-a_{1,1}^{\ntn -j}}{1-a_{1,1}} \EE(M_{j,1}^2)
            + n_2^{-1} \sum_{j=1}^{\ntk\wedge\ntn} a_{1,1}^{\ntk -j} \cov(M_{j,1}, M_{j,2})\\
   &=  \frac{a_{2,1}}{1-a_{1,1}} n_2^{-1} \sum_{j=1}^{\ntk\wedge\ntn} a_{1,1}^{\ntk -j}(1-a_{1,1}^{\ntn -j})\big(v_{1,1}^{(0)} + \EE(X_{j-1,1}) v_{1,1}^{(1)}\big)\\
   &\phantom{=\,} + n_2^{-1} \sum_{j=1}^{\ntk\wedge\ntn} a_{1,1}^{\ntk -j} \Big( v_{1,2}^{(0)} + \EE(X_{j-1,1}) v_{1,2}^{(1)} + \EE(X_{j-1,2}) v_{1,2}^{(2)} \Big).
 \end{align*}
Since \ $a_{1,2}=0$, \ we have \ $\PP(\xi_{1,1,2,1}=0)=1$, \ yielding \ $v^{(2)}_{1,2}=0$, \ and using \eqref{help11},
 we have for all \ $t_1,t_2\in\RR_+$,
 \begin{align*}
   &\cov\big( X_{\ntk,1}, n_2^{-1} X_{\ntn,2}\big) \\
   &= \frac{a_{2,1}}{1-a_{1,1}} n_2^{-1} \sum_{j=1}^{\ntk\wedge\ntn} a_{1,1}^{\ntk -j} (1 - a_{1,1}^{\ntn -j}) \left(v_{1,1}^{(0)} + v_{1,1}^{(1)} b_1
              \frac{1 - a_{1,1}^{j-1}}{1-a_{1,1}} \right)\\
   &\phantom{=} + n_2^{-1} \sum_{j=1}^{\ntk\wedge\ntn} a_{1,1}^{\ntk -j}
                 \left(v_{1,2}^{(0)} + v_{1,2}^{(1)} b_1 \frac{1 - a_{1,1}^{j-1}}{1-a_{1,1}} \right)\\
   &\leq \frac{a_{2,1}}{1-a_{1,1}} \left(v_{1,1}^{(0)} + \frac{v_{1,1}^{(1)} b_1}{1-a_{1,1}}  \right)
            n_2^{-1} \sum_{j=1}^{\ntk\wedge\ntn} a_{1,1}^{\ntk -j}\\
   &\phantom{\leq}
          + \left( v_{1,2}^{(0)} + v_{1,2}^{(1)} \frac{b_1}{1-a_{1,1}} \right) n_2^{-1} \sum_{j=1}^{\ntk\wedge\ntn} a_{1,1}^{\ntk -j} \\
   &=\left( \frac{a_{2,1}}{1-a_{1,1}} +1 \right)  \left(v_{1,2}^{(0)} + v_{1,2}^{(1)} b_1 \frac{1 - a_{1,1}^{j-1}}{1-a_{1,1}} \right)
           n_2^{-1} a_{1,1}^{\ntk - \ntk\wedge\ntn} \frac{1 - a_{1,1}^{\ntk\wedge\ntn} }{1-a_{1,1}},
 \end{align*}
 which yields \eqref{help14_cov2_mod} (due to \ $a_{1,1}\in[0,1)$ \ and \ $X_{0,1}=0$).
\ Using again \ $a_{1,1}\in[0,1)$, \ we have
 \begin{align*}
    \sup_{t_1,t_2\in\RR_+} \sup_{n_1\in\NN}
            \Big\vert \cov(X_{\lfloor n_1t_1\rfloor, 1}, n_2^{-1}  X_{\lfloor n_2t_2\rfloor, 2} )\Big\vert
          \leq  \left( \frac{a_{2,1}}{1-a_{1,1}} +1 \right)  \left(v_{1,2}^{(0)} + v_{1,2}^{(1)} b_1 \frac{1 - a_{1,1}^{j-1}}{1-a_{1,1}} \right)
                \frac{1}{n_2}\to 0
 \end{align*}
 as \ $n_2\to\infty$, \ which yields \eqref{help14_cov2}, as desired.
\proofend

\vspace*{5mm}

\appendix

\vspace*{5mm}

\noindent{\bf\Large Appendices}

\section{Moments}
\label{app:moments}

In the proof of the results, we will use some formulae and estimates for the first,
 second and fourth order moments of the coordinates of the processes \ $(\bX_k)_{k\in\ZZ_+}$ \ and
 \ $(\bM_k)_{k\in\ZZ_+}$.

\begin{Lem}\label{Moments}
Let \ $(\bX_k)_{k\in\ZZ_+}$ \ be a $p$-type Galton-Watson process with immigration such that
 \ $\bX_0 = \bzero$ \ and the moment condition \eqref{help4} holds.
Then for all \ $k \in \NN$, \ we have \ $\EE(\bX_k \mid \cF_{k-1}^\bX) = \bA \bX_{k-1} + \bb$ \
 and
 \begin{gather}
  \EE(\bX_k)
   = \sum_{j=0}^{k-1} \bA^j \bb, \label{mean} \\
  \var\bigl(\bX_k \mid \cF_{k-1}^\bX\bigr)
    = \var\bigl(\bM_k \mid \cF_{k-1}^\bX\bigr)
	= \EE(\bM_k \bM_k^\top \mid \cF_{k-1}^\bX)
     = \bV^{(0)} + \sum_{i=1}^p X_{k-1,i} \bV^{(i)}, \label{condVar}\\
   \var\bigl(\bX_k\bigr)
      = \sum_{j=0}^{k-1} \bA^j \EE(\bM_{k-j} \bM_{k-j}^\top) (\bA^\top)^j, \label{Var}\\
  \EE(\bM_k \bM_k^\top)
  = \bV^{(0)} + \sum_{i=1}^p \EE(X_{k-1,i}) \bV^{(i)} , \label{Cov}
 \end{gather}
 where \ $\bV^{(0)}$ \ and \ $\bV^{(i)}$, \ $i=1,\ldots,p$, \ are given in Section \ref{section_multi_branching}.
\end{Lem}

Lemma \ref{Moments} is a special case of Lemma A.1 in Isp\'any and Pap \cite{IspPap2} for
 $p$-type Galton-Watson processes with immigration starting from \ $\bzero$.
\ For completeness, we note that Lemma A.1 in Isp\'any and Pap \cite{IspPap2} is stated only for critical
 $p$-type Galton-Watson processes with immigration, but its proof readily shows that it holds not only in the critical case.

\begin{Lem}\label{EEX1}
Let \ $(\bX_k)_{k\in\ZZ_+}$ \ be a critical decomposable 2-type Galton-Watson process with immigration such that \ $\bX_0=\bzero$, \
 the moment conditions  \ $\EE(\|\bxi_i\|^4) < \infty$, \ $i=1,2$, \ and \ $\EE (\|\bvare\|^4) < \infty$ \ hold and its offspring mean matrix \ $\bA$ \ satisfies \textup{(1)}
 of \eqref{tablazat_esetek}.
Then we have
 \begin{align*}
 & \EE(X_{k,1}) = \OO(k) , \qquad
  \EE(X_{k,2}) = \OO(k) , \qquad
  \EE(|M_{k,1}|) = \OO(k^{1/2}), \qquad
  \EE(|M_{k,2}|) = \OO(k^{1/2}) ,\\
 & \EE(M_{k,1}^2) = \OO(k), \qquad
   \EE(M_{k,2}^2) = \OO(k), \qquad
   \EE(X_{k,1}^2) = \OO(k^2), \qquad
   \EE(X_{k,2}^2) = \OO(k^2), \\
  & \EE(M_{k,1}^4) = \OO(k^2), \qquad
   \EE(M_{k,2}^4) = \OO(k^2)
 \end{align*}
 for \ $k\in\NN$.
\end{Lem}

\noindent \textbf{Proof.}
By \eqref{mean} and \eqref{m^ell}, we obtain
 \[
   \begin{bmatrix} \EE(X_{k,1}) \\ \EE(X_{k,2}) \end{bmatrix}
   = \sum_{j=0}^{k-1}
      \begin{bmatrix}
       1 & 0 \\
       0 & 1
      \end{bmatrix}
      \bb
   = \begin{bmatrix}
      b_1 k \\
      b_2 k
     \end{bmatrix} , \qquad k\in\ZZ_+,
 \]
 and we conclude the first two statements.

By \eqref{Cov},
 \begin{align*}
  \EE(M_{k,1}^2)
  &= v^{(0)}_{1,1}
     + \EE(X_{k-1,1}) v^{(1)}_{1,1} , \qquad k\in\NN,\\
  \EE(M_{k,2}^2)
  &= v^{(0)}_{2,2}
     + \EE(X_{k-1,2}) v^{(2)}_{2,2}, \qquad k\in\NN,
 \end{align*}
 since \ $a_{1,2} = a_{2,1} = 0$ \ implies \ $\xi_{1,1,2,1}\ase 0$ \ and \ $\xi_{1,1,1,2}\ase 0$,
 \ yielding \ $v^{(2)}_{1,1} = v^{(1)}_{2,2} = 0$.
\ This together with the first two statements yield \ $\EE(M_{k,i}^2) = \OO(k)$ \ for \ $k\in\NN$ \ and \ $i=1,2$.
\ Consequently, using the inequalities \ $\EE(\vert M_{k,i}\vert)\leq \sqrt{\EE(M_{k,i}^2)}$, \ $i=1,2$, \
 we also have \ $\EE(\vert M_{k,i}\vert) = \OO(k^{1/2})$ \ for \ $k\in\NN$ \ and \ $i=1,2$.

Further, \ $\EE(X_{k,i}^2) = \var(X_{k,i}) + (\EE(X_{k,i}))^2$, \ $i=1,2$, \ and, by \eqref{Var} and \eqref{Cov}, we have
 \begin{align*}
   \var(\bX_k) & = \sum_{j=0}^{k-1} \bA^j \EE(\bM_{k-j} \bM_{k-j}^\top) (\bA^\top)^j
                 = \sum_{j=0}^{k-1} \bI_2^j \EE(\bM_{k-j} \bM_{k-j}^\top) \bI_2^j \\
                & = \sum_{j=0}^{k-1} \left( \bV^{(0)} + \sum_{i=1}^2 \EE(X_{k-j-1,i}) \bV^{(i)} \right),
                \qquad k\in\ZZ_+.
 \end{align*}
By the first two statements,
  \begin{align*}
  \Vert \var(\bX_k)\Vert
   & \leq \sum_{j=0}^{k-1} \left( \Vert \bV^{(0)} \Vert + \sum_{i=1}^2 \EE(X_{k-j-1,i}) \Vert \bV^{(i)} \Vert \right)\\
   &  = \sum_{j=0}^{k-1} \left( \Vert \bV^{(0)} \Vert + \sum_{i=1}^2 \OO(k-j-1) \Vert \bV^{(i)} \Vert \right)
     = \OO(k^2),\qquad k\in\NN,
 \end{align*}
 where \ $\Vert \bB\Vert$ \ denotes the operator norm of a matrix \ $\bB\in\RR^{2\times 2}$ \ defined by
  \ $\Vert \bB\Vert:=\sup_{\Vert \bx\Vert=1,\, \bx\in\RR^2} \Vert \bB\bx\Vert$.
\ This together with the first two statements yield \ $\EE(X_{k,i}^2) = \OO(k^2)$ \ for \ $k\in\NN$ \ and \ $i=1,2$.
\ Finally, the relations \ $\EE(M_{k,i}^4) = \OO(k^2)$ \ for \ $k\in\NN$, \ $i=1,2$, \ follow in the same way as in the proof of
 Lemma A.2 in Isp\'any and Pap \cite{IspPap2} (at this part the authors do not use that the critical multi-type Galton-Watson
 process with immigration that they consider is primitive, they only use the fact that the second moments of the
 coordinates of the branching process in question at \ $k$ \ is of \ $\OO(k^2)$, \ $k\in\NN$).
\proofend

\begin{Lem}\label{EEX2}
Let \ $(\bX_k)_{k\in\ZZ_+}$ \ be a critical decomposable 2-type Galton-Watson process with immigration such that \ $\bX_0=\bzero$, \
 the moment condition \eqref{help4} holds and its offspring mean matrix \ $\bA$ \ satisfies \textup{(2)}  of \eqref{tablazat_esetek}.
Then we have
 \begin{align*}
 & \EE(X_{k,1}) = \OO(k) , \qquad
  \EE(X_{k,2}) = \OO(k^2) , \qquad
  \EE(|M_{k,1}|) = \OO(k^{1/2}), \qquad
  \EE(|M_{k,2}|) = \OO(k) ,\\
 & \EE(M_{k,1}^2) = \OO(k), \qquad \EE(M_{k,2}^2) = \OO(k^2)
 \end{align*}
 for \ $k\in\NN$.
\end{Lem}

\noindent
\textbf{Proof.}
By \eqref{mean} and \eqref{m^ell}, we obtain
 \[
   \begin{bmatrix} \EE(X_{k,1}) \\ \EE(X_{k,2}) \end{bmatrix}
   = \sum_{j=0}^{k-1}
      \begin{bmatrix}
       1 & 0 \\
       a_{2,1} j & 1
      \end{bmatrix}
      \bb
   = \begin{bmatrix}
      b_1 k \\
      \frac{1}{2} a_{2,1} b_1 k(k-1)
      + b_2 k
     \end{bmatrix} , \qquad k\in\ZZ_+,
 \]
 and we conclude the first two statements.

By \eqref{Cov},
 \begin{align}\label{help5}
  \begin{split}
  \EE(M_{k,1}^2)
  &= v^{(0)}_{1,1}
     + \EE(X_{k-1,1}) v^{(1)}_{1,1} , \qquad k\in\NN,\\
  \EE(M_{k,2}^2)
  &= v^{(0)}_{2,2}
     + \EE(X_{k-1,1}) v^{(1)}_{2,2}
     + \EE(X_{k-1,2}) v^{(2)}_{2,2}, \qquad k\in\NN,
  \end{split}
 \end{align}
 since \ $a_{1,2} = 0$ \ implies \ $\xi_{1,1,2,1}\ase 0$, \ yielding \ $v^{(2)}_{1,1} = 0$.
\ Note that for deriving \eqref{help5} we did not use that \ $\bA$ \ satisfies (2), we only used \ $a_{1,2}=0$,
 \ so \eqref{help5} holds if \ $\bA$ \ has the form (1), (3), (4) or (5) as well.
Using \eqref{help5} together with the first two statements we have the last two statements.
Finally, using the inequalities \ $\EE(\vert M_{k,i}\vert)\leq \sqrt{\EE(M_{k,i}^2)}$, \ $i=1,2$, \
 and the last two statements, we have the third and fourth statements.
\proofend

We note that the statements \ $\EE(X_{k,1}) = \OO(k)$, \ $k\in\NN$, \ and \ $\EE(X_{k,2}) = \OO(k^2)$, \ $k\in\NN$, \ in Lemma \ref{EEX2}
 are in accordance with the corresponding ones in Theorem 3 in Foster and Ney \cite{FosNey2}.

\begin{Lem}\label{EEX3}
Let \ $(\bX_k)_{k\in\ZZ_+}$ \ be a critical decomposable 2-type Galton-Watson process with immigration such that \ $\bX_0=\bzero$, \
 the moment condition \eqref{help4} holds and its offspring mean matrix \ $\bA$ \ satisfies \textup{(3)}  of \eqref{tablazat_esetek}.
Then we have
 \begin{align*}
 & \EE(X_{k,1}) = \OO(k) , \qquad
  \EE(X_{k,2}) = \OO(1) , \qquad
  \EE(|M_{k,1}|) = \OO(k^{1/2}), \qquad
  \EE(|M_{k,2}|) = \OO(1) ,\\
 & \EE(M_{k,1}^2) = \OO(k), \qquad \EE(M_{k,2}^2) = \OO(1)
 \end{align*}
 for \ $k\in\NN$.
\end{Lem}

\noindent
\textbf{Proof.}
By \eqref{mean} and \eqref{m^ell}, we obtain
 \[
   \begin{bmatrix} \EE(X_{k,1}) \\ \EE(X_{k,2}) \end{bmatrix}
   = \sum_{j=0}^{k-1}
      \begin{bmatrix}
       1 & 0 \\
       0 & a_{2,2}^j
      \end{bmatrix}
      \bb
   = \begin{bmatrix}
      b_1 k \\
      \frac{1 - a_{2,2}^k}{1-a_{2,2}} b_2
     \end{bmatrix} , \qquad k\in\ZZ_+,
 \]
 and we conclude the first two statements.

Using \eqref{help5} (which holds in case of \textup{(3)} as well, since only the fact that \ $a_{1,2}=0$ \ was used for deriving it),
 the first two statements, and that \ $a_{2,1}=0$ \ implies \ $v^{(1)}_{2,2}=0$, \ we have the last two statements.
Finally, using the inequalities \ $\EE(\vert M_{k,i}\vert)\leq \sqrt{\EE(M_{k,i}^2)}$, \ $i=1,2$, \
 and the last two statements, we have the third and fourth statements.
\proofend

\begin{Lem}\label{EEX4}
Let \ $(\bX_k)_{k\in\ZZ_+}$ \ be a critical decomposable 2-type Galton-Watson process with immigration such that \ $\bX_0=\bzero$, \
 the moment conditions \ $\EE(\|\bxi_i\|^4) < \infty$, \ $i=1,2$, \ and \ $\EE (\|\bvare\|^4) < \infty$ \
 hold and its offspring mean matrix \ $\bA$ \ satisfies \textup{(4)}  of \eqref{tablazat_esetek}.
Then we have
 \begin{gather*}
  \EE(X_{k,1}) = \OO(k) , \qquad
  \EE(X_{k,2}) = \OO(k) , \qquad
  \EE(|M_{k,1}|) = \OO(k^{1/2}), \qquad
  \EE(|M_{k,2}|) = \OO(k^{1/2}) , \\
  \EE(M_{k,1}^2) = \OO(k), \qquad
  \EE(M_{k,2}^2) = \OO(k) , \qquad
  \EE(X_{k,1}^2) = \OO(k^2) , \qquad
  \EE(X_{k,2}^2) = \OO(k^2) , \\
  \EE(M_{k,1}^4) = \OO(k^2) , \qquad
  \EE(M_{k,2}^4) = \OO(k^2) , \qquad
  \EE(V_{k,1}^4) = \OO(k^2) , \qquad
  \EE(V_{k,2}^4) = \OO(k^2)
 \end{gather*}
 for \ $k\in\NN$, \ where
 \[
   V_{k,i} := \sum_{j=1}^k a_{2,2}^{k-j} M_{j,i} , \qquad k \in \NN , \quad i \in \{1, 2\} .
 \]
\end{Lem}

\noindent
\textbf{Proof.}
By \eqref{mean} and \eqref{m^ell}, we obtain
 \[
   \begin{bmatrix} \EE(X_{k,1}) \\ \EE(X_{k,2}) \end{bmatrix}
   = \sum_{j=0}^{k-1}
      \begin{bmatrix}
       1 & 0 \\
       a_{2,1} \frac{1-a_{2,2}^j}{1-a_{2,2}}
        & a_{2,2}^j
      \end{bmatrix}
      \bb
   = \begin{bmatrix}
      k b_1 \\
      \frac{a_{2,1} b_1}{1-a_{2,2}}
      \left(k - \frac{1-a_{2,2}^k}{1-a_{2,2}}\right)
      + \frac{1-a_{2,2}^k}{1-a_{2,2}} b_2
     \end{bmatrix} , \qquad k\in\ZZ_+,
 \]
 and we conclude the first two statements.

The next four statements follow from \eqref{help5} (which holds in case of (4) as well as it was explained earlier)
 using the inequalities \ $\EE(\vert M_{k,i}\vert) \leq \sqrt{\EE(M_{k,i}^2)}$, \ $i=1,2$.

Further, \ $\EE(X_{k,i}^2) = \var(X_{k,i}) + (\EE(X_{k,i}))^2$, \ $i=1,2$, \ and, by \eqref{Var} and \eqref{Cov},
 \begin{align*}
   \var(\bX_k) & = \sum_{j=0}^{k-1} \bA^j \EE(\bM_{k-j} \bM_{k-j}^\top) (\bA^\top)^j
                = \sum_{j=0}^{k-1} \bA^j \left( \bV^{(0)} + \sum_{i=1}^2 \EE(X_{k-j-1,i}) \bV^{(i)} \right)(\bA^\top)^j
 \end{align*}
 for \ $k\in\ZZ_+$.
\ Hence, using the first two statements, we have
 \begin{align*}
  \Vert \var(\bX_k)\Vert
   & \leq \sum_{j=0}^{k-1} \Vert \bA^j \Vert^2  \left( \Vert \bV^{(0)} \Vert + \sum_{i=1}^2 \EE(X_{k-j-1,i}) \Vert \bV^{(i)} \Vert \right)
     = \sum_{j=0}^{k-1} \Vert \bA^j \Vert^2 \OO(k-j-1) \\
   & = \left(\sum_{j=0}^{k-1} \Vert \bA^j \Vert^2 \right)\OO(k),\qquad k\in\NN.
 \end{align*}
Using the continuity of the operator norm function \ $\Vert \cdot \Vert$, \ we have
 \begin{align*}
   \lim_{j\to\infty}  \Vert \bA^j \Vert
     = \lim_{j\to\infty}
         \left\Vert  \begin{bmatrix}
                       1 & 0 \\
                       a_{2,1}\frac{1-a_{2,2}^j}{1- a_{2,2}} & a_{2,2}^j \\
                     \end{bmatrix}
           \right\Vert
     = \left\Vert  \begin{bmatrix}
                       1 & 0 \\
                       \frac{a_{2,1}}{1-a_{2,2}} & 0 \\
                     \end{bmatrix}
           \right\Vert
     <\infty,
 \end{align*}
 yielding that \ $c_\bA := \sup_{j\in\NN}\Vert \bA^j \Vert <\infty$.
\ Hence \ $\Vert \var(\bX_k)\Vert \leq c_\bA^2 \left(\sum_{j=0}^{k-1} 1\right)\OO(k) = \OO(k^2)$, \ $k\in\NN$.
\ This together with the first two statements yield \ $\EE(X_{k,i}^2) = \OO(k^2)$ \ for \ $k\in\NN$ \ and \ $i=1,2$.
\ The relations \ $\EE(M_{k,i}^4) = \OO(k^2)$ \ for \ $k\in\NN$, \ $i=1,2$, \ follow in the same way as in the proof of
 Lemma A.2 in Isp\'any and Pap \cite{IspPap2} (at this part the authors do not use that the critical multi-type Galton-Watson
 process with immigration that they consider is primitive, they only use the fact that the second moments of the
 coordinates of the branching process in question at \ $k$ \ is of \ $\OO(k^2)$, \ $k\in\NN$).

Finally, we prove \ $\EE(V_{k,i}^4) = \OO(k^2)$ \ for \ $k\in\NN$, \ $i=1,2$, \ using induction in \ $k$.
Since \ $\EE(M_{k,i}^4) = \OO(k^2)$ \ for \ $k\in\NN$ \ and \ $i=1,2$, \ there exists
  \ $\tC\in\RR_{++}$ \ such that \ $\EE(M_{k,i}^4)\leq \tC k^2$ \ for each \ $k\in\NN$ \ and \ $i=1,2$.
\ For each \ $k\in\NN$, \ let
 \[
  C_k:=\left(\sum_{j=0}^{k-1}a_{2,2}^j\right)^4\tC=\left(\frac{1-a_{2,2}^k}{1-a_{2,2}}\right)^4\tC.
 \]
Then \ $C_k\leq\frac{\tC}{(1-a_{2,2})^4}<\infty$, \ $k\in\NN$ \ (due to \ $a_{2,2}\in[0,1)$),
 \ and \ $C_{k+1}^{1/4} = a_{2,2}C_{k}^{1/4}+\tC^{1/4}$, \ $k\in\NN$, \ since
 \[
   a_{2,2}C_{k}^{\frac{1}{4}}+\tC^{\frac{1}{4}}
     = a_{2,2} \frac{1-a_{2,2}^k}{1-a_{2,2}} \tC^{\frac{1}{4}} + \tC^{\frac{1}{4}}
     = \frac{1-a_{2,2}^{k+1}}{1-a_{2,2}} \tC^{\frac{1}{4}}
     = C_{k+1}^{\frac{1}{4}}.
 \]
Since \ $V_{1,i}=M_{1,i}$, \ $i=1,2$, \ we get \ $\EE(V_{1,i}^4) = \EE(M_{1,i}^4) \leq \tC = C_1$, \ $i=1,2$.
\ Now assume that for some \ $k_0\in\NN$, \ the following inequalities hold
 \begin{align}\label{help18}
   \EE(V_{\ell,i}^4)\leq C_\ell\, \ell^2, \qquad \ell=1,\ldots,k_0.
 \end{align}
By the decompositions \ $V_{k,i} = a_{2,2} V_{k-1,i} + M_{k,i}$, \ $k\in\NN$, $i=1,2$, \
 where \ $V_{0,i}:=0$, $i=1,2$ \ (see also \eqref{help2} and \eqref{help3}),
 and the triangular inequality for the \ $L_4$-norm, we have that
 \[
   (\EE(V_{k_0+1,i}^4))^{\frac{1}{4}}  \leq a_{2,2} (\EE(V_{k_0,i}^4))^{\frac{1}{4}} + (\EE(M_{k_0,i}^4))^{\frac{1}{4}}, \qquad i=1,2.
 \]
Consequently, using also the induction hypothesis \eqref{help18}, we get that
 \begin{align*}
  (\EE(V_{k_0+1,i}^4))^{\frac{1}{4}}\leq a_{2,2}C_{k_0}^{\frac{1}{4}}k_0^{\frac{1}{2}}+\tC^{\frac{1}{4}}k_0^{\frac{1}{2}}
      =C_{k_0+1}^{\frac{1}{4}}k_0^{\frac{1}{2}}
      \leq C_{k_0+1}^{\frac{1}{4}}(k_0+1)^{\frac{1}{2}},\qquad i=1,2,
 \end{align*}
 yielding that \ $\EE(V_{k_0+1,i}^4) \leq C_{k_0+1}(k_0+1)^2$, \ $i=1,2$.
\ Hence  \ $\EE(V_{k,i}^4) \leq C_k k^2\leq \frac{\tC}{(1-a_{2,2})^4} k^2$ \ for each \ $k\in\NN$ \ and \ $i=1,2$,
 \ which implies that \ $\EE(V_{k,i}^4) = \OO(k^2)$ \ for each \ $k\in\NN$, \ $i=1,2$.
\proofend

\begin{Lem}\label{EEX5}
Let \ $(\bX_k)_{k\in\ZZ_+}$ \ be a critical decomposable 2-type Galton-Watson process with immigration such that \ $\bX_0=\bzero$, \
 the moment conditions \ $\EE(\|\bxi_i\|^4) < \infty$, \ $i=1,2$, \ and \ $\EE (\|\bvare\|^4) < \infty$ \ hold
 and its offspring mean matrix \ $\bA$ \ satisfies \textup{(5)} of \eqref{tablazat_esetek}.
Then we have
 \begin{align*}
 & \EE(X_{k,1}) = \OO(1) , \qquad
  \EE(X_{k,2}) = \OO(k) , \qquad
  \EE(|M_{k,1}|) = \OO(1), \qquad
  \EE(|M_{k,2}|) = \OO(k^{1/2}) ,\\
 & \EE(M_{k,1}^2) = \OO(1), \qquad
   \EE(M_{k,2}^2) = \OO(k), \qquad
   \EE(X_{k,1}^2) = \OO(1), \qquad
   \EE(X_{k,2}^2) = \OO(k^2), \\
 & \EE(M_{k,1}^4) = \OO(1), \qquad
   \EE(M_{k,2}^4) = \OO(k^2), \qquad
   \EE(\widetilde V_{k,1}^2) = \OO(1)
 \end{align*}
 for \ $k\in\NN$, where
 \[
   \widetilde V_{k,1}:=\sum_{j=1}^k a_{1,1}^{k-j} M_{j,1},\qquad k\in\NN.
 \]
\end{Lem}

\noindent
\textbf{Proof.}
By \eqref{mean} and \eqref{m^ell}, we obtain
 \begin{align}\label{help11}
   \begin{bmatrix} \EE(X_{k,1}) \\ \EE(X_{k,2}) \end{bmatrix}
   = \sum_{j=0}^{k-1}
      \begin{bmatrix}
       a_{1,1}^j & 0 \\
       a_{2,1}\frac{1-a_{1,1}^j}{1-a_{1,1}}  & 1
      \end{bmatrix}
      \bb
   = \begin{bmatrix}
      b_1 \frac{1-a_{1,1}^k}{1-a_{1,1}}  \\[1mm]
      \frac{a_{2,1} b_1}{1-a_{1,1}} \left(k- \frac{1-a_{1,1}^k}{1-a_{1,1}} \right) + b_2 k
     \end{bmatrix} , \qquad k\in\ZZ_+,
 \end{align}
 and we conclude the first two statements.

The next four statements follow from \eqref{help5} (which holds in case of (5) as well as it was explained earlier)
 using the inequalities \ $\EE(\vert M_{k,i}\vert) \leq \sqrt{\EE(M_{k,i}^2)}$, \ $i=1,2$.

Further, \ $\EE(X_{k,1}^2) = \var(X_{k,1}) + (\EE(X_{k,1}))^2$, \ and using that \ $(X_{k,1})_{k\in\ZZ_+}$ \ is a single-type
 Galton-Watson process with immigration starting from \ $0$ \ (explained at the beginning of Section \ref{Section_conv_results}),
 by \eqref{Var} and \eqref{help5}, we have
 \begin{align*}
   \var(X_{k,1}) & = \sum_{j=0}^{k-1} a_{1,1}^{2j} \EE(M_{k-j,1}^2)
                   = \sum_{j=0}^{k-1} a_{1,1}^{2j} \left( v_{1,1}^{(0)} + \EE(X_{k-j-1,1}) v_{1,1}^{(1)} \right) \\
                 & = \sum_{j=0}^{k-1} a_{1,1}^{2j} \left( v_{1,1}^{(0)} + b_1 \frac{1-a_{1,1}^{k-j-1}}{1-a_{1,1}} v_{1,1}^{(1)} \right)
                  \leq \left( v_{1,1}^{(0)} + \frac{b_1}{1-a_{1,1}} v_{1,1}^{(1)} \right)
                       \sum_{j=0}^{k-1} a_{1,1}^{2j}\\
                 & \leq \left( v_{1,1}^{(0)} + \frac{b_1}{1-a_{1,1}} v_{1,1}^{(1)} \right) \frac{1}{1-a_{1,1}^2}
                    = \OO(1), \qquad  k\in\NN.
 \end{align*}
This together with \ $\EE(X_{k,1}) = \OO(1)$ \ for \ $k\in\NN$ \ yield \ $\EE(X_{k,1}^2) = \OO(1)$ \ for \ $k\in\NN$.
\ The relation \ $\EE(M_{k,1}^4) = \OO(1)$ \ for \ $k\in\NN$ \ follows in the same way as in the proof of
 Lemma A.2 in Isp\'any and Pap \cite{IspPap2} taking into account the fact that \ $\EE(X_{k,1}^2) = \OO(1)$ \ for \ $k\in\NN$.

Next, we check that \ $\EE(X_{k,2}^2) = \OO(k^2)$ \ for \ $k\in\NN$.
\ By \eqref{Var}, \eqref{Cov} and the first two statements, we have
 \begin{align*}
  \Vert \var(\bX_k)\Vert
   & = \left\Vert \sum_{j=0}^{k-1} \bA^j \left(\bV^{(0)} + \sum_{i=1}^2 \EE(X_{k-j-1,i})\bV^{(i)} \right) (\bA^\top)^j \right\Vert \\
   & \leq \sum_{j=0}^{k-1} \Vert \bA^j \Vert^2  \left( \Vert \bV^{(0)} \Vert + \sum_{i=1}^2 \EE(X_{k-j-1,i}) \Vert \bV^{(i)} \Vert \right)\\
   & = \sum_{j=0}^{k-1} \Vert \bA^j \Vert^2 \OO(k-j-1)
      = \left(\sum_{j=0}^{k-1} \Vert \bA^j \Vert^2 \right)\OO(k),\qquad k\in\NN.
 \end{align*}
Using the continuity of the norm function \ $\Vert \cdot \Vert$, \ we have
 \begin{align*}
   \lim_{j\to\infty}  \Vert \bA^j \Vert
     = \lim_{j\to\infty}
         \left\Vert  \begin{bmatrix}
                       a_{1,1}^j & 0 \\
                       a_{2,1}\frac{1-a_{1,1}^j}{1- a_{1,1}} & 1 \\
                     \end{bmatrix}
           \right\Vert
     = \left\Vert  \begin{pmatrix}
                       0 & 0 \\
                       \frac{a_{2,1}}{1-a_{1,1}} & 1 \\
                     \end{pmatrix}
           \right\Vert
     <\infty,
 \end{align*}
 yielding that \ $c_\bA := \sup_{j\in\NN}\Vert \bA^j \Vert <\infty$.
\ Hence \ $\Vert \var(\bX_k)\Vert \leq c_\bA^2 \left(\sum_{j=0}^{k-1} 1\right)\OO(k) = \OO(k^2)$, \ $k\in\NN$.
 This together with \ $\EE(X_{k,2}) = \OO(k)$ \ for \ $k\in\NN$ \ yields \ $\EE(X_{k,2}^2) = \OO(k^2)$ \ for \ $k\in\NN$, \ as desired.
Note that the above estimation for \ $\Vert \var(\bX_k)\Vert$ \ also yields the crude estimation
 \ $\EE(X_{k,1}^2) = \OO(k^2)$ \ for \ $k\in\NN$, \ but, as we already showed, \ $\EE(X_{k,1}^2) = \OO(1)$ \ for \ $k\in\NN$ \
 holds.

The relation \ $\EE(M_{k,2}^4) = \OO(k^2)$ \ for \ $k\in\NN$ \ follows in the same way as in the proof of
 Lemma A.2 in Isp\'any and Pap \cite{IspPap2} taking into account the fact that \ $\EE(X_{k,2}^2) = \OO(k^2)$ \ for \ $k\in\NN$.

Finally, we prove \ $\EE(\widetilde V_{k,1}^2) = \OO(1)$ \ for \ $k\in\NN$.
Note that, for each \ $k\in\NN$, \ we have that \ $\EE(M_{j,1} M_{\ell,1})=0$ \ for \ $j\ne \ell$, \ $j,\ell = 1,\ldots,k$.
\ Consequently, using also that \ $\EE(M_{k,1}^2) = \OO(1)$ \ for \ $k\in\NN$ \ (which was proved before), we get that
 \begin{align*}
  \EE(\widetilde V_{k,1}^2)&=\EE\left(\left(\sum_{j=1}^ka_{1,1}^{k-j}M_{j,1}\right)^2\right)
                           =\sum_{j=1}^ka_{1,1}^{2(k-j)}\EE(M_{j,1}^2)
                           =\left(\sum_{j=1}^ka_{1,1}^{2(k-j)}\right)\OO(1)=\OO(1).
 \end{align*}
 for \ $k\in\NN$, \ since
 \ $\Big\vert\sum_{j=1}^k a_{1,1}^{2(k-j)}\Big\vert \leq \sum_{j=1}^\infty a_{1,1}^{2j} = \frac{a_{1,1}^2}{1-a_{1,1}^2}<\infty$, \ $k\in\NN$,
 \ due to \ $a_{1,1}\in[0,1)$.
\proofend

\section{Asymptotic behaviour of a single-type Galton-Watson process with immigration in the subcritical and critical cases}
\label{app:singletype_GWI}

Let \ $(X_k)_{k\in\ZZ_+}$ \ be a single-type Galton-Watson process with immigration, i.e.,
 \ $X_k = \sum_{j=1}^{X_{k-1}} \xi_{k,j} + \vare_k$, \ $k\in\NN$, \ where \ $\{X_0, \xi_{k,j}, \vare_k : k,j\in\NN\}$ \
 are supposed to be independent, \ $\{\xi_{k,j} : k,j\in\NN\}$ \ and \ $\{\vare_k : k\in\NN\}$ \ are supposed to consist of
 identically distributed \ $\ZZ_+$-valued random variables.
Let \ $\xi$ \ and \ $\vare$ \ be random variables such that \ $\xi \distre \xi_{1,1}$ \ and \ $\vare \distre \vare_1$.
\ If \ $a= \EE(\xi) \in [0, 1)$ \ and \ $\sum_{\ell=1}^\infty \log(\ell) \PP(\vare = \ell) < \infty$, \
 then the Markov chain \ $(X_k)_{k\in\ZZ_+}$ \ admits a unique stationary distribution \ $\mu$ \
 with a generator function
 \begin{align}\label{stac_generator}
  \prod_{j=0}^\infty H(G_{(j)}(z)) , \qquad z \in D:=\{ z \in \CCC : |z| \leq 1\} ,
 \end{align}
 where
 \[
   G(z) := \EE(z^{\xi}) ,\qquad  z\in D, \qquad\text{and} \qquad
   H(z) := \EE(z^{\vare}),
   \qquad z\in D,
 \]
 are the generator functions of \ $\xi$ \ and \ $\vare$, \ respectively,
 \ $G_{(0)}(z) := z$, \ $G_{(1)}(z) := G(z)$, \ and \ $G_{(k+1)}(z) := G_{(k)}(G(z))$, \ $z \in D$,
 \ $k \in \NN$, \ see, e.g., Quine \cite{Q}.
Note also that if \ $a \in [0, 1)$ \ and \ $\PP(\vare = 0) = 1$, \ then \ $\sum_{\ell=1}^\infty \log(\ell) \PP(\vare = \ell) = 0$ \ and
  \ $\mu$ \ is the Dirac measure \ $\delta_0$ \ concentrated at the point \ $0$.
\ In fact, \ $\mu = \delta_0$ \ if and only if \ $\PP(\vare = 0) = 1$.
\ Moreover, if \ $a = 0$ \ (which is equivalent to \ $\PP(\xi = 0) = 1$), \ then \ $\mu$ \ is the distribution of \ $\vare$.

The next result is about the asymptotic behaviour of single-type subcritical Galton-Watson processes with immigration
 satisfying first order moment conditions, which may be known, but we could not address any reference for it, so we provide a proof as well.

\begin{Lem}\label{lemma:subcritical}
Let \ $(X_k)_{k\in\ZZ_+}$ \ be a single-type Galton-Watson process with immigration
 such that \ $a =\EE(\xi) \in [0, 1)$, \ $\EE(\vare) < \infty$ \ and \ $\EE(X_0)<\infty$.
\ Then
 \[
   (X_\nt)_{t\in\RR_{++}} \distrf (\cX_t)_{t\in\RR_{++}} \qquad
  \text{as \ $n \to \infty$,}
 \]
 where \ $(\cX_t)_{t\in\RR_{++}}$ \ is an i.i.d. process such that for each \ $t \in \RR_{++}$,
 \ the distribution of \ $\cX_t$ \ is \ $\mu$ \ having generator function given in \eqref{stac_generator}.
\end{Lem}

\noindent
\textbf{Proof.}
Under the assumptions, the Markov chain \ $(X_k)_{k\in\ZZ_+}$ \ admits a unique stationary distribution \ $\mu$ \
 with expectation \ $\frac{b}{1-a}$ \ (where \ $b=\EE(\vare)$) \ and a generator function given in \eqref{stac_generator}.
Let \ $(Y_k)_{k\in\ZZ_+}$ \ be a Galton-Watson processes with immigration with the same offspring and immigration variables as
 \ $(X_k)_{k\in\ZZ_+}$, \ but let the distribution of \ $Y_0$ \ be \ $\mu$.
\ Hence the Markov chain \ $(Y_k)_{k\in\ZZ_+}$ \ is strongly stationary.
By induction with respect to \ $m \in \NN$, \ first we check that for each \ $t_1, \ldots, t_m \in \RR_{++}$ \ with \ $t_1 < \ldots < t_m$, \ we have
 \begin{equation}\label{stationary_convergence}
  (Y_{\lfloor nt_1\rfloor}, \ldots, Y_{\lfloor nt_m\rfloor})
  \distr (\cX_{t_1}, \ldots, \cX_{t_m}) \qquad \text{as \ $n \to \infty$.}
 \end{equation}
For \ $m = 1$, \ we have \ $Y_{\lfloor nt_1\rfloor} \distr \cX_{t_1}$ \ as \ $n \to \infty$, \ since for each \ $n\in\NN$, \
 the distribution of \ $Y_{\lfloor nt_1\rfloor}$ \ is \ $\mu$ \ which coincides with the distribution of \ $\cX_{t_1}$.
\ Suppose that \eqref{stationary_convergence} holds for some \ $m \in \NN$, \ $t_1, \ldots, t_m \in \RR_{++}$ \ with \ $t_1 < \ldots < t_m$,
 \ and let \ $t_{m+1}\in\RR_{++}$ \ such that \ $t_m < t_{m+1}$.
\ The strongly stationary Markov chain \ $(Y_k)_{k\in\ZZ_+}$ \ is strongly mixing, i.e.,
 \[
   \alpha_\ell := \sup_{i\in\ZZ_+} \sup_{A\in\cF_{0,i}^Y,\; B\in\cF_{i+\ell, \infty}^Y} |\PP(A \cap B) - \PP(A) \PP(B)| \to 0 \qquad \text{as \ $\ell \to \infty$,}
 \]
 where \ $\cF_{0,i}^Y := \sigma(Y_0,\ldots, Y_i)$ \ and \ $\cF_{i,\infty}^Y := \sigma(Y_i, Y_{i+1}, \ldots)$ \ for \ $i \in \ZZ_+$,
 \ see Barczy et al. \cite[Lemma F.1]{BarNedPap} or  Basrak et al.\ \cite[Remark 3.1]{BasKulPal}.
The strong stationarity of \ $(Y_k)_{k\in\ZZ_+}$ \ implies that
 \[
   |\EE(U V) - \EE(U) \EE(V)| \leq 4 C_1 C_2 \alpha_m
 \]
 for any \ $\cF_{0,j}^Y$-measurable (real-valued) random variable \ $U$ \ and any \ $\cF_{j+m,\infty}^Y$-measurable random variable \ $V$ \
 with \ $j, m \in \NN$, \ $|U| \leq C_1$ \ and \ $|V| \leq C_2$ \ (see, e.g., Lemma 1.2.1 in Lin and Lu \cite{LinLu}).
Hence, for any \ $\cF_{0,j}^Y$-measurable complex-valued random variable \ $U$ \ and any \ $\cF_{j+m,\infty}^Y$-measurable complex-valued
 random variable \ $V$ \ with \ $j, m \in \NN$, \ $|U| \leq C_1$ \ and \ $|V| \leq C_2$, \ we get
 \begin{align*}
  &|\EE(U V) - \EE(U) \EE(V)| \\
  &\leq |\EE(\Re(U) \Re(V)) - \EE(\Re(U)) \EE(\Re(V))|
       + |\EE(\Im(U) \Im(V)) - \EE(\Im(U)) \EE(\Im(V))| \\
  &\phantom{\leq}
       + |\EE(\Re(U) \Im(V)) - \EE(\Re(U)) \EE(\Im(V))|
       + |\EE(\Im(U) \Re(V)) - \EE(\Im(U)) \EE(\Re(V))|
  \leq 16 C_1 C_2 \alpha_m ,
 \end{align*}
 where for a complex number \ $z\in\CCC$, \ $\Re(z)$ \ and \ $\Im(z)$ \ denote the real and imaginary part of \ $z$, \ respectively.
Consequently, for each \ $u_1, \ldots, u_m, u_{m+1} \in \RR$, \ we obtain
 \begin{align}\label{help6}
  \begin{split}
  & |\EE(\ee^{\ii(u_1Y_{\lfloor nt_1\rfloor}+\cdots+u_{m+1}Y_{\lfloor nt_{m+1}\rfloor})})
    - \EE(\ee^{\ii(u_1Y_{\lfloor nt_1\rfloor}+\cdots+u_mY_{\lfloor nt_m\rfloor})})
      \EE(\ee^{\ii u_{m+1}Y_{\lfloor nt_{m+1}\rfloor}})|\\
  &\leq 16 \alpha_{\lfloor nt_{m+1}\rfloor-\lfloor nt_m\rfloor} .
 \end{split}
 \end{align}
We have \ $\alpha_{\lfloor nt_{m+1}\rfloor-\lfloor nt_m\rfloor} \to 0$ \ as \ $n \to \infty$, \ since
 \[
   \lfloor nt_{m+1}\rfloor-\lfloor nt_m\rfloor \geq nt_{m+1} - 1 - n t_m = (t_{m+1} - t_m) n - 1 \to \infty \qquad \text{as \ $n \to \infty$.}
 \]
Since \eqref{stationary_convergence} holds for \ $m \in \NN$, \ $t_1,\ldots,t_m\in\RR_{++}$ \ with \ $t_1<\ldots<t_m$, \ by the continuity theorem, we have
 \[
   \EE(\ee^{\ii(u_1Y_{\lfloor nt_1\rfloor}+\cdots+u_mY_{\lfloor nt_m\rfloor})})
   \to \EE(\ee^{\ii(u_1\cX_{t_1}+\cdots+u_m\cX_{t_m})}) \qquad \text{as \ $n \to \infty$.}
 \]
Moreover, by using \eqref{stationary_convergence} with \ $m=1$, \ and the continuity theorem,
 we also have \ $\EE(\ee^{\ii u_{m+1}Y_{\lfloor nt_{m+1}\rfloor}}) \to \EE(\ee^{\ii u_{m+1}\cX_{t_{m+1}}})$ \ as \ $n \to \infty$, \ hence, by \eqref{help6}, we conclude for all \ $u_1,\ldots, u_m, u_{m+1}\in\RR$,
 \begin{align*}
  \EE(\ee^{\ii(u_1Y_{\lfloor nt_1\rfloor}+\cdots+u_{m+1}Y_{\lfloor nt_{m+1}\rfloor})})
  &\to \EE(\ee^{\ii(u_1\cX_{t_1}+\cdots+u_m\cX_{t_m})})
       \EE(\ee^{\ii u_{m+1}\cX_{t_{m+1}}}) \\
  &= \EE(\ee^{\ii(u_1\cX_{t_1}+\cdots+u_{m+1}\cX_{t_{m+1}})})
   \qquad \text{as \ $n \to \infty$,}
 \end{align*}
 by the independence of \ $(\cX_{t_1}, \ldots, \cX_{t_m})$ \ and \ $\cX_{t_{m+1}}$.
\ Again by the continuity theorem, we get
 \[
   (Y_{\lfloor nt_1\rfloor}, \ldots, Y_{\lfloor nt_{m+1}\rfloor})
   \distr (\cX_{t_1}, \ldots, \cX_{t_{m+1}}) \qquad \text{as \ $n \to \infty$,}
 \]
 which is \eqref{stationary_convergence} with \ $m$ \ replaced by \ $m + 1$, \ as desired.

Next, using a coupling argument, we show that for each \ $m\in\NN$, \ $t_1, \ldots, t_m \in \RR_{++}$ \ with \ $t_1 < \ldots < t_m$, \ we have
 \begin{equation}\label{stoch_convergence_difference}
  (X_{\lfloor nt_1\rfloor}, \ldots, X_{\lfloor nt_m\rfloor})
    - (Y_{\lfloor nt_1\rfloor}, \ldots, Y_{\lfloor nt_m\rfloor})
  \stoch (0, \ldots, 0) \qquad \text{as \ $n \to \infty$,}
 \end{equation}
 which, together with \eqref{stationary_convergence} and Slutsky's lemma, yields that
 for each \ $m\in\NN$, \ $t_1, \ldots, t_m \in \RR_{++}$ \ with \ $t_1 < \ldots < t_m$, \ we have
 \begin{equation*}
  (X_{\lfloor nt_1\rfloor}, \ldots, X_{\lfloor nt_m\rfloor})
  \distr (\cX_{t_1}, \ldots, \cX_{t_m}) \qquad \text{as \ $n \to \infty$,}
 \end{equation*}
 as desired.
Observe that
 \[
   X_1 = \sum_{j=1}^{X_0} \xi_{1,j} + \vare_1 , \qquad Y_1 = \sum_{j=1}^{Y_0} \xi_{1,j} + \vare_1
 \]
 implies
 \[
   |X_1 - Y_1| = \sum_{j=X_0\land Y_0}^{X_0\lor Y_0} \xi_{1,j} ,
 \]
 where \ $x\land y:=\min(x,y)$ \ and \ $x\lor y:=\max(x,y)$ \ for \ $x,y\in\RR$.
\ Thus
 \[
   \EE(|X_1 - Y_1| \mid X_0, Y_0)
   = ((X_0 \lor Y_0) - (X_0 \land Y_0)) a
   = a |X_0 - Y_0|  ,
 \]
 and hence
 \[
   \EE(|X_1 - Y_1|) = a \EE(|X_0 - Y_0|) \leq a (\EE(X_0) + \EE(Y_0)) .
 \]
In a similar way, by recursion, we obtain
 \[
   \EE(|X_n - Y_n|) = a \EE(|X_{n-1} - Y_{n-1}|) \leq a^n (\EE(X_0) + \EE(Y_0)) , \qquad n \in \NN.
 \]
Hence
 \[
   \EE(|X_n - Y_n|) \to 0 \qquad \text{as $n\to\infty$,}
 \]
 yielding
 \[
   |X_n - Y_n| \stoch 0 \qquad \text{as $n\to\infty$.}
 \]
Since \ $X_n$ \ and \ $Y_n$ \ are (nonnegative) integer-valued random variables, we conclude
 \[
   \PP(X_n = Y_n) \to 1 \qquad \text{as \ $n\to\infty$.}
 \]
If \ $X_N = Y_N$ \ is satisfied for some \ $N \in \NN$, \ then, by the definition of \ $(Y_k)_{k\in\ZZ_+}$, \ we have
 $X_n = Y_n$ \ is satisfied for all \ $n \geq N$, \ thus
 \[
   \PP(\text{$X_n = Y_n$ for all $n \geq N$}) \to 1 \qquad \text{as \ $N \to\infty$.}
 \]
For each \ $N \in \NN$, \ let \ $\Omega_N := \{\text{$X_n = Y_n$ for all $n \geq N$}\}$.
For each \ $\delta \in \RR_{++}$ \ and \ $n, N \in \NN$ \ with \ $\lfloor nt_1\rfloor > N$, \ we have
 \begin{align*}
  &\PP\bigl( \Vert (X_{\lfloor nt_1\rfloor}, \ldots, X_{\lfloor nt_m\rfloor})
           - (Y_{\lfloor nt_1\rfloor}, \ldots, Y_{\lfloor nt_m\rfloor}) \Vert > \delta \bigr) \\
  &= \PP\bigl(\bigl\{ \Vert (X_{\lfloor nt_1\rfloor}, \ldots, X_{\lfloor nt_m\rfloor})
    - (Y_{\lfloor nt_1\rfloor}, \ldots, Y_{\lfloor nt_m\rfloor}) \Vert
                   > \delta \bigr\} \cap \Omega_N\bigr) \\
  &\quad
     + \PP\bigl(\bigl\{ \Vert (X_{\lfloor nt_1\rfloor}, \ldots, X_{\lfloor nt_m\rfloor})
                   - (Y_{\lfloor nt_1\rfloor}, \ldots, Y_{\lfloor nt_m\rfloor}) \Vert > \delta \bigr\} \cap \Omega_N^\mathrm{c}\bigr) \\
  &\leq \PP( \Vert (0, \ldots, 0) \Vert > \delta ) + \PP(\Omega_N^\mathrm{c})
   = \PP(\Omega_N^\mathrm{c}),
 \end{align*}
 where \ $\Omega_N^\mathrm{c}$ \ denotes the complement of \ $\Omega_N$.
\ Letting \ $N\to\infty$, \ we obtain \eqref{stoch_convergence_difference}, as desired.
\proofend

The following result is about the asymptotic behaviour of a single-type critical Galton-Watson process with immigration
 due to Wei and Winnicki \cite[Theorem 2.1]{WW}.

\begin{Thm}\label{thm:critical}
Let \ $(X_k)_{k\in\ZZ_+}$ \ be a single-type Galton-Watson process with immigration such that
 \ $\EE(\xi^2) < \infty$, \ $\EE(\vare^2) < \infty$, \ $\EE(\xi) = 1$ \ (critical case) and \ $\EE(X_0^2) < \infty$.
\ Then
 \begin{equation*}
  (n^{-1} X_\nt)_{t\in\RR_+} \distr (\cX_t)_{t\in\RR_+} \qquad
  \text{as \ $n \to \infty$,}
 \end{equation*}
 where the limit process \ $(\cX_t)_{t\in\RR_+}$ \ is the pathwise unique strong solution of the SDE
 \begin{equation*}
  \dd \cX_t
  = \EE(\vare) \, \dd t
    + \sqrt{\var(\xi) \, \cX_t^+} \, \dd \cW_t , \qquad t\in\RR_+,
  \qquad \cX_0 = 0 ,
 \end{equation*}
 where \ $(\cW_t)_{t\in\RR_+}$ \ is a standard Wiener process.
\end{Thm}

\section{A version of the continuous mapping theorem}
\label{CMT}

A function \ $f : \RR_+ \to \RR^d$ \ is called \emph{c\`adl\`ag} if it is right
 continuous with left limits.
\ Let \ $\DD(\RR_+, \RR^d)$ \ and \ $\CC(\RR_+, \RR^d)$ \ denote the space of
 all $\RR^d$-valued c\`adl\`ag and continuous functions on \ $\RR_+$,
 \ respectively.
Let \ $\cB(\DD(\RR_+, \RR^d))$ \ denote the Borel $\sigma$-algebra on
 \ $\DD(\RR_+, \RR^d)$ \ for the metric defined in Jacod and Shiryaev
 \cite[Chapter VI, (1.26)]{JacShi} (with this metric \ $\DD(\RR_+, \RR^d)$ \ is a
 complete and separable metric space and the topology induced by this metric is
 the so-called Skorokhod topology).
Note that \ $\CC(\RR_+,\RR^d)\in\cB(\DD(\RR,\RR^d))$, \ see, e.g., Ethier and Kurtz \cite[Problem 3.11.25]{EthKur}.
For $\RR^d$-valued stochastic processes \ $(\bcY_t)_{t \in \RR_+}$ \ and
 \ $(\bcY^{(n)}_t)_{t \in \RR_+}$, \ $n \in \NN$, \ with c\`adl\`ag paths we write
 \ $\bcY^{(n)} \distr \bcY$ \ if the distribution of \ $\bcY^{(n)}$ \ on the
 space \ $(\DD(\RR_+, \RR), \cB(\DD(\RR_+, \RR^d)))$ \ converges weakly to the
 distribution of \ $\bcY$ \ on the space
 \ $(\DD(\RR_+, \RR), \cB(\DD(\RR_+, \RR^d)))$ \ as \ $n \to \infty$.
\ If \ $\xi$ \ and \ $\xi_n$, \ $n \in \NN$, \ are random elements with values in a metric space \ $(E, d)$,
 \ then we denote by \ $\xi_n \distr \xi$ \ the weak convergence of the
 distribution of \ $\xi_n$ \ on the space \ $(E, \cB(E))$ \ towards the
 distribution of \ $\xi$ \ on the space \ $(E, \cB(E))$ \ as \ $n \to \infty$,
 \ where \ $\cB(E)$ \ denotes the Borel \ $\sigma$-algebra on \ $E$ \ induced by
 the given metric \ $d$.

The following version of the continuous mapping theorem can be found for
 example in Theorem 3.27 of Kallenberg \cite{Kal}.

\begin{Lem}\label{lemma:kallenberg}
Let \ $(S, d_S)$ \ and \ $(T, d_T)$ \ be metric spaces and
 \ $(\xi_n)_{n \in \NN}$, \ $\xi$ \ be random elements with values in \ $S$
 \ such that \ $\xi_n \distr \xi$ \ as \ $n \to \infty$.
\ Let \ $f : S \to T$ \ and \ $f_n : S \to T$, \ $n \in \NN$, \ be measurable
 mappings and \ $C \in \cB(S)$ \ such that \ $\PP(\xi \in C) = 1$ \ and
 \ $\lim_{n \to \infty} d_T(f_n(s_n), f(s)) = 0$ \ if
 \ $\lim_{n \to \infty} d_S(s_n,s) = 0$ \ and \ $s \in C$, \ $s_n\in S$, \ $n\in\NN$.
\ Then \ $f_n(\xi_n) \distr f(\xi)$ \ as \ $n \to \infty$.
\end{Lem}

For functions \ $f$ \ and \ $f_n$, \ $n \in \NN$, \ in \ $\DD(\RR_+, \RR^d)$,
 \ we write \ $f_n \lu f$ \ if \ $(f_n)_{n \in \NN}$ \ converges to \ $f$
 \ locally uniformly, i.e., if \ $\sup_{t \in [0,T]} \|f_n(t) - f(t)\| \to 0$
 \ as \ $n \to \infty$ \ for all \ $T \in \RR_{++}$.
\ For measurable mappings \ $\Phi : \DD(\RR_+, \RR^d) \to \DD(\RR_+, \RR^q)$
 \ and \ $\Phi_n : \DD(\RR_+, \RR^d) \to \DD(\RR_+, \RR^q)$, \ $n \in \NN$,
 \ we will denote by \ $C_{\Phi, (\Phi_n)_{n \in \NN}}$ \ the set of all functions
 \ $f \in \CC(\RR_+, \RR^d)$ \ for which \ $\Phi_n(f_n) \lu \Phi(f)$ \ whenever
 \ $f_n \lu f$ \ with \ $f_n \in \DD(\RR_+, \RR^d)$, \ $n \in \NN$.

One can formulate the following consequence of Lemma \ref{lemma:kallenberg}.

\begin{Lem}\label{Conv2Funct}
Let \ $d,\,q\in\NN$.
\ Let \ $(\bcU_t)_{t \in \RR_+}$ \ and \ $(\bcU^{(n)}_t)_{t \in \RR_+}$, \ $n \in \NN$,
 \ be \ $\RR^d$-valued stochastic processes with c\`adl\`ag paths such that
 \ $\bcU^{(n)} \distr \bcU$ \ as \ $n\to\infty$.
\ Let \ $\Phi : \DD(\RR_+, \RR^d) \to \DD(\RR_+, \RR^q)$ \ and
 \ $\Phi_n : \DD(\RR_+, \RR^d) \to \DD(\RR_+, \RR^q)$, \ $n \in \NN$, \ be
 measurable mappings such that there exists
 \ $C \subset C_{\Phi,(\Phi_n)_{n\in\NN}}$ \ with \ $C \in \cB(\DD(\RR_+, \RR^d))$ \ and \ $\PP(\bcU \in C) = 1$.
\ Then \ $\Phi_n(\bcU^{(n)}) \distr \Phi(\bcU)$ \ as \ $n\to\infty$.
\end{Lem}

\section{Convergence of random step processes}
\label{section_conv_step_processes}

We recall a result about convergence of random step processes towards a diffusion process, see Isp\'any and Pap \cite{IspPap}.

\begin{Thm}\label{Conv2DiffThm}
Let \ $\bbeta : \RR_+ \times \RR^d \to \RR^d$ \ and \ $\bgamma : \RR_+ \times \RR^d \to \RR^{d \times r}$ \ be continuous
 functions.
Assume that uniqueness in the sense of probability law holds for the SDE
 \begin{equation}\label{SDE}
  \dd \, \bcU_t
  = \bbeta (t, \bcU_t) \, \dd t + \bgamma (t, \bcU_t) \, \dd \bcW_t ,
  \qquad t \in \RR_+,
 \end{equation}
 with initial value \ $\bcU_0 = \bu_0$ \ for all \ $\bu_0 \in \RR^d$, \ where
 \ $(\bcW_t)_{t\in\RR_+}$ \ is an $r$-dimensional standard Wiener process.
Let \ $(\bcU_t)_{t\in\RR_+}$ \ be a solution of \eqref{SDE} with initial value
 \ $\bcU_0 = \bzero \in \RR^d$.

For each \ $n \in \NN$, \ let \ $(\bU^{(n)}_k)_{k\in\ZZ_+}$ \ be a sequence of $d$-dimensional random vectors adapted to a filtration
 \ $(\cF^{(n)}_k)_{k\in\ZZ_+}$ \ (i.e., \ $\bU^{(n)}_k$ \ is \ $\cF^{(n)}_k$-measurable) such that \ $\EE(\Vert \bU^{(n)}_k\Vert^2)<\infty$ \ for each \ $n,k \in \NN$.
\ Let
 \[
   \bcU^{(n)}_t := \sum_{k=0}^{\nt}  \bU^{(n)}_k \, ,
   \qquad t \in \RR_+, \quad n \in \NN .
 \]
Suppose that \ $\bcU^{(n)}_0 = \bU^{(n)}_0 \distr \bzero$ \ as \ $n\to\infty$ \ and that for each \ $T \in \RR_{++}$,
 \begin{enumerate}
  \item[\textup{(i)}]
   $\sup\limits_{t\in[0,T]}
     \biggl\|\sum\limits_{k=1}^{\nt}
              \EE\bigl(\bU^{(n)}_k \mid \cF^{(n)}_{k-1}\bigr)
             - \int_0^t \bbeta(s,\bcU^{(n)}_s) \dd s\biggr\|
   \stoch 0$ \ as \ $n\to\infty$,  \\
  \item[\textup{(ii)}]
   $\sup\limits_{t\in[0,T]}
     \biggl\|\sum\limits_{k=1}^{\nt}
              \var\bigl(\bU^{(n)}_k \mid \cF^{(n)}_{k-1}\bigr)
             - \int_0^t
                \bgamma(s,\bcU^{(n)}_s) \bgamma(s,\bcU^{(n)}_s)^\top
                \dd s\biggr\|
   \stoch 0$ \ as \ $n\to\infty$, \\
  \item[\textup{(iii)}]
   $\sum\limits_{k=1}^{\lfloor nT \rfloor}
     \EE\bigl(\|\bU^{(n)}_k\|^2 \bbone_{\{\|\bU^{(n)}_k\| > \theta\}}
             \bmid \cF^{(n)}_{k-1}\bigr)
   \stoch 0$ \ as \ $n\to\infty$ \ for all \ $\theta \in \RR_{++}$.
 \end{enumerate}
Then \ $\bcU^{(n)} \distr \bcU$ \ as \ $n \to \infty$.
\end{Thm}

Note that in (ii) of Theorem \ref{Conv2DiffThm}, \ $\|\cdot\|$ \ denotes
 an operator norm, while in (i) it denotes a vector norm.

\section*{Acknowledgements}
This paper was finished after the sudden death of the third author Gyula Pap in October 2019.
We would like to thank the referees for their comments that helped us improve the paper.

\end{document}